\pdfoutput=1
\documentclass[11pt,reqno]{amsart}
\usepackage[T1]{fontenc}
\usepackage[utf8]{inputenc}
\usepackage{lmodern}
\usepackage{amsmath,amssymb,amsthm,mathtools}
\usepackage[a4paper,left=29mm,right=29mm,top=30mm,bottom=30mm]{geometry}
\usepackage{microtype}
\usepackage{xcolor}
\usepackage{booktabs,array,enumitem}
\usepackage[hyphens]{url}
\definecolor{citeblue}{RGB}{0,62,168}
\definecolor{linkred}{RGB}{168,24,36}
\definecolor{urlgreen}{RGB}{0,104,74}
\usepackage[colorlinks=true,citecolor=citeblue,linkcolor=linkred,urlcolor=urlgreen,
  bookmarksnumbered=true,pdflang={en},
  pdftitle={Simultaneous Fermat pseudoprimes with two prime factors: rank-dependent logarithmic savings and sifted base-3 proportions},
  pdfauthor={Santiago Pineda Montoya, Johan Heriberto Rúa Muñoz, Jaiberth Porras Barrera, Edisson Gallego González}]{hyperref}
\newcommand{\doi}[1]{\href{https://doi.org/#1}{doi:\nolinkurl{#1}}}

\newtheorem{theorem}{Theorem}[section]
\newtheorem{lemma}[theorem]{Lemma}
\newtheorem{prop}[theorem]{Proposition}
\newtheorem{cor}[theorem]{Corollary}
\newtheorem{introtheorem}{Theorem}
\theoremstyle{remark}
\newtheorem{remark}[theorem]{Remark}

\DeclareMathOperator{\ord}{ord}
\DeclareMathOperator{\lcm}{lcm}
\newcommand{\Z}{\mathbb Z}
\newcommand{\Q}{\mathbb Q}
\newcommand{\F}{\mathbb F}
\newcommand{\Pp}{\mathbb P}
\newcommand{\A}{\mathcal A}
\newcommand{\Qq}{\mathcal Q}
\newcommand{\B}{\mathcal B}
\newcommand{\dd}{\delta_{\!F}}
\newcommand{\wt}{\mathcal W}
\newcommand{\Ocal}{\mathcal O}
\newcommand{\Ctw}{C_{\rm tw}}
\newcommand{\Nrm}{\operatorname{N}}
\newcommand{\Gal}{\operatorname{Gal}}
\newcommand{\Gw}{\mathcal G}
\begin{document}

\title[Simultaneous Fermat pseudoprimes with two prime factors]{Simultaneous Fermat pseudoprimes with two prime factors: rank-dependent logarithmic savings and sifted base-3 proportions}

\author{Santiago Pineda Montoya}
\author{Johan Heriberto Rúa Muñoz}
\author{Jaiberth Porras Barrera}
\author{Edisson Gallego González}
\renewcommand{\shortauthors}{S. Pineda Montoya, J. H. Rúa Muñoz, J. Porras Barrera and E. Gallego González}

\date{November 29, 2026}

\subjclass[2020]{Primary 11A51; Secondary 11N25, 11N36}
\keywords{Fermat pseudoprimes, semiprimes, multiplicative order, residual index, power reciprocity, sieve methods}

\begin{abstract}
We count products of two distinct primes up to $X$ that pass the Fermat test in prescribed bases. For base 2 and for bases 2 and 3 we obtain rank-dependent logarithmic savings at the respective powers $X^{3/4}$ and $X^{2/3}$. Under the generalized Riemann hypothesis for the associated Kummer fields, the two-base bound improves to $O_\varepsilon(X^{13/20+\varepsilon})$ for every $\varepsilon>0$. A power saving below $3/4$ for the base-2 count, or for overpseudoprimes to base 2, is equivalent to the following estimate for some $\delta,\eta>0$: the number of primes $p\le y$ for which 2 has multiplicative order at most $p^{1/2+\delta}$ is $O(y^{1-\eta})$. On five fixed ratios of $q-1$ to $p-1$, we give exact base-3 criteria for prime pairs passing base 2. On every fixed reduced ratio $b/a>1$ with $a\le2$, the associated local conditions are asymptotically independent when $q$ is prime and its companion has no small prime factor. For the five specified ratios we obtain the same proportions for prime pairs under explicit prime-value conjectures. In the factored catalog up to $10^{16}$, six ratios account for $315\,497$ of the $394\,131$ simultaneous semiprimes.
\end{abstract}

\maketitle

\section{Introduction}\label{sec:intro}

For a composite integer $n$ and an integer base $t$ coprime to $n$, the Fermat test accepts $n$ if $t^{n-1}\equiv1\pmod n$. The products $341=11\cdot31$ and $31\,621=103\cdot307$ both pass the test in base 2. Base 3 rejects $341$ and accepts $31\,621$. In both factorizations the smaller prime $p$ and the larger prime $q$ satisfy $q-1=3(p-1)$. We address two questions.
\begin{enumerate}[label=\textbf{Question~\arabic*.},ref=\arabic*,leftmargin=*,widest=2]
\item\label{q:local} For prime pairs with a fixed ratio $(q-1)/(p-1)$ that pass base 2, what additional condition determines whether their product passes base 3?
\item\label{q:count} Given a product bound $X$, how many pairs $p<q$ with $pq\le X$ pass both tests as the ratio varies freely?
\end{enumerate}

\subsection{Background}

For a unit $t$ modulo $m>1$ (a residue class coprime to $m$), let $\ord_m(t)$ denote its multiplicative order, the least $h\ge1$ with $t^h\equiv1\pmod m$. At a prime $p$, a unit of order $p-1$ is a primitive root and $(p-1)/\ord_p(t)$ is its residual index. For an odd composite $n$, write $n-1=2^s e$ with $e$ odd. The strong test in base $t$ accepts if $t^e\equiv1\pmod n$ or $t^{2^j e}\equiv-1\pmod n$ for some $0\le j<s$. The notation $f\ll_A g$ or $f=O_A(g)$ means $|f|\le C_A g$ for a constant depending only on $A$. $f\asymp_A g$ means bounds in both directions, $f\sim g$ means $f/g\to1$, and $f=o(g)$ means $f/g\to0$. Throughout, $\log$ denotes the natural logarithm.

Lehmer \cite[p.~348]{Lehmer} constructed Fermat pseudoprimes with two prime factors, and Erdős \cite{E49} proved that infinitely many squarefree pseudoprimes have exactly $k$ prime factors for each fixed $k\ge2$. Primitive-divisor theorems of Zsigmondy and Schinzel \cite{BZ,Schinzel}, in the forms recorded by Somer \cite{Somer}, provide primes at which a base has prescribed multiplicative order. Together with Bateman's asymptotic for the number of integers with totient at most a given cutoff \cite{Bateman}, these theorems yield our lower bounds for two disjoint families in base 2, those that pass the strong test and those that fail it.

A Fermat test at a composite modulus can be expressed in terms of multiplicative orders at its prime factors. Pomerance, Selfridge and Wagstaff \cite{PSW} studied the classical pseudoprime tests and tabulated examples through $25\cdot10^9$. Jaeschke \cite[Proposition 1]{Jaeschke} characterized simultaneous strong pseudoprimality by the powers of 2 dividing the orders at the prime factors. Pomerance \cite{Pom81} also obtained an upper bound for the number of pseudoprimes over all composite moduli. Wagstaff \cite{Wagstaff} related pseudoprimes to the residual indices in Artin's primitive-root problem, which Hooley \cite{Hooley} treated under the generalized Riemann hypothesis (GRH). Bounds for primes at which a base has small order, including those of Erdős and Murty \cite{EM}, are relevant to our counting argument.

For several bases, one can track the subgroup they generate and the order of its image at each prime. Matthews, Gupta and Murty, and Pappalardi \cite{Matt,GM,Pappa} studied how this order varies with the prime. Sorenson and Webster \cite[Lemma 3.4 and \S3.4.1]{SW} used Pappalardi's reciprocal-order bound and the restriction to an arithmetic progression in their algorithm for strong pseudoprimes. Section~\ref{ssec:main} explains what these two ingredients give for ordinary Fermat pseudoprimes and where Ford's shifted-prime divisor bound \cite{Ford} enters.

For base 2, Galway \cite{Galway} proposed a global asymptotic of the form $C\sqrt X/(\log X)^2$, where $X$ bounds the product and $C$ is a constant. The families in this model are indexed by the reduced ratio $(q-1)/(p-1)=b/a$, which we call a \emph{ray}. Granville and Pomerance \cite[\S8 and note added in proof, p.~908]{GP} explained the power-residue conditions on a ray, used Chebotarev for the individual base-2 densities, and corrected the conjectural constant. Feitsma and Galway's factored catalog \cite{FG}, extending to $2^{64}$, provides finite data with which to compare counts on individual rays and the total count. Passing from finitely many rays to a global asymptotic requires uniform control as the ratio varies.

The closest quantitative comparison comes from the study of exceptional moduli in RSA factorization. In the window $z<p<q\le2z$, Dieulefait and Jiménez Urroz \cite[Theorem 2.3]{DJU10} obtained an upper bound of $O(z^2/(\log z)^3)$ for the exceptional base-2 moduli. Their bound applies because $pq\mid2^g-1$ for $g=\gcd(p-1,q-1)$, which forces $g>\log z$ for large $z$. Their later work \cite{DJU} further developed the relation between malleability and factorization. Luca, Ralaivaosaona and Jiménez Urroz \cite[Theorems 2.1 and 3.1]{LRU} proved the bounds
\[
 O\!\left(\frac{z^{8/5}}{(\log z)^2}\right)
 \qquad\text{and}\qquad
 O\!\left(\frac{X^{4/5}}{(\log X)^{2/5}}\right)
\]
for the window and the product cutoff, respectively. Their Propositions 3.1 and 3.2 give a family with $p<q<2p$ and a global lower bound of logarithmic size. Li and Pomerance \cite[Theorem 1]{LiP} study a complementary uniformity problem. For each integer base $2\le a\le X$, uniformly in $a$, they bound the number of composite integers $n\le X$ satisfying $a^n\equiv a\pmod n$ with the upper bound
\[
 X\exp\!\left(-\frac{\log X\,\log\log\log X}{2\log\log X}\right).
\]

A different way to bound the base-2 count is to estimate the number of primes at which 2 has a given exact order $t$. Conjecture 2 of Kovač and Luca \cite{KL} would bound the number of such primes by $O(\log(2t))$. We combine this hypothesis with a sieve for two linear forms and Pollack's mean-value estimate for multiplicative functions on sifted sets \cite{Pollack}.

On a fixed ray, power-residue characters give exact tests for individual pairs. Zhang \cite[(4.1) and \S5]{Zhang01} gave the ordinary Fermat equivalence for the ratio $4$ and used the cubic character for the ratio $3$ in his study of strong pseudoprimes. Zhang and Tang \cite[\S3, (3.1)--(3.2)]{ZT03} gave the ordinary Fermat equivalences for the ratios $3/2$ and $6$ before applying the strong test to their candidates. For the octic condition on the ratio $8$ we use Western's criterion as recorded by Hudson \cite{Hudson}. Section~\ref{sec:rayos} specializes these character conditions to base 3 and expresses them in the coordinates supplied by base 2.

We count the resulting local conditions with a sifted companion using the Chebotarev density theorem and the distribution theorem of Murty and Murty \cite{MM}, in the form given by Thorner \cite{Thorner}. Requiring both factors to be prime leads to hypotheses on simultaneous prime values. For linear forms these are the conjectures of Dickson and Hardy--Littlewood \cite{Dickson,HL}. For polynomial pairs we use Bateman--Horn \cite{BH} and its multivariable formulation due to Destagnol and Sofos \cite[Conjecture A.3 and Remark A.4]{DS}.

Monier \cite{Monier} counted the bases accepted by probabilistic primality tests. Zhang and Tang \cite[(1.2)]{ZT03} used the fraction accepted by the strong test to select useful families for their search, and Ishmukhametov, Mubarakov and Rubtsova \cite[Theorem 4]{IMR} gave its explicit numerator for semiprimes. Erdős and Pomerance \cite{EP} studied how such counts vary over composite moduli. We average fixed moments of the fractions over products whose two prime factors lie in one dyadic window, using Maynard's theorem on prime values of linear forms \cite{Maynard} for the lower bound. Equal orders at all prime factors link our results to the overpseudoprimes of Shevelev and collaborators \cite{Shev12}, whose number up to $X$ Shevelev \cite[Theorem~3]{Shev08} bounded by $(1+o(1))X^{3/4}$.

\subsection{Main results}\label{ssec:main}

Theorem~\ref{intro:rank} below answers Question~\ref{q:count} of the introduction for any fixed group of bases, Theorem~\ref{intro:rays} answers Question~\ref{q:local}, and Theorem~\ref{intro:over} treats a related question about products whose factors give 2 the same order, for example $2047=23\cdot89$.

A fixed finitely generated multiplicative group $\Gamma\subset\Q^\times$ has \emph{rank} $r$ if its quotient by its finite torsion subgroup is isomorphic to $\Z^r$. A prime is of \emph{good reduction} for $\Gamma$ if it divides none of the numerators or denominators of a fixed generating set in lowest terms. An integer is of good reduction if each of its prime factors is. At a good-reduction modulus $n$, a rational $u/v$ in lowest terms represents the unit $uv^{-1}\bmod n$. Let $\Qq_\Gamma(X)$ count the products of distinct odd good-reduction primes $p<q$ with $pq\le X$ for which every element of $\Gamma$ passes Fermat. Let $\A_\Gamma(z)$ count the same products with $z<p<q\le2z$. A subscript $2$ means $\Gamma=\langle2\rangle$, and $23$ means $\Gamma=\langle2,3\rangle$.

For sufficiently large $y$, set
\[
 L_y=\log y,\quad H_y=\log\log y,
 \quad \dd=1-\frac{1+\log\log2}{\log2},\quad \beta=\frac1{r+1},
 \quad \wt_r(y)=\frac{1}{L_y^{\dd\beta}H_y^{3\beta/2}}.
\]
The constant $\dd=0.0860713\ldots$ is the exponent in Ford's theorem on integers with a divisor in a given interval \cite{Ford}.

\begin{introtheorem}[Rank-dependent bounds]\label{intro:rank}
For every fixed $\Gamma$ of rank $r\ge1$, unconditionally,
\[
 \A_\Gamma(z)\ll_\Gamma \frac{z^{1+\beta}\wt_r(z)}{L_z^2},\qquad
 \Qq_\Gamma(X)\ll_\Gamma \frac{X^{(1+\beta)/2}\wt_r(X)}{L_X^2}.
\]
In particular,
\[
 \Qq_2(X)\ll \frac{X^{3/4}}{L_X^{2+\dd/2}H_X^{3/4}},
 \qquad
 \Qq_{23}(X)\ll \frac{X^{2/3}}{L_X^{2+\dd/3}\sqrt{H_X}}.
\]
\end{introtheorem}

Theorem~\ref{intro:rank} consists of Theorem~\ref{th:main} and Corollary~\ref{cor:23}, which are proved in
Section~\ref{ssec:rect} from the estimates of Section~\ref{ssec:small}. In base 2, the power of $z$ drops from $8/5$ to $3/2$ and the power of $X$ from $4/5$ to $3/4$ relative to \cite[Theorems 2.1 and 3.1]{LRU}. We have not found these exponents stated in the literature, but their proofs are short. Section~\ref{ssec:elem} deduces $\A_2(z)\ll z^{3/2}/L_z^2$ from the Brun--Titchmarsh inequality, and Section~\ref{ssec:compare} deduces the global powers $3/4$ in base 2 and $2/3$ in bases 2 and 3 from Pappalardi's reciprocal-order bound \cite[Corollary~2.4]{Pappa} combined with the restriction to an arithmetic progression that Sorenson and Webster \cite[Lemma~3.4]{SW} use in an algorithm for strong pseudoprimes. The displayed product-cutoff estimates for ordinary semiprime pseudoprimes are deductions from these inputs. The main technical content of Theorem~\ref{intro:rank} is the factor $\wt_r$. It comes from Ford's theorem on shifted primes with a divisor in a given interval, applied near the order scale $q^{r/(r+1)+o(1)}$ for a prime $q$. Brun--Titchmarsh attaches the weight $h_q/\varphi(h_q)$ to each such prime, where $h_q=|\Gamma\bmod q|$ is the order of the image of $\Gamma$ in $\F_q^\times$, and bounding this weight pointwise outside a sparse exceptional set would cost a factor $\log\log\log X$. Lemma~\ref{le:ford} avoids this loss: it bounds the moments of $(q-1)/\varphi(q-1)$ over the primes in Ford's estimate, and Hölder's inequality combines them with the small-order estimate. In base 2 the factor $\wt_1$ improves the elementary bound $z^{3/2}/L_z^2$ by $L_z^{\dd/2}H_z^{3/4}$, which tends to infinity. For bases 2 and 3 together the exponents are $4/3$ in the window and $2/3$ globally. Theorem~\ref{th:momento} characterizes power bounds on the global count in terms of the second moment of the multiplicities of the reduction orders, and Corollary~\ref{co:erdos} in terms of the number of primes with small reduction order. Under GRH for the Dedekind zeta functions in Corollary~\ref{co:grh}, the two-base global exponent improves to $13/20+\epsilon$ for each $\epsilon>0$.

Question~\ref{q:local}, stated at the beginning of the introduction, concerns what base 3 adds to base 2 for a given pair of primes, and here the ratio of the factors matters. On a fixed ray both tests are power-residue conditions at the two primes, and reciprocity relates them. Theorem~\ref{intro:rays} shows that on five rays the answer is an explicit congruence, that base 3 is automatic on the ray $(1,2)$, and that, in the population with a prime larger factor and a sifted companion, the two local indicators defined in Section~\ref{sec:criba} are asymptotically independent.

On a ray we write $p=ag+1$, $q=bg+1$, with coprime positive integers $a<b$ and $g=\gcd(p-1,q-1)$. For prime factors greater than 3, let $Q_{a,b}^{(2)}(X)$ count products $pq\le X$ on the ray that pass base 2, and let $Q_{a,b}^{(2,3)}(X)$ count those that pass both bases. The superscripts identify the tested bases.

\begin{introtheorem}[Base-3 acceptance on five rays]\label{intro:rays}
On each of the five rays in the table, suppose that $p,q>3$ are prime and $pq$ passes base 2. Write $q$ in the indicated form with integers $h$ and $k>0$. Then $pq$ passes base 3 exactly when the additional condition in its row holds.
\[
\begin{array}{c|l|l|c}
(a,b)&q&\text{additional condition}&\rho_{a,b}\\ \hline
(1,3)&(3h+1)^2+27k^2&3\mid k&1/3\\
(2,3)&(6h+1)^2+108k^2&3\mid h,\ 3\mid k&1/6\\
(1,4)&(4h+1)^2+64k^2&3\mid k&1/4\\
(1,6)&(6h+1)^2+432k^2&3\mid k&1/3\\
(1,8)&(8h+1)^2+256k^2&8h+1\equiv c\pmod3&1/4
\end{array}
\]
In the last row, $q=c^2+2d^2$ with $c\equiv1\pmod4$ and $d>0$. On the additional ray $(1,2)$, every prime pair with factors greater than 3 passes base 3, and it passes base 2 precisely when $4\mid g$.

For each of the five rows, $\rho_{a,b}$ is an unconditional limiting proportion in the prime-with-sifted-companion problem specified in Theorem~\ref{t:criba}. Under the explicit prime-value hypotheses of Theorem~\ref{t:pares}, it is also the limiting ratio
\[
 Q_{a,b}^{(2,3)}(X)/Q_{a,b}^{(2)}(X)\longrightarrow\rho_{a,b}.
\]
More generally, for every ray $(a,b)$ with $a\le2$ the sifted proportion exists and equals $1/b$, $2/b$ or $1/(2b)$, according to the parity of $b$ and its residue modulo $3$.
\end{introtheorem}

Theorem~\ref{intro:rays} consists of Theorem~\ref{t:exacto}, which gives the criteria in coordinates, Proposition~\ref{pr:autom}, which treats the ray $(1,2)$, Theorem~\ref{t:criba} and Proposition~\ref{pr:allrays}, which give the sifted proportions, and Theorem~\ref{t:pares}, which gives the proportions for prime pairs. Character conditions of this kind appear for the ratios $4$ and $3$ in the work of Zhang \cite{Zhang01} and for the ratios $3/2$ and $6$ in that of Zhang and Tang \cite{ZT03}, who used them to search for strong pseudoprimes. On $(1,8)$ the congruence itself is Western's classical criterion as recorded by Hudson \cite[criterion (a)]{Hudson}. We specialize it to a prime pair whose base-2 coordinates satisfy $A=8h+1$, and we use it alongside the four other coordinate tests. Our distribution results count these conditions after sifting the companion. In Proposition~\ref{pr:allrays}, $\rho_{1,b}=1/b$ if $b$ is odd or $b\equiv1\pmod3$ and $\rho_{1,b}=2/b$ otherwise, while $\rho_{2,b}=1/b$ if $b\equiv1\pmod3$ and $\rho_{2,b}=1/(2b)$ otherwise; the six rays above are instances of this formula. In the sifted problem the larger factor $q$ is prime, and the companion $a(q-1)/b+1$ is required to have no prime factor at most the sieve cutoff. Theorem~\ref{t:criba} specifies the range of the cutoff and proves both the limiting proportions and asymptotic independence of the two local conditions. Corollary~\ref{co:seis} combines the six rays under its stated prime-value hypotheses to give a lower limit of $3.7842146471\ldots$ for $\Qq_{23}(X)(\log X)^2/\sqrt X$.

For fixed $k\ge2$ let $\mathcal E_k(Z)$ count the products $n\le Z$ of $k$ distinct odd primes at which $2$ has a common order, and let $\mathrm{Ov}_2(Z)$ count the odd composite $n\le Z$ with $\ord_e(2)=\ord_n(2)$ for every divisor $e>1$ of $n$. These are the overpseudoprimes of \cite{Shev12}, and every product counted by $\mathcal E_2(Z)$ is also counted by $\Qq_2(Z)$.

\begin{introtheorem}[Equal-order products and overpseudoprimes]\label{intro:over}
For each fixed $k\ge2$ and all sufficiently large $Z$,
\[
 \mathcal E_k(Z)\ll_k\frac{Z^{(k+1)/(2k)}}{L_Z^{k+\dd}H_Z^{3/2}},
 \qquad
 \mathrm{Ov}_2(Z)\ll\frac{Z^{3/4}}{L_Z^{2+\dd}H_Z^{3/2}}.
\]
\end{introtheorem}

Shevelev \cite[Theorem~3 and Remark~1]{Shev08} proved $\mathrm{Ov}_2(Z)\le Z^{3/4}(1+O(H_Z^2/L_Z))$, and his argument gives $O_k(Z^{(k+1)/(2k)}H_Z^{2k})$ for the overpseudoprimes with exactly $k\ge3$ prime factors. Theorem~\ref{intro:over} saves more than $L_Z^2$ in the total count, and for each $k$ it replaces the factor $H_Z^{2k}$ by the saving $L_Z^{k+\dd}H_Z^{3/2}$. For $k=2$ most of the saving already follows from Theorem~\ref{intro:rank}, because $\mathcal E_2(Z)\le\Qq_2(Z)$. Theorem~\ref{intro:over} consists of Theorem~\ref{r15ps:equal} and Corollary~\ref{cor:ovtotal}. In the other direction, the strong family of Theorem~\ref{th:lower} consists of products of two primes of equal order, so that $\mathcal E_2(Z)\ge(0.467+o(1))\log Z$.

The exponent $3/4$ is common to Theorem~\ref{intro:rank} for one base, to Theorem~\ref{intro:over} and to Shevelev's bound. Corollary~\ref{co:erdos} shows what a smaller exponent requires: each of $\Qq_2(X)$, $\mathcal E_2(X)$ and $\mathrm{Ov}_2(X)$ is $O(X^{3/4-\kappa})$ for some $\kappa>0$ if and only if the number of primes $p\le y$ with $\ord_p(2)\le p^{1/2+\delta}$ is $O(y^{1-\eta})$ for some $\delta,\eta>0$. This condition controls the exceptional set at one fixed exponent above $1/2$. It is related to the conjecture of Erdős \cite{Erdos76} that $\ord_p(2)>p^{c}$ for almost all $p$ and every $c<1$, in which the exponent $1/2$ of \cite{EM} is the best known. The same corollary shows that Pomerance's conjecture $\mathrm{Ov}_2(X)=o(X^{1/2+\epsilon})$, recorded in \cite[Remark~2]{Shev08}, implies Erdős's conjecture, with at most $O(y^{(1+c)/2+\epsilon})$ exceptional primes up to $y$.

Theorem~\ref{th:lower} gives explicit unconditional lower bounds of logarithmic size for base-2 products that pass the strong test and for those that fail it. Luca, Ralaivaosaona and Jiménez Urroz \cite[Propositions 3.1 and 3.2]{LRU} gave a global lower bound of logarithmic size; ours has explicit constants and separates the two families. Under the hypothesis on primes of a given exact order implied by Kovač and Luca's conjecture, Theorem~\ref{th:KLsharp} gives upper bounds of order $z(\log z)^\lambda$ in the window and $\sqrt X(\log X)^\lambda$ globally, with $\lambda=0.9105498551\ldots$. Remark~\ref{re:KLartin} shows that this hypothesis already implies Artin's conjecture for the base $2$ in its qualitative form, so these bounds rest on a strong input. Combined with Proposition~\ref{pr:lowerbox} and Corollary~\ref{co:ventdoble}, which assume prime-value hypotheses, this yields the conditional growth exponents $1$ in the window and $1/2$ globally for each of the one-base and two-base counts. The global exponent $1/2$ is the one in Galway's conjectured asymptotic $C\sqrt X/(\log X)^2$ \cite{Galway,GP}. Since $\mathcal E_2\le\Qq_2$, the same hypothesis gives Pomerance's conjectured bound $o(X^{1/2+\epsilon})$, recorded in \cite[Remark~2]{Shev08}, for the overpseudoprimes with two prime factors.

Theorem~\ref{th:random} treats a related sampling problem. For a uniformly chosen modulus in the window $z<p<q\le2z$, any fixed positive number of independently and uniformly chosen unit bases on that modulus pass with mean probability $z^{-1+o(1)}$. The exponent holds for both the Fermat and strong tests.

\subsection{Outline of the proofs}

We describe the main steps of the proofs in the order of the sections.

The Chinese remainder theorem gives the order conditions at the two prime factors in Section~\ref{sec:criterio}. In Section~\ref{sec:superior}, multiplicative relations among the generators bound the number of primes with small reduction order. The order at the larger prime then restricts the smaller prime to a progression. Ford's shifted-prime divisor estimate saves logarithms near the cutoff between these two estimates. Summing over dyadic rectangles under $pq=X$ gives the global bound, including pairs of very different sizes.

Section~\ref{sec:momentos} uses Cauchy--Schwarz on a matrix indexed by the two orders to control the semiprime count by the second moment of the order multiplicities. Pairs with equal orders give the reverse bound. To transfer a residual-index moment to the semiprime count, we balance the moment bound against the small-order estimate. Under GRH, the degrees and discriminants of the relevant Kummer fields control the contribution of small indices through the prime ideal theorem, while the small-order bound controls the large indices. Section~\ref{sec:fibras} uses a coefficient-uniform sieve for pairs of linear forms and a moment bound for the divisor function to optimize the logarithmic exponent under the exact-order hypothesis. Section~\ref{sec:equal} bounds products of primes of equal order by the method of Section~\ref{sec:superior}, now with one order shared by all factors, and the same sieve bounds equal-order pairs with restricted indices. Section~\ref{sec:inferior} uses primitive divisors with prescribed orders to construct disjoint strong and nonstrong families, and Bateman's inverse-totient asymptotic gives their lower-bound constants.

On the five rays in Section~\ref{sec:rayos}, cubic, quartic and octic reciprocity turn the order conditions into congruences on the coordinates of prime representations. The ray $(1,2)$ follows from a direct calculation. Section~\ref{sec:criba} expresses the local conditions on the other five rays as Frobenius classes in fixed fields. Their distribution in arithmetic progressions allows us to apply the fundamental lemma of sieve theory to the companion. For a companion that is itself prime, Section~\ref{sec:pares} uses explicit instances of Bateman--Horn and compares the local factors. The octic case uses a norm from $\Q(\zeta_8)$, where $\zeta_m=e^{2\pi i/m}$. Section~\ref{sec:constantes} computes the constants from the singular series and the volumes determined by the product cutoff and the dyadic window.

Section~\ref{sec:rsa} first normalizes the fixed-base counts by the size of the sample of prime pairs. It then counts the unit bases accepted by the two tests at a fixed modulus. Averaging those counts uses a two-form sieve for the upper bound and admissible linear forms with at least two prime values for the lower bound. Under Hardy--Littlewood, the uniform upper bound lets us sum the leading terms over all slopes. Section~\ref{sec:datos} checks the criteria and the constants by finite computation.

\section{Common orders and counting conventions}\label{sec:criterio}

Every result of the paper starts from one observation: whether $pq$ passes the Fermat test depends only on the multiplicative orders of the bases modulo $p$ and modulo $q$. This section records that criterion and its analogue for the strong test (Lemmas~\ref{le:crt} and~\ref{le:fuerte}), and fixes the counting conventions used throughout.

We retain the group $\Gamma$, the rank $r\ge1$ and the good-reduction convention from Section~\ref{sec:intro}. Throughout this section $p<q$ are distinct odd primes of good reduction and $n=pq$. In this notation the counting functions of Section~\ref{sec:intro} are
\begin{align}
 \A_\Gamma(z)&=\#\{z<p<q\le2z:\ t^{pq-1}\equiv1\pmod{pq}\ \ (\forall t\in\Gamma)\},\label{eq:A}\\
 \Qq_\Gamma(X)&=\#\{p<q,\ pq\le X:\ t^{pq-1}\equiv1\pmod{pq}\ \ (\forall t\in\Gamma)\}.\label{eq:Q}
\end{align}
It suffices to test a generating set. For $\Gamma=\langle t\rangle$ we write $\A_t$ and $\Qq_t$, extending the subscripts $2$ and $23$ of Section~\ref{sec:intro}. For a ray $(a,b)$, let $N_{a,b}(X)$ count all prime products $(ag+1)(bg+1)\le X$ with $g>0$ even, $\gcd(a,b)=1$ and both primes greater than 3, and let $Q_{a,b}^{(3)}(X)$ count those that pass base 3; together with $Q_{a,b}^{(2)}(X)$ and $Q_{a,b}^{(2,3)}(X)$ from Section~\ref{sec:intro}, these are the restrictions of $N_{a,b}(X)$ to acceptance in the indicated bases. In particular, the superscript $(2)$ means ``base 2'', not a square, and $a,b$ are ray coefficients, not test bases.

We write $\Pp$ for the set of primes and use $\varphi(n)=\#(\Z/n\Z)^\times$ for Euler's totient, $\tau(n)$ for the number of positive divisors of $n$, $\omega(n)$ for the number of distinct prime factors, and $\Omega(n)$ for their number with multiplicity. The Möbius function $\mu(n)$ is zero if a prime square divides $n$, and is $(-1)^{\omega(n)}$ otherwise. The polynomial $\Phi_m(T)$ is the $m$th cyclotomic polynomial. We write $\pi(y)$ for the number of primes at most $y$, $\pi(y;k,a)$ for the number of primes $p\le y$ with $p\equiv a\pmod k$, and $\operatorname{Li}(y)=\int_2^y dt/\log t$ for the logarithmic integral. For an odd prime $p$ and an integer $t$ with $p\nmid t$, $(t/p)$ denotes the Legendre symbol: $1$ for a square modulo $p$ and $-1$ otherwise. For each prime power $u$, let $\F_u$ denote the field with $u$ elements. The notation $(\F_p^\times)^k$ denotes the subgroup of $k$th powers. For a prime $\ell$ and a nonzero rational $u$, $v_\ell(u)$ denotes its $\ell$-adic valuation: the exponent of $\ell$ in its numerator minus that in its denominator. As usual, $(u,v)=\gcd(u,v)$ for integers $u,v$, and $[u,v]=\lcm(u,v)$. Write $\Gamma\bmod p$ for the image of $\Gamma$ in $\F_p^\times$. Put
\[
 h_p=|\Gamma\bmod p|,\qquad I_p=(p-1)/h_p,\qquad g=\gcd(p-1,q-1).
\]
We call $I_p$ the residual index of the group. For a single base $t$, we have $h_p=\ord_p(t)$.

\begin{lemma}[Classical common-order criterion]\label{le:crt}
The product $pq$ passes Fermat for every base in $\Gamma$ if and only if $h_p\mid q-1$ and $h_q\mid p-1$. Equivalently, both $h_p$ and $h_q$ divide $g$. In the single-base case $\Gamma=\langle t\rangle$ with integer $t\ge2$, this is equivalent to $pq\mid t^g-1$. With the additional size condition $pq\le X$, these are exactly the products counted in \eqref{eq:Q}.
\end{lemma}
\begin{proof}
By the Chinese remainder theorem, $t^{pq-1}\equiv1\pmod{pq}$ if and only if the congruence holds modulo $p$ and
modulo $q$. Since $pq-1=q(p-1)+(q-1)$ and $t^{p-1}\equiv1\pmod p$, the congruence modulo $p$ is equivalent to
$t^{q-1}\equiv1\pmod p$. The group $\Gamma\bmod p$ is a subgroup of the cyclic group $\F_p^\times$, hence cyclic of
order $h_p$, so the congruence holds for every $t\in\Gamma$ if and only if $h_p\mid q-1$. As $h_p\mid p-1$, this is
equivalent to $h_p\mid g$. The same argument applies modulo $q$. For an integer base $t$ the two conditions say that
$t^g\equiv1$ modulo $p$ and modulo $q$, that is, $pq\mid t^g-1$.
\end{proof}
The criterion is classical; for base 2 it is stated and derived in \cite[\S1, equation (1)]{LRU}. We record it
with its short proof, without claiming novelty.

For $\Gamma=\langle2,3\rangle$ it gives
\[
 h_p=\lcm(\ord_p(2),\ord_p(3)),\qquad
 pq\mid\gcd(2^g-1,3^g-1).
\]
If $p=ag+1$ and $q=bg+1$ with $\gcd(a,b)=1$, the two subgroup conditions become $a\mid I_p$ and $b\mid I_q$. Indeed, the subgroup of order $g$ in $\F_p^\times$ has index $a$, and it contains $\Gamma\bmod p$ precisely when $a\mid I_p$; the second condition is identical with $p,q$ exchanged. Globally we count the possible orders and partners. On a fixed ray we test membership in these subgroups by characters and reciprocity.

For clarity, we spell out three elementary size consequences for an \emph{integer test base} $t\ge2$. They do not constitute a separate new estimate. If $pq$ passes this base, Lemma~\ref{le:crt} gives
\[
 1\le pq\le t^g-1,\qquad g\le p-1,\qquad q-p\le q=\frac{pq}{p}.
\]
Write $\log_t y=\log y/\log t$. Taking logarithms in the first inequality and substituting into the other two yields
\[
 g\ge\log_t(pq+1),\qquad
 p\ge1+\log_t(pq+1),\qquad
 q-p\le\frac{pq}{1+\log_t(pq+1)}.
\]
If $t$ is a primitive root modulo the larger prime, meaning $\ord_q(t)=q-1$, the criterion excludes the pair because $q-1>p-1$.

\begin{lemma}[Classical strong-test criterion]\label{le:fuerte}
Let $t$ be a unit modulo $n=pq$ that passes Fermat. Set $s=v_2(n-1)$ and $e=(n-1)/2^s$, so $s\ge1$ and $e$ is odd. The strong probable-prime test accepts if $t^e\equiv1\pmod n$ or if $t^{2^j e}\equiv-1\pmod n$ for some integer $0\le j<s$. It accepts if and only if
\[
 v_2(\ord_p(t))=v_2(\ord_q(t)).
\]
\end{lemma}
\begin{proof}
Put $a_p=v_2(\ord_p(t))$ and $a_q=v_2(\ord_q(t))$. Since $t$ passes Fermat, $\ord_p(t)$ and $\ord_q(t)$ divide
$n-1=2^se$, so $a_p,a_q\le s$ and the odd parts of both orders divide $e$. Hence $t^e$ has order $2^{a_p}$ modulo
$p$ and order $2^{a_q}$ modulo $q$, and $t^e\equiv1\pmod n$ exactly when $a_p=a_q=0$. In the cyclic group
$\F_p^\times$ the only element of order $2$ is $-1$, so for $0\le j<s$ we have $t^{2^je}\equiv-1\pmod p$ exactly
when $a_p=j+1$, and similarly modulo $q$. By the Chinese remainder theorem, $t^{2^je}\equiv-1\pmod n$ exactly when
$a_p=a_q=j+1$. As $a_p,a_q\le s$, the test accepts if and only if $a_p=a_q$.
\end{proof}
This is the one-base, two-prime specialization of Jaeschke's signature criterion \cite[Proposition 1]{Jaeschke}; see also the definition of the strong test in \cite[p.~1004]{PSW}. The integer $s$ is the exponent of 2 in $n-1$ and has no relation to the ray coefficients $a,b$. We use Lemma~\ref{le:crt} in the global bounds and the ray criteria, and Lemma~\ref{le:fuerte} in the truncated-fiber construction, the two lower-bound families, and the count of accepting unit bases in Section~\ref{sec:rsa}.

\subsection{Relation to conventions allowing squares or nonunit bases}\label{sec:conventions}

Some results in the literature count all semiprimes or use the congruence $t^n\equiv t\pmod n$. This subsection shows how our counts differ from those conventions, so that the comparisons in Sections~\ref{sec:superior} and~\ref{sec:equal} are made on the same footing.

The restriction to distinct odd factors is part of our definition, so it must be taken into account when comparing a result stated for all semiprimes or with the extended Fermat congruence $t^n\equiv t\pmod n$. For a unit $t$ the extended congruence is equivalent to $t^{n-1}\equiv1\pmod n$. The remaining possibilities are prime squares and factors dividing the base.

For an odd prime $p\nmid t$, the square $p^2$ passes Fermat precisely when
\[
 t^{p-1}\equiv1\pmod{p^2}.
\]
This is the Wieferich condition for the base $t$. To check the equivalence directly, write $t^{p-1}\equiv1+cp\pmod{p^2}$, which is possible by Fermat's little theorem. Since $p^2-1=(p-1)(p+1)$,
\[
 t^{p^2-1}\equiv(1+cp)^{p+1}\equiv1+cp\pmod{p^2}.
\]
Thus the square passes exactly when $p\mid c$. For bases 2 and 3 simultaneously, a square with $p>3$ must satisfy both $2^{p-1}\equiv1\pmod{p^2}$ and $3^{p-1}\equiv1\pmod{p^2}$. Allowing all prime squares adds at most $\pi(\sqrt X)=O(\sqrt X/\log X)$ terms. This does not change the upper bounds in Section~\ref{sec:superior}. It is not negligible a priori at the conjectural scale $\sqrt X/(\log X)^2$, so the fixed-ray asymptotics and all numerical tables retain the distinct-factor convention.

For base 2 the extended convention introduces no even semiprimes. If $n=2q$ with $q$ an odd prime, Fermat's little theorem gives $2^{2q}\equiv4\pmod q$, whereas acceptance would require $2^{2q}\equiv2\pmod q$. This is impossible. Also $2^4\not\equiv2\pmod4$. These observations specify exactly how our counts compare with the base-2 convention of \cite{LRU}; no even-factor contribution is hidden in that comparison.

\section{Upper bounds for a fixed group of bases}\label{sec:superior}

This section proves Theorem~\ref{intro:rank}, which answers Question~\ref{q:count} of the introduction with an upper bound valid for every ratio of the factors and every fixed group of bases. We use the notation $L_y$, $H_y$, $\dd$, $\beta$ and
$\wt_r$ introduced before that theorem. Section~\ref{ssec:elem} gives an elementary argument for one base.
Section~\ref{ssec:small} collects the estimates for primes with small reduction order. Section~\ref{ssec:rect}
states Theorem~\ref{intro:rank} as Theorem~\ref{th:main} and Corollary~\ref{cor:23} and proves them, and
Section~\ref{ssec:compare} compares the bounds with others.

\subsection{An elementary argument for one base}\label{ssec:elem}

Before the proof of Theorem~\ref{intro:rank} we show what an elementary argument gives. In the window it already improves the exponent $8/5$ of \cite{LRU} to $3/2$; Theorem~\ref{intro:rank} adds the logarithmic factor coming from Ford's theorem.

For a fixed integer $t\ge2$ let $M_z(d)$ count the primes $z<p\le2z$ that are $1\bmod d$ and divide $t^d-1$.
Then
\[
 \A_t(z)\le\sum_{d\le 2z}\binom{M_z(d)}{2},\quad
 M_z(d)\le \frac{d\log t}{\log z},\quad
 M_z(d)\le\frac{2z}{d}.
\]
For the first inequality, every pair contributes at least at its greatest common divisor. The second follows
by multiplying the distinct primes, and the third by counting the positive integers $(p-1)/d$. Splitting
at $d=\sqrt z$ already gives the exponent $3/2$. More precisely, for $\sqrt z<d\le z^{3/4}$ the Brun--Titchmarsh
inequality \cite{MV} gives $M_z(d)\ll z/(\varphi(d)L_z)$. The sum over $d>\!z^{3/4}$ is $O(z^{5/4})$. Since
$\sum_{d>Y}\varphi(d)^{-2}\ll 1/Y$, this proves $\A_t(z)\ll_t z^{3/2}/L_z^2$. To justify the last tail estimate,
observe that $(d/\varphi(d))^2$ has bounded mean. Expand it over squarefree divisors, with local coefficient
$(1-1/\ell)^{-2}-1=O(1/\ell)$, and use the convergent product of $1+O(1/\ell^2)$. Dyadic summation gives
the claim. The logarithmic improvement below uses the divisor $h_q\mid q-1$ in the transition range and
the cross-divisibility $h_q\mid p-1$ to count the companion primes.

\subsection{Small orders and divisors of shifted primes}\label{ssec:small}

For a single base the transition occurs near the square root of the prime. In rank $r$ the transition occurs near $q^{r/(r+1)}$, up to logarithmic factors. Lemma~\ref{le:small} bounds the number of primes with small reduction order, Lemma~\ref{le:ford} is the input
from Ford's work, weighted by powers of $(q-1)/\varphi(q-1)$, and Proposition~\ref{pr:rec} combines them into a bound for $\sum1/\varphi(h_q)$. The proof
of Theorem~\ref{th:main} in Section~\ref{ssec:rect} sums this bound over dyadic rectangles.

Throughout this subsection $Q<q\le2Q$ runs over primes of good reduction. Write $L=L_Q$, $H=H_Q$
and $\alpha=r/(r+1)=1-\beta$.

\begin{lemma}\label{le:small}

For $T\ge1$ we have
\[
 \sum_{h_q\le T}\log q\ll_\Gamma T^{1+1/r},\qquad
 \#\{q:h_q\le T\}\ll_\Gamma T^{1+1/r}/L .
\]
\end{lemma}

\begin{proof}

The small-reduction-group estimate goes back to Matthews and to Gupta and Murty \cite{Matt,GM}. Choose multiplicatively
independent elements $g_j=A_j/B_j$ in lowest terms. For $K=\lceil T^{1/r}\rceil$ the box $\{0,\ldots,K\}^r$
contains more than $T$ points. If $h_q\le T$, two monomials have the same image. Hence some nonzero vector
$v$ with $|v_j|\le K$ satisfies
\[
 q\mid D_v,\qquad
 D_v=\prod_j A_j^{v_j^+}B_j^{v_j^-}
       -\prod_j A_j^{v_j^-}B_j^{v_j^+},
 \qquad v_j^\pm=\max(\pm v_j,0).
\]
By independence each $D_v$ is a nonzero integer, and $\log |D_v|\ll_\Gamma K$. The product over the $O_r(K^r)$
possible vectors is divisible by every $q$ counted. Taking logarithms and using $q>Q$ proves both
bounds.

\end{proof}

\begin{lemma}\label{le:ford}

Fix an integer $k\ge0$ and $0<c<1/2$. For $Q^c\le D\le Q^{1-c}$ we have
\[
 \sum_{\substack{Q<q\le2Q\\ \exists f\mid q-1:\ D<f\le2D}}
   \Bigl(\frac{q-1}{\varphi(q-1)}\Bigr)^{k}
 \ll_{k,c} \frac{Q}{L^{1+\dd}H^{3/2}} ,
\]
where $q$ runs over primes. In particular, $\#\{q:D<h_q\le2D\}\ll_c Q/(L^{1+\dd}H^{3/2})$.

\end{lemma}

\begin{proof}

We first need Ford's estimate for shifted primes with a coefficient. Fix $0<c'<1/2$. For integers $e\ge1$ and
$x^{c'}\le y\le x^{1-c'}$, let $\mathcal N_e(x,y)$ count the integers $x/2<n\le x$ such that $en+1$ is prime and
$n$ has a divisor in $(y,2y]$. We claim that, uniformly in $e$,
\begin{equation}
 \mathcal N_e(x,y)\ll_{c'}\frac{e}{\varphi(e)}\cdot\frac{x}{(\log x)^{1+\dd}(\log\log x)^{3/2}} .
 \label{eq:forde}
\end{equation}
It suffices to treat $y\le\sqrt x$. Indeed, if $y>\sqrt x$, a divisor of $n$ in $(y,2y]$ has its complementary
divisor in $(x/(4y),x/y]$. This interval is covered by two intervals $(y',2y']$ with $x^{c'/2}\le y'\le\sqrt x$
for large $x$.

For $e=1$ the integers $n$ are shifted primes $q-1$, and \eqref{eq:forde} is Ford's Theorem 6 with $\lambda=-1$
combined with his Corollary 2 \cite{Ford}. Ford's Theorem 6 bounds the number of shifted primes with a divisor in
$(y,2y]$ by the number of all integers $n\le x$ with such a divisor, divided by $\log x$. Ford's Corollary 2 shows that
the latter number has order $x/((\log Y)^{\dd}(\log\log Y)^{3/2})$ with $Y=\min(y,x/y)$, and here
$\log Y\asymp\log x$.

For general $e$ we follow the proof of Theorem 6 in \cite[\S14]{Ford}. Since $y\ge x^{c'}$, we are in the case
$z>\exp((\log x)^{9/10})$ of the proof of Ford's Lemma 14.1 \cite{Ford}, with $z=2y$. In that proof the primality of the
shifted integer is used in two places. The first is the passage to the integers of a given interval whose shift
is prime. We replace them by the integers $n$ with $en+1$ prime. The other steps of the proof bound sets of
integers defined without reference to primes, so their contribution does not change. The second place is a sieve
bound for a cofactor. Write $P^+(m)$ and $P^-(m)$ for the largest and the smallest prime factor of an integer
$m>1$, and put $P^+(1)=1$ and $P^-(1)=\infty$. One writes $n=arb$ with $r$ prime, $P^+(a)<r$ and $P^-(b)>r$, fixes $a$ and $r$, and counts
the admissible $b$ in an interval of length $N$. Ford's parameter in this step, which we call $Q'$, satisfies $r^{1/2}\le Q'\le N$.

We sieve only by odd primes $\ell$. Put $U={Q'}^{1/10}$ and $V=\min(r,Q')^{1/10}$. We remove the class
$b\equiv0\bmod\ell$ for $\ell\le V$, which is allowed since $P^-(b)>r>V$. We also remove the root of
$earb+1\equiv0\bmod\ell$ for $\ell\le U$, which is allowed since $earb+1=en+1>x/2$ is a prime larger than $U$.
The linear form has one root modulo $\ell$ if $\ell\nmid ear$, and that root is not zero. It has no root if
$\ell\mid ear$. For $e=1$ these are Ford's classes. Let $\nu_\ell\le2$ be the number of removed classes, and put
\[
 \mathcal Z=\prod_{2<\ell\le U}\Bigl(1-\frac{\nu_\ell}{\ell}\Bigr)^{-1}=\sum_m \psi(m),\qquad
 \psi(m)=\mu^2(m)\prod_{\ell\mid m}\frac{\nu_\ell}{\ell-\nu_\ell},
\]
where $m$ runs over products of odd primes $\ell\le U$. The inequality $(1-2/\ell)^{-1}\ge(1-1/\ell)^{-2}$ and
Mertens' theorem give
\[
 \mathcal Z\gg \log U\,\log V\,\frac{\varphi(ear)}{ear}.
\]
Moreover $\sum_m \psi(m)\log m=\mathcal Z\sum_{2<\ell\le U}\nu_\ell(\log\ell)/\ell\le(1/5+o(1))\mathcal Z\log Q'$. Put $R=\sqrt{Q'}$.
Since $\log m>\log R$ for $m>R$, the sum $G(R)=\sum_{m\le R}\psi(m)$ satisfies $G(R)\ge \mathcal Z/2$ for large $Q'$; bounded
values of $Q'$ only affect the constant. The arithmetic large sieve \cite{MV} bounds the number of remaining $b$
by $(N+R^2)/G(R)\le4N/\mathcal Z$. Since $\log U\asymp\log Q'$, $\log V\asymp\log r$ and
$ear/\varphi(ear)\le(e/\varphi(e))(a/\varphi(a))\,r/(r-1)$, the number of $b$ is
\[
 \ll\frac{e}{\varphi(e)}\cdot\frac{a}{\varphi(a)}\cdot\frac{N}{\log r\,\log Q'}.
\]
For $e=1$ this is the bound in Ford's proof, and for general $e$ it carries the additional factor $e/\varphi(e)$.
The rest of the proof of Ford's Lemma 14.1, his deduction of Theorem 6 from it and his Corollary 2 \cite{Ford} add these bounds over
$a$ and $r$ with nonnegative weights that do not depend on $e$. This proves \eqref{eq:forde}.

We now insert the weight. Let $w_k$ be the multiplicative function with $(m/\varphi(m))^k=\sum_{e\mid m}w_k(e)$.
It is supported on squarefree integers, and $0\le w_k(\ell)=(1-1/\ell)^{-k}-1\le C_k/\ell$ at primes $\ell$.
The sum $S$ in the lemma therefore equals $\sum_e w_k(e)S_e$, where $S_e$ counts the primes $Q<q\le2Q$ with
$e\mid q-1$ for which $q-1$ has a divisor in $(D,2D]$. Put $E=Q^{c/4}$. The integers $e>E$ contribute at most
\[
 \sum_{e>E}w_k(e)\frac{2Q}{e}\le\frac{2Q}{E^{1/2}}\sum_{e\ge1}\frac{w_k(e)}{e^{1/2}}\ll_k Q^{1-c/8},
\]
since the last series has a convergent Euler product. Now let $e\le E$ and write $q-1=en$. If $f\mid q-1$ with
$D<f\le2D$ and $d=\gcd(f,e)$, then $f/d$ divides $n$ and lies in $(D/d,2D/d]$. The integer $n$ lies in $(x/4,x]$
with $x=2Q/e$, and for large $Q$ the endpoint $D/d$ lies between $(x/2)^{c/2}$ and $(x/2)^{1-c/2}$. We apply
\eqref{eq:forde} with $c'=c/2$ on the intervals $(x/4,x/2]$ and $(x/2,x]$ and sum over the divisors $d$ of $e$.
Since $\log x\asymp\log Q$, this gives
\[
 S_e\ll_c \tau(e)\,\frac{e}{\varphi(e)}\cdot\frac{Q}{e\,L^{1+\dd}H^{3/2}}
 =\frac{\tau(e)}{\varphi(e)}\cdot\frac{Q}{L^{1+\dd}H^{3/2}}.
\]
The series $\sum_e w_k(e)\tau(e)/\varphi(e)$ converges, since its Euler factors are $1+O_k(\ell^{-2})$. This proves
the bound for $S$. For the last assertion take $k=0$ and use $h_q\mid q-1$.

\end{proof}

\begin{lemma}\label{le:smooth}

For every $1\le m\le 2Q$ we have $m/\varphi(m)\ll H$.

\end{lemma}

\begin{proof}

This is the classical maximal order of $m/\varphi(m)$. Separate the prime divisors of $m$ at $\log(2Q)$. By
Mertens' theorem the primes below the cutoff contribute a factor $O(H)$. There are at most
$\log(2Q)/\log\log(2Q)$ prime divisors above the cutoff, and their Euler factors have bounded product.

\end{proof}

\begin{prop}\label{pr:rec}

For fixed $\Gamma$ we have
\[
  \sum_{Q<q\le2Q}\frac1{\varphi(h_q)}
       \ll_\Gamma \frac{Q^\beta}{L}\wt_r(Q).
\]
We shall also use the simpler bound $\sum_{Q<q\le2Q}1/h_q\ll_\Gamma Q^\beta/L$.

\end{prop}

\begin{proof}

Put $a_0=1+1/r$, so that $1-1/a_0=\beta$. For $T\ge1$ we use the following geometric sum over powers of
two $D$.
\begin{equation}
 \sum_D \frac{\min(D^{a_0},T)}D\ll_r T^\beta .
 \label{eq:geom}
\end{equation}
Below $T^{1/a_0}$ we sum $D^{a_0-1}$, and above it we sum $T/D$. With Lemma~\ref{le:small} and $\pi(2Q)\ll Q/L$,
this proves the second assertion on taking $T=Q$. An initial block covers any primes with $h_q=1$. No such
primes occur for large $Q$.

For the first assertion put $c=\beta/4$. By Lemma~\ref{le:smooth} and $h_q\mid q-1$, the ratio $h_q/\varphi(h_q)$
is uniformly $O(H)$. The blocks $D<Q^c$ contribute $O_\Gamma(H Q^{c/r}/L)$ by Lemma~\ref{le:small}. The
blocks $D>Q^{1-c}$ contribute $O(HQ^c/L)$ by the prime number theorem. Each saves a power of $Q$ relative
to the claimed bound.

In the intermediate blocks write $F=Q/(L^{\dd}H^{3/2})$ and $k=r+2$. Let $N_D$ be the number of $q$ with
$D<h_q\le2D$. Lemma~\ref{le:small} and Lemma~\ref{le:ford} with exponent $0$ give $N_D\ll_\Gamma\min(D^{a_0},F)/L$.
Since $h_q$ divides $q-1$, we have $h_q/\varphi(h_q)\le(q-1)/\varphi(q-1)$, and Lemma~\ref{le:ford} with exponent
$k$ bounds the sum of $((q-1)/\varphi(q-1))^k$ over these $q$ by $O_r(F/L)$. By Hölder's inequality,
\[
 \sum_{D<h_q\le2D}\frac{h_q}{\varphi(h_q)}
 \le N_D^{1-1/k}\Bigl(\sum_{D<h_q\le2D}\Bigl(\frac{q-1}{\varphi(q-1)}\Bigr)^{k}\Bigr)^{1/k}
 \ll_\Gamma \frac{F^{1/k}\min(D^{a_0},F)^{1-1/k}}{L}.
\]
Since $1/\varphi(h_q)<(h_q/\varphi(h_q))/D$ in the block, the block contributes
\[
 \ll_\Gamma \frac{F^{1/k}\min(D^{a_0},F)^{1-1/k}}{LD}.
\]
For $D^{a_0}\le F$ this is $F^{1/k}D^{a_0(1-1/k)-1}/L$, and
the exponent $a_0(1-1/k)-1=1/(r(r+2))$ is positive. The sum over these $D$ is therefore $O_r(F^{1-1/a_0}/L)$,
the size of its largest term. For $D^{a_0}>F$ the contribution is $F/(LD)$, and these blocks also give
$O(F^{1-1/a_0}/L)$. Since $1-1/a_0=\beta$ and $F^\beta=Q^\beta\wt_r(Q)$, this proves the first assertion.

\end{proof}

\subsection{The bound for a fixed group and its proof}\label{ssec:rect}

We now state Theorem~\ref{intro:rank} in the form in which we prove it, for an arbitrary fixed group of rank $r$, and prove it by summing Proposition~\ref{pr:rec} over dyadic rectangles.

\begin{theorem}[Theorem~\ref{intro:rank}, general form]\label{th:main}

For each fixed $\Gamma$ as in \eqref{eq:A}--\eqref{eq:Q} we have, unconditionally,
\[
 \A_\Gamma(z)\ll_\Gamma \frac{z^{1+\beta}}{L_z^2}\wt_r(z),
 \qquad
 \Qq_\Gamma(X)\ll_\Gamma \frac{X^{(1+\beta)/2}}{L_X^2}\wt_r(X).
\]
In particular, for every fixed integer base $t\ge2$, the factors multiplying $z^{3/2}$ and $X^{3/4}$ are,
in their respective arguments, $1/(L_x^{2+\dd/2}H_x^{3/4})$.

\end{theorem}

\begin{proof}[Proof of Theorem~\ref{th:main}]
Let $\B_\Gamma(P,Q)$ count admissible pairs with $P<p\le2P$, $Q<q\le2Q$ and $p<q$ that satisfy the two conditions
of Lemma~\ref{le:crt}. Assume $2\le P\le Q$. For fixed $q$ we use the necessary condition $p\equiv1\bmod h_q$.
Its solutions with $p\ge2$ are of the form $1+kh_q$ with $k\ge1$. Thus even for $h_q>P$ the number of choices
is $\le2P/h_q$. Proposition~\ref{pr:rec} gives
\begin{equation}
 \B_\Gamma(P,Q)\ll_\Gamma P Q^\beta/\log Q .
 \label{eq:rectrough}
\end{equation}
With $\alpha=r/(r+1)$, put $\eta=(1+\alpha)/2$ and $\theta=(\alpha+\eta)/2$. When $P\ge Q^\eta$ and $h_q\le Q^\theta$, Brun--Titchmarsh
up to $2P$ gives $O_r(P/(\varphi(h_q)\log Q))$ choices for $p$. For $h_q>Q^\theta$, counting integers gives
a total of $O(PQ^{1-\theta}/\log Q)$, which saves a power of $Q$ relative to the central bound since $1-\theta<\beta$.
Consequently
\begin{equation}
 \B_\Gamma(P,Q)\ll_\Gamma P Q^\beta\,\wt_r(Q)/L_Q^2
 \qquad (Q^\eta\le P\le Q).
 \label{eq:rect}
\end{equation}
Taking $P=Q=z$ proves the first part of Theorem~\ref{th:main}.

For the second part take $P,Q$ to be powers of two. If $pq\le X$ then $PQ\le X$. In boxes with $P<Q^\eta$
we have
\[
 P Q^\beta \le (PQ)^{(\eta+\beta)/(1+\eta)}
 \le X^{(\eta+\beta)/(1+\eta)},\qquad
 \frac{\eta+\beta}{1+\eta}<\frac{1+\beta}{2}.
\]
This follows by raising $P/Q^\eta\le1$ to the power $(1-\beta)/(1+\eta)$. Summing \eqref{eq:rectrough} over
the $O((\log X)^2)$ boxes of this type preserves a power saving. In every box with $PQ\le X^{1/2}$ we have
$P Q^\beta\le(PQ)^{(1+\beta)/2}\le X^{(1+\beta)/4}$. In the remaining boxes $Q>X^{1/4}$ and $Q\le X$, so
$\wt_r(Q)/L_Q^2\asymp_r
\wt_r(X)/L_X^2$. To sum \eqref{eq:rect}, observe that
\[
 \sum_{\substack{P,Q\ \text{dyadic}\\P\le Q,\ PQ\le X}}P Q^\beta
 \ll \sum_{Q\le\sqrt X}Q^{1+\beta}
        + X\sum_{\sqrt X<Q\le X}Q^{\beta-1}
 \ll_r X^{(1+\beta)/2}.
\]
The last two sums run over powers of two. This proves the global bound, including pairs with factors of
unequal sizes.
\end{proof}

\begin{cor}[Theorem~\ref{intro:rank} for bases 2 and 3]\label{cor:23}

Let $D(x)=L_x^{2+\dd/3}H_x^{1/2}$. For sufficiently large $X,z$,
\[
 \Qq_{23}(X)\ll X^{2/3}/D(X),\qquad \A_{23}(z)\ll z^{4/3}/D(z).
\]
In particular, $\Qq_{23}(X)\ll X^{2/3}/\log X$.

\end{cor}

\begin{proof}

The identity $2^i3^j=1$ over the rationals forces $i=j=0$ by the valuations at 2 and 3. Thus $r=2$ and $\beta=1/3$.
Substituting in Theorem~\ref{th:main} gives both estimates. The last assertion follows from
$D(X)/\log X\to\infty$.

\end{proof}

In rank one this proof gives the global exponent $3/4$. In rank three it gives $5/8$. The implied constant
depends on the fixed group.

\subsection{Comparison with other bounds}\label{ssec:compare}

We place Theorem~\ref{th:main} among the other bounds for the same counts.

For base 2 the changes of exponent $8/5\to3/2$ in the dyadic window and $4/5\to3/4$ under the product cutoff
are relative to Theorems 2.1 and 3.1 of \cite{LRU}. Each entry in the table is the right-hand side of an upper bound, up to an implied constant.

\begin{center}\small
\begin{tabular}{@{}lll@{}}
\toprule
Base-2 bounds & Window $\A_2(z)$ & Product $\Qq_2(X)$\\
\midrule
Luca, Ralaivaosaona and Jiménez Urroz \cite{LRU}
 & $z^{8/5}/L_z^2$ & $X^{4/5}/L_X^{2/5}$\\
Theorem \ref{th:main}
 & $z^{3/2}\wt_1(z)/L_z^2$ & $X^{3/4}\wt_1(X)/L_X^2$\\
Theorem \ref{th:KLsharp}, under the order-fiber hypothesis
 & $zL_z^\lambda$ & $\sqrt X L_X^\lambda$\\
\bottomrule
\end{tabular}
\end{center}

\paragraph{Comparison with the reciprocal-order bound.}

For a group generated by the first $r$ primes, Pappalardi \cite[Corollary 2.4]{Pappa} proved
\[
 \sum_{\substack{q\le Y\\ q\ {\rm of\ good\ reduction}}}\frac1{h_q}
 \ll_\Gamma \frac{Y^\beta}{(\log Y)^{1+\tau_P\beta}},
 \qquad \tau_P=0.0306.
\]
Sorenson and Webster use this form in \cite[Lemma 3.4]{SW}. Counting at most $2P/h_q$ candidates for the other factor in each rectangle gives
\[
 \B_\Gamma(P,Q)\ll_\Gamma \frac{PQ^\beta}{L_Q^{1+\tau_P\beta}} .
\]
Taking $P=Q=z$ or summing the dyadic rectangles with $PQ\le X$ gives
$\A_\Gamma(z)\ll_\Gamma z^{1+\beta}/L_z^{1+\tau_P\beta}$ and
$\Qq_\Gamma(X)\ll_\Gamma X^{(1+\beta)/2}/L_X^{1+\tau_P\beta}$, respectively. In the latter sum,
the rectangles with $PQ\le X^{1/2}$ save a power, and all remaining $Q$ have $\log Q\asymp\log X$.
In particular, the powers $3/4$ and $2/3$ for the global counts in bases $2$ and $2,3$ already follow from the reciprocal-order estimate.

One may add Brun--Titchmarsh on $P\ge Q^\eta$ for $h_q\le Q^\theta$, with the same thresholds as in
\eqref{eq:rect}. Using $h_q/\varphi(h_q)\ll H_Q$ then gives a factor $H_Q/L_Q^{2+\tau_P\beta}$
in the rectangle bound. The complementary ranges save a power as above, so this factor also holds at the
relevant argument $z$ or $X$ in the final counts. The factor $\wt_r/L^2$ in Theorem~\ref{th:main} improves
this comparison since $\dd>\tau_P$.

\section{Order collisions and moments of the residual index}\label{sec:momentos}

Theorem~\ref{intro:rank} stops at the exponent $(1+\beta)/2$, which is $3/4$ for one base. This section explains why. A smaller exponent is equivalent to a bound for pairs of primes with the same reduction order (Theorem~\ref{th:momento}) and to a bound for primes with small reduction order (Corollary~\ref{co:erdos}); for base 2 the latter is a power-saving form of a conjecture of Erdős. In rank at least two, the generalized Riemann hypothesis gives such a bound (Corollary~\ref{co:grh}).

\subsection{A necessary and sufficient condition for a power saving}

We first relate the exponent of $\Qq_\Gamma$ to the number of pairs of primes with the same reduction order.

For a fixed group $\Gamma$ as in Theorem~\ref{th:main}, let $\mathcal P_Y$ be the set of its good-reduction primes
in $(Y,2Y]$ and put
\[
 r_Y(s)=\#\{p\in \mathcal P_Y:h_p=s\},\qquad E_\Gamma(Y)=\sum_{s\ge1}r_Y(s)^2 .
\]
The sum counts ordered pairs of primes at which the reduction has the same order. The proof compares these equal-order counts with counts in which the two reduction orders have prescribed quotients.

\begin{theorem}\label{th:momento}

For every $\alpha\ge1$ the following estimates are equivalent.
\[
 \Qq_\Gamma(X)\ll X^{\alpha/2+o(1)},\qquad
 E_\Gamma(Y)\ll Y^{\alpha+o(1)} .
\]
In the right-to-left implication it suffices to assume the estimate for dyadic $Y$. The notation $o(1)$
allows an arbitrarily small fixed loss in the exponent.

\end{theorem}

\begin{proof}

If $h_p=h_q=s$, then $s$ divides both $p-1$ and $q-1$, so Lemma~\ref{le:crt} applies. This gives $E_\Gamma(Y)\le |\mathcal P_Y|+2\Qq_\Gamma(4Y^2)$.
The diagonal is absorbed when $\alpha\ge1$. Conversely, define the finite matrix
\[
 B_Y(s,t)=\#\{p\in \mathcal P_Y:h_p=s,\ t\mid p-1\}.
\]
The number of ordered pairs in $\mathcal P_P\times\mathcal P_Q$ satisfying the divisibilities, including the diagonal when
present, is $\sum_{s,t}B_P(s,t)B_Q(t,s)$. Cauchy--Schwarz bounds it by $\|B_P\|_2\|B_Q\|_2$, where the norm
is the square root of the sum of squares of the entries. Since $B_Y(s,t)\le r_Y(s)$,
\[
 \|B_Y\|_2^2\le\sum_{p\in \mathcal P_Y}r_Y(h_p)\tau(p-1)
 \le\max_{n\le 2Y}\tau(n)\,E_\Gamma(Y)=Y^{o(1)}E_\Gamma(Y).
\]
Here $\tau(n)$ is the divisor function. To obtain the last bound, fix $\epsilon>0$. At primes $\ell\ge2^{1/\epsilon}$
one has $j+1\le2^j\le\ell^{\epsilon j}$. For each of the remaining primes the maximum of $(j+1)/\ell^{\epsilon j}$
is finite, and the product of these maxima depends only on $\epsilon$. Thus $\tau(n)\ll_\epsilon n^\epsilon$.
Applying the hypothesis to both norms, each dyadic box with $PQ\le X$ contributes $\ll_\epsilon (PQ)^{\alpha/2+\epsilon}$.
There are $O((\log X)^2)$ boxes. Since $\epsilon$ is arbitrary, the result follows.

\end{proof}

For the two bases, an exponent $\alpha<4/3$ in the second moment would improve the power $2/3$ in Corollary~\ref{cor:23}.
The equivalence applies to all rays, including pairs with different local orders.

Theorem~\ref{th:momento} concerns exponents. The next corollary makes both directions explicit at the critical
exponent $(1+\beta)/2$ and expresses them through the number of primes with small reduction order. With
$\alpha=r/(r+1)=1-\beta$ as in Section~\ref{ssec:small}, for $\delta>0$ and $y\ge2$ let
\[
 \mathcal N_\Gamma(y;\delta)=\#\{p\le y:\ p\ \text{of good reduction},\ h_p\le p^{\alpha+\delta}\},
\]
and let $\Qq^{=}_\Gamma(X)\le\Qq_\Gamma(X)$ count the products $pq\le X$ of distinct odd good-reduction primes with
$h_p=h_q$. For $\Gamma=\langle2\rangle$ this is $\mathcal E_2(X)$.

\begin{cor}[Power savings and small reduction orders]\label{co:erdos}
Let $\Gamma$ have rank $r\ge1$.
\begin{enumerate}[label=\textup{(\roman*)},leftmargin=2.4em,itemsep=1pt]
\item If $0<\delta<\beta$, $0<\eta\le1$ and $\mathcal N_\Gamma(y;\delta)\ll y^{1-\eta}$, then
$\Qq_\Gamma(X)\ll X^{(1+\beta-\kappa)/2}$ with $\kappa=\min(\delta,\beta\eta)$.
\item For all $z\ge2$ and $T\ge1$, the number of good-reduction primes $p\in(z,2z]$ with $h_p\le T$ is at most
\[
 2T+2\bigl(T\,\Qq^{=}_\Gamma(4z^2)\bigr)^{1/2}.
\]
Hence, if $\kappa>0$ and $\Qq^{=}_\Gamma(X)\ll X^{(1+\beta)/2-\kappa}$, then $\mathcal N_\Gamma(y;\delta)\ll y^{1-\eta}$
with $\eta=\min(\kappa-\delta/2,\ \beta-\delta)$, for every $0<\delta<\min(2\kappa,\beta)$.
\end{enumerate}
In particular, $\Qq_\Gamma(X)\ll X^{(1+\beta)/2-\kappa}$ for some $\kappa>0$ if and only if the same holds for
$\Qq^{=}_\Gamma$, and if and only if $\mathcal N_\Gamma(y;\delta)\ll y^{1-\eta}$ for some $\delta,\eta>0$.
\end{cor}

\begin{proof}
(i) Let $2\le P\le Q$ be powers of two. As in Section~\ref{sec:superior},
$\B_\Gamma(P,Q)\le\sum_{Q<q\le2Q}2P/h_q$. The primes $q$ with $h_q>Q^{\alpha+\delta}$ contribute at most
$2PQ^{1-\alpha-\delta}=2PQ^{\beta-\delta}$. The remaining primes satisfy $h_q\le Q^{\alpha+\delta}<q^{\alpha+\delta}$, so
there are at most $\mathcal N_\Gamma(2Q;\delta)\ll Q^{1-\eta}$ of them. By Lemma~\ref{le:small}, the number of them with
$D<h_q\le2D$ is $\ll_\Gamma\min(D^{1+1/r},Q^{1-\eta})$. Divide by $D$ and sum over powers of two $D$. The first term
increases geometrically and the second decreases geometrically, and at the transition $D_1^{1+1/r}=Q^{1-\eta}$ we
have $D_1^{1/r}=Q^{\beta(1-\eta)}$. Hence these primes contribute $\ll PQ^{\beta-\beta\eta}$, and
$\B_\Gamma(P,Q)\ll PQ^{\beta-\kappa}$. The summation over dyadic rectangles at the end of
Section~\ref{sec:superior}, with $\beta-\kappa$ in place of $\beta$, gives $\Qq_\Gamma(X)\ll X^{(1+\beta-\kappa)/2}$.

(ii) Let $\mathcal S$ be the set of primes in question, and for $s\ge1$ let $r(s)$ be the number of $p\in\mathcal S$
with $h_p=s$. At most $T$ values of $s$ occur, so Cauchy--Schwarz gives $\sum_sr(s)^2\ge|\mathcal S|^2/T$. The number
of unordered pairs of distinct primes of $\mathcal S$ with equal reduction order is therefore at least
$|\mathcal S|^2/(2T)-|\mathcal S|/2$. Different pairs give different products, each of them is at most $4z^2$, and
each is counted by $\Qq^{=}_\Gamma(4z^2)$ by Lemma~\ref{le:crt}. Thus
$|\mathcal S|^2\le T|\mathcal S|+2T\,\Qq^{=}_\Gamma(4z^2)$, and one of the two terms on the right is at least
$|\mathcal S|^2/2$. This gives the displayed bound. For the second assertion take $T=(2z)^{\alpha+\delta}$. Every
good-reduction prime $p\in(z,2z]$ with $h_p\le p^{\alpha+\delta}$ belongs to $\mathcal S$, and since $\alpha+\beta=1$,
\[
 |\mathcal S|\ll z^{\alpha+\delta}+z^{(\alpha+\delta)/2+(1+\beta)/2-\kappa}=z^{1-(\beta-\delta)}+z^{1-(\kappa-\delta/2)} .
\]
Summing over $z=y/2^j$ with $j\ge1$ gives the bound for $\mathcal N_\Gamma(y;\delta)$. The last assertion of the
corollary follows from (i), from (ii) and from $\Qq^{=}_\Gamma\le\Qq_\Gamma$.
\end{proof}

For $\Gamma=\langle2\rangle$ we have $\mathcal E_2\le\mathrm{Ov}_2$, and Corollary~\ref{cor:ovtotal} gives
$\mathrm{Ov}_2(X)\le\mathcal E_2(X)+O_\varepsilon(X^{2/3+\varepsilon})$. Hence $\Qq_2$, $\mathcal E_2$ and
$\mathrm{Ov}_2$ admit a bound $O(X^{3/4-\kappa})$ with some $\kappa>0$ simultaneously, and they do so exactly when
$\#\{p\le y:\ord_p(2)\le p^{1/2+\delta}\}\ll y^{1-\eta}$ for some $\delta,\eta>0$. Erdős \cite{Erdos76} conjectured
that $\ord_p(2)>p^{c}$ for almost all primes $p$ and every $c<1$. The best known result is that of Erdős and
Murty \cite{EM}, with the exponent $1/2+\varepsilon(p)$ for any positive function $\varepsilon(p)$ tending to $0$. Lowering the exponent $3/4$ in Theorem~\ref{th:main}, in
Theorem~\ref{r15ps:equal} or in Shevelev's bound \cite{Shev08} is thus equivalent to a power-saving form of this
problem at an exponent $1/2+\delta$. In the opposite direction, take $T=(2z)^{c}$ in part~(ii). If
$\mathrm{Ov}_2(X)=o(X^{1/2+\epsilon})$ for every $\epsilon>0$, as conjectured by Pomerance
\cite[Remark~2]{Shev08}, then for every $0<c<1$ and $\epsilon>0$
\[
 \#\{p\le y:\ \ord_p(2)\le p^{c}\}\ll_{c,\epsilon}y^{(1+c)/2+\epsilon},
\]
because $\mathcal E_2\le\mathrm{Ov}_2$. Thus Pomerance's conjecture implies Erdős's conjecture for the base 2, with
a power saving in the number of exceptional primes. We have not found this implication in the literature.

\subsection{From an index moment to the semiprime count}

In rank at least two, a moment bound for the residual index gives a power saving, and the generalized Riemann hypothesis for the Kummer fields provides such a moment bound.

\begin{prop}[Transfer from an index moment]\label{pr:moment}

Let $s$ be fixed with $0<s<1$. Assume that $\sum_{Q<q\le2Q}I_q^s\ll_{\Gamma,s} Q/\log Q$ for every sufficiently large $Q$. With
$b_s=(1-s)/(r+1-rs)$ we have
\[
 \A_\Gamma(z)\ll_{\Gamma,s} z^{1+b_s},\qquad
 \Qq_\Gamma(X)\ll_{\Gamma,s} X^{(1+b_s)/2}.
\]

\end{prop}

\begin{proof}

In the block $D<h_q\le2D$ the residual index is $\gg Q/D$. The number of $q$ in the block is
\[
 \ll_{\Gamma,s} L_Q^{-1}\min(D^{1+1/r},\,Q^{1-s}D^s).
\]
Divide by $D$ and sum dyadically. The first contribution increases geometrically and the second decreases
geometrically on the two sides of $D_0=Q^{(1-s)/(1+1/r-s)}$. At the transition point $D_0^{1/r}=Q^{b_s}$.
Hence $\sum_q 1/h_q\ll Q^{b_s}/L_Q$. The bound $2P/h_q$ and the geometric sum over rectangles in Section~\ref{ssec:rect}, with $b_s$ in place of $\beta$, even save a logarithm in the boxes that contribute
at the leading scale.

\end{proof}

\begin{cor}[Residual-index moment under GRH]\label{co:grh}

Let $r\ge2$. Assume GRH for the Dedekind zeta functions of $K_m=\Q(\zeta_m,\Gamma^{1/m})$, $m\ge1$, where
one adjoins $m$th roots of the generators. If $0<s<(r-1)/(r+2)$, the hypothesis of Proposition~\ref{pr:moment}
holds. Consequently, for every $\epsilon>0$,
\begin{equation}
 \A_\Gamma(z)\ll_{\Gamma,\epsilon} z^{1+3/(4r+2)+\epsilon},\qquad
 \Qq_\Gamma(X)\ll_{\Gamma,\epsilon} X^{1/2+3/(8r+4)+\epsilon}.
 \label{eq:grh}
\end{equation}
In particular, for $\Gamma=\langle2,3\rangle$ the second exponent is $13/20+\epsilon$.

\end{cor}

The same range for weighted residual-index moments is treated over number fields in \cite{PinedaRI}.
The proof below establishes the rational case used in the transfer to the semiprime count.

\begin{proof}

For a good-reduction prime $q$, the
condition $m\mid I_q$ means that $q\equiv1\bmod m$ and that every generator is an $m$th power modulo $q$.
This follows by working in the cyclic group of order $q-1$. In that case $q\nmid m$, and the condition is
equivalent to complete splitting of $q$ in $K_m$, since all the roots of $X^m-g$ reduce to the residue field
for each generator $g$. In the converse direction we restrict to unramified primes dividing none of the
chosen numerators or denominators.

Write $\operatorname{disc}(K)$ for the discriminant of a number field $K$. Put $n_m=[K_m:\Q]$. We first establish the uniform estimates
\[
 n_m\gg_\Gamma \varphi(m)m^r,\qquad
 \frac{\log|\operatorname{disc}(K_m)|}{n_m}\ll_\Gamma \log(2m).
\]
For the first inequality choose $r$ positive multiplicatively independent elements $g_i$ of $\Gamma$, taking
squares if necessary. The integer matrix $(v_\ell(g_i))_{\ell,i}$ has a minor of size $r$ with determinant
$D\ne0$. Put $F=\Q(\zeta_m)$. If $x=\prod g_i^{t_i}>0$ is an $m$th power in $F$, its positive real root
$x^{1/m}$ belongs to $F$. Set $d=\gcd(m,(v_\ell(x))_\ell)$ and $e=m/d$. If $x\ne1$, there is a rational
$u>0$ with $x=u^d$ and $\gcd(e,(v_\ell(u))_\ell)=1$. The degree $h$ of $u^{1/e}$ is exactly $e$. Indeed,
taking norms in $(u^{1/e})^e=u$ gives $e\mid h v_\ell(u)$ for every $\ell$, hence $e\mid h$, and the polynomial
$T^e-u$ gives $h\le e$. If $e>2$, the conjugate roots include a nonreal one. The field generated by the
positive root is real and fails to be normal, whereas every subfield of $F/\Q$ is normal. Thus $e\le2$ and
$m\mid 2v_\ell(x)$ for every $\ell$. This divisibility holds directly for $x=1$. Multiplying by the adjugate
of the minor yields $m\mid2Dt_i$ for every $i$. Hence the kernel of $(\Z/m)^r\to F^\times/F^{\times m}$
has size at most $(2|D|)^r$. Kummer theory gives the lower bound on the degree.

For the second inequality, ramification is confined to primes dividing $m$ or in the fixed support of $\Gamma$.
In a Galois extension of degree $n_m$ with ramification index $e_\ell$, the discriminant exponent at $\ell$
is at most $n_m(1-1/e_\ell+v_\ell(e_\ell))$ by the local bound for the different \cite[Chapter~III]{SerreLocal}. Since
$e_\ell\mid n_m$, the logarithmically weighted sum of the last terms, after division by the degree, is at
most $\log n_m$. Adjoining one root for each of $k$ generators gives $n_m\le \varphi(m)m^k$. The remaining
terms contribute at most $\log m+O_\Gamma(1)$, as claimed.

For a number field $K$, let $\Ocal_K$ be its ring of integers and put $\psi_K(x)=\sum_{\mathfrak p^j:\,\Nrm(\mathfrak p)^j\le x}\log\Nrm(\mathfrak p)$, summed over nonzero prime ideals $\mathfrak p$ of $\Ocal_K$ and integers $j\ge1$, where $\Nrm(\mathfrak p)=|\Ocal_K/\mathfrak p|$. Under GRH the prime ideal theorem \cite{GMo} bounds the absolute value of $\psi_{K_m}(x)-x$ by $O(\sqrt x\log(2x)(\log |\operatorname{disc}K_m|+n_m\log(2x)))$,
with an absolute implied constant. Every completely split prime $q$ in $(Q,2Q]$ contributes at least $n_m\log Q$
to $\psi_{K_m}(2Q)$. In particular,
\[
 \#\{Q<q\le2Q:m\mid I_q\}
 \ll_\Gamma \frac{Q}{n_m\log Q}+\sqrt Q\log(2mQ).
\]
This upper bound is all that the argument uses.

Choose $\eta$ such that
\[
 \frac1{r+1-rs}<\eta<\frac1{2(s+1)},\qquad Y=Q^\eta .
\]
This interval is nonempty precisely in the stated range. For $I_q\le Y$ group the primes according to $m=I_q$
and bound the equality condition by divisibility. Since $\sum_{m\ge1}m^s/(\varphi(m)m^r)<\infty$ for $r\ge2$
and $s<1$, the contribution to the moment is
\[
 \ll_{\Gamma,s} \frac{Q}{\log Q}+\sqrt Q\,Y^{s+1}\log(2QY)
 \ll_{\Gamma,s} Q/\log Q.
\]
For $T<I_q\le2T$ with $T\ge Y$ we have $h_q\le 2Q/T$. Lemma~\ref{le:small} bounds this block by $O_{\Gamma,s}(Q^{1+1/r}T^{s-1-1/r}/\log Q)$.
Summing the blocks from $Y$ onwards gives
\[
 \sum_{\substack{Q<q\le2Q\\I_q>Y}} I_q^s
 \ll_{\Gamma,s} \frac{Q^{1+1/r}Y^{s-1-1/r}}{\log Q}.
\]
With $Y=Q^\eta$, the inequality $\eta>1/(r+1-rs)$ makes the exponent of $Q$ strictly less than one. Finally,
$b_s=(1-s)/(r+1-rs)$ decreases to $3/(4r+2)$ as $s$ increases to $(r-1)/(r+2)$. Proposition~\ref{pr:moment}
now gives \eqref{eq:grh}.

\end{proof}

\section{Order fibers of base 2}\label{sec:fibras}

For one base Corollary~\ref{co:grh} gives nothing, and Corollary~\ref{co:erdos} shows that a power saving requires new information on the orders of $2$. This section assumes the strongest available form of such information, a bound for the number of primes of each exact order that follows from a conjecture of Kovač and Luca. It gives an upper bound of the order $\sqrt X$ of Galway's conjecture, up to a power of $\log X$ (Theorem~\ref{th:KLsharp}), and Remark~\ref{re:KLartin} shows how strong the hypothesis is.

In base 2 put $\nu(t)=\#\{q\in\Pp:q>2,\ \ord_q(2)=t\}$. If $t\mid p-1$, the primes in this fiber are potential
companions of $p$. A single value of $\nu(t)$ controls the total number of these candidates even when $q$
is much larger than $p$.

\subsection{A sieve estimate uniform in the coefficients}

The proofs in this section, in Remark~\ref{re:trunc} and in Section~\ref{sec:rsa} count pairs of primes of the form $ut+1$, $vt+1$, and they need an upper bound uniform in the coefficients $u,v$.

\begin{lemma}\label{le:sieve}

For integers $1\le u<v$ put $\Delta=uv(v-u)$ and $G(n)=(n/\varphi(n))^2$. If $T\ge2$, then, with an absolute
implied constant,
\[
 \#\{1\le t\le T:ut+1,\ vt+1\ \text{primes}\}
 \ll \frac{T}{\log^2(2T)}\,G(\Delta).
\]
Moreover, $\sum_{u\sim U,\ u<v\le2Ru}G(uv(v-u))\ll RU^2$ for $R,U\ge1$, where $u\sim U$ means $U\le u<2U$,
and $\sum_{1\le u<v\le Y}G(uv(v-u))\ll Y^2$ for $Y\ge1$.

\end{lemma}

\begin{proof}

Let $\rho(\ell)$ be the number of roots of $(ut+1)(vt+1)$ in $\F_\ell$. It is at most two and is always
less than $\ell$ since $t=0$ avoids all the roots. It equals two except at primes dividing $\Delta$. For
squarefree $d$, extend $\rho$ multiplicatively. The exact residue-class count has main term $T\rho(d)/d$
and error at most $\rho(d)$.

Apply the Selberg upper-bound sieve \cite[Chapter~7]{FI} with $D=T^{1/10}$, $w=D^c$ and $P(w)=\prod_{\ell\le w}\ell$,
taking $c>0$ to be a sufficiently small absolute constant. The sifted count up to $T$ is bounded by
\[
 \frac{T}{\mathcal J}+
 O\!\left(\sum_{d\le D^2}6^{\omega(d)}\right),\qquad
 \mathcal J=\sum_{\substack{d\mid P(w)\\d\le D}}
       \prod_{\ell\mid d}\frac{\rho(\ell)}{\ell-\rho(\ell)}.
\]
Indeed, upon squaring the Selberg weights, supported up to $D$ and of absolute value $\le1$, each least
common multiple $d$ occurs at most $3^{\omega(d)}$ times, with remainder $\le2^{\omega(d)}$. The error is
$O(D^2\log^5(2D))$ by comparison with the sixfold divisor function. Normalize the full Euler product associated
with $\mathcal J$ to a probability measure on divisors of $P(w)$. Each prime $\ell$ is included independently with
probability $\rho(\ell)/\ell$. The mean of $\log d$ is at most $2\sum_{\ell\le w}(\log\ell)/\ell\ll\log(2w)$.
For sufficiently small $c$ and large $T$, Markov's inequality leaves mass $\ge1/2$ on $d\le D$. Hence
\[
 \mathcal J^{-1}\le2\prod_{\ell\le w}(1-\rho(\ell)/\ell)
 \ll \frac{G(\Delta)}{\log^2(2w)} .
\]
In the last bound compare each factor with $(1-1/\ell)^2$. At the good primes the ratio is $\le1$, and at
the divisors of $\Delta$ it is $\le(1-1/\ell)^{-2}$. If one of the two prime values equals a sieving prime,
then $t\le w$. These exceptions and the remainder are absorbed in the stated bound. Increasing the constant
handles bounded $T$.

For the averages, use $G(uv(v-u))\le G(u)G(v)G(v-u)$. Since $\sum_{n\le T}G(n)^3\ll T$, Hölder separates
the three factors. To prove the mean estimate expand $(n/\varphi(n))^6$ over squarefree divisors. Its coefficient
at a prime is $(1-1/\ell)^{-6}-1=O(1/\ell)$, so a convergent Euler product bounds the mean. In the first
rectangle each $u$ occurs $O(RU)$ times, and each $v$ or difference $v-u$ occurs $O(U)$ times in $(0,4RU]$.
The three sums of cubes are $O(RU^2)$. In the square of side $Y$ each integer or difference occurs at most
$Y$ times. The same inequality gives $O(Y^2)$.

\end{proof}

\subsection{A cyclotomic hypothesis and its logarithmic exponent}

We now combine the hypothesis on $\nu(t)$ with this sieve.

\begin{prop}\label{pr:KL}

If $\nu(t)\ll\log(2t)$ for every $t$, then
\[
 \A_2(z)\ll z\log z,\qquad \Qq_2(X)\ll\sqrt X\log X.
\]
Conjecture 2 of Kovač and Luca \cite{KL} implies this hypothesis.

\end{prop}

\begin{proof}

Let $\Phi_t$ be the $t$th cyclotomic polynomial. Every prime of order $t$ divides $\Phi_t(2)$, so $\nu(t)\le\omega(\Phi_t(2))$.
Here $\omega$ counts distinct prime factors. A prime dividing the cyclotomic index can contribute to $\omega$
without having order $t$. The conjecture asserts that $\omega(\Phi_t(2))\le10\log t$ apart from finitely
many $t$, which can be absorbed into the implied constant. For each prime $p$, Lemma \ref{le:crt} gives
\[
 \#\{q>p:pq\text{ is a base-2 pseudoprime}\}
 \le \sum_{t\mid p-1}\nu(t)\ll \tau(p-1)\log(2p).
\]
We use $\sum_{p\le y}\tau(p-1)\ll y$. To see this, pair the divisors of $p-1$ on either side of its square
root and apply the Brun--Titchmarsh inequality to obtain
\[
 \sum_{p\le y}\tau(p-1)
 \le 2\sum_{d\le\sqrt y}\pi(y;d,1)
 \ll \frac y{\log y}\sum_{d\le\sqrt y}\frac1{\varphi(d)}
 \ll y .
\]
The last sum is $O(\log y)$ by partial summation from the mean value of $d/\varphi(d)$. Summing over $p\le2z$
and $p\le\sqrt X$, respectively, gives the two conclusions.

\end{proof}

\begin{theorem}\label{th:KLsharp}

Under the hypothesis on $\nu(t)$ in Proposition \ref{pr:KL}, put
\[
 \lambda=\min_{1\le s\le2}\frac{2^s-1}{2s-1}
          =0.9105498551\ldots .
\]
Then $\A_2(z)\ll z L_z^\lambda$ and $\Qq_2(X)\ll \sqrt X L_X^\lambda$.

\end{theorem}

\begin{proof}

Fix $1<s<2$, $\gamma=1-1/s$ and $b=1+\gamma+(2^s-2)/s$. We need the divisor moment
\begin{equation}
 \sum_{P<p\le2P}\tau(p-1)^s\ll_s P(\log P)^{2^s-2}.
 \label{eq:pollack}
\end{equation}
It follows directly from Pollack's multiplicative sieve \cite[Theorem 1.1]{Pollack} applied to $f(m)=\tau(m)^s$,
which satisfies $f(\ell^j)=(j+1)^s\le2^{sj}$ and $f(m)\ll_{s,\epsilon}m^\epsilon$. If $m=p-1$ and $P<p\le2P$,
then $m$ avoids $-1\bmod\ell$ for each $\ell\le\sqrt P$. The exponent in the exponential in Pollack's bound
is $\sum_{\ell\le2P}f(\ell)/\ell-\sum_{\ell\le\sqrt P}1/\ell
=(2^s-1)\log\log P+O_s(1)$. Multiplication
by the leading factor $P/\log P$ gives \eqref{eq:pollack}.

Write $L=L_X$ and put $R=X/P^2$ on a dyadic block for the smaller prime $P<p\le2P$, where $P\le\sqrt X$. If
$P<\sqrt X/L^{(b+2)/2}$, summing over all such blocks gives $O_s(\sqrt X L^{-b/2})$ by the elementary bound
in Proposition \ref{pr:KL}. In the remaining blocks $R\ge1$, $R\le L^{b+2}$ and $\log P\asymp L$. For $Y\ge1$,
if $I_p>Y$, the hypothesis on $\nu$ bounds the number of such primes by $O(PL/Y)$, on summing $\nu(t)$ over
$t\le 2P/Y$. By the individual estimate in Proposition \ref{pr:KL}, \eqref{eq:pollack} and Hölder's inequality,
the number of their partners is
\begin{equation}
 \ll_s L\,(P L^{2^s-2})^{1/s}(PL/Y)^\gamma
 =P L^{b}Y^{-\gamma}.
 \label{eq:KLbad}
\end{equation}
If $I_p\le Y$, the representation $p=gu+1,\ q=gv+1$ satisfies $u\le Y$. Since $q\le X/P$ and $p>P$, we have
$u<v\le 2Ru$ for large $P$. On a block $U\le u<2U$, we have $g\le2P/U$. When $Y$ is of polylogarithmic size
in $X$, the logarithm of this bound is comparable to $L$. Lemma \ref{le:sieve} bounds the number for each
pair $u,v$ by
\[
 \ll\frac{P}{U L^2}\,G(uv(v-u)),
\]
Averaging by Lemma \ref{le:sieve} gives $\sum G(uv(v-u))\ll RU^2$ in this rectangle. Summing over $U$ gives
$O(PRY/L^2)$ for $I_p\le Y$.

Choose $Y=(L^{b+2}/R)^{1/(1+\gamma)}\ge1$. Together the two estimates give, on each block,
\[
 \ll_s P R^{\gamma/(1+\gamma)} L^{(b-2\gamma)/(1+\gamma)}
 =\sqrt X\left(\frac{P}{\sqrt X}\right)^{(1-\gamma)/(1+\gamma)}
   L^{(2^s-1)/(2s-1)} .
\]
The sum over $P$ is geometric. The objective takes the value $1$ at both endpoints and is less than $1$
at $s=3/2$, so its minimum lies in the interior. Its derivative vanishes when $2^s((2s-1)\log2-2)+2=0$,
and the derivative of the left-hand side is $2^s(2s-1)(\log2)^2>0$. Fixing this value of $s$ proves the
global estimate. The dyadic bound follows from $\A_2(z)\le\Qq_2(4z^2)$.

\end{proof}

\begin{remark}[Strength of the exact-order hypothesis]\label{re:KLartin}
The hypothesis $\nu(t)\ll\log(2t)$ of Proposition~\ref{pr:KL} and Theorem~\ref{th:KLsharp} is strong. Every odd prime
$p\le x$ with $\ord_p(2)=t$ is counted by $\nu(t)$, and $t<x$, so the hypothesis gives
\[
 \sum_{2<p\le x}\frac{1}{\ord_p(2)}\le\sum_{t<x}\frac{\nu(t)}{t}\ll(\log x)^2 .
\]
Murty and Srinivasan \cite[Theorem~4]{MurtySrinivasan} showed that a bound $O(x^{1/4})$ for this sum already implies
that $2$ is a primitive root modulo infinitely many primes. Hence Conjecture~2 of Kovač and Luca \cite{KL} implies
Artin's conjecture for the base $2$ in its qualitative form. This implication does not need the full strength of the
hypothesis. Erdős conjectured the weaker bound $\nu(t)\ll_\varepsilon t^{\varepsilon}$ for every $\varepsilon>0$, and
Felix and Murty \cite{FM} studied what it implies for Artin's conjecture. That bound already gives
$\sum_{2<p\le x}1/\ord_p(2)\ll_\varepsilon x^{\varepsilon}$, so the criterion of Murty and Srinivasan applies under it
as well. The same count gives
$\#\{p\le y:\ \ord_p(2)\le y^{c}\}\le\sum_{t\le y^{c}}\nu(t)\ll y^{c}\log y$, so the hypothesis also implies that
$\ord_p(2)>p^{c}$ for almost all primes $p$, for every fixed $c<1$, far beyond the exponent $\frac12$ of Erdős and
Murty \cite{EM}. What Theorem~\ref{th:KLsharp} uses beyond the conjecture of Erdős is the logarithmic size of
$\nu(t)$, which gives the exponent $\lambda=0.9105498551\ldots$ in the bound for the semiprime count.
\end{remark}

\subsection{Remarks on the order fibers}

We end this section by relating the order fibers to the equal-order products of Section~\ref{sec:equal} and to the lower bound of Section~\ref{sec:inferior}.

The hypothesis in Theorem \ref{th:KLsharp} controls $\nu(t)$ for each $t$ without restricting the index; fibers with bounded indices are counted in Remark~\ref{re:trunc}.
Near the critical rank-one cutoff, the order and the index can both be of size $Q^{1/2+o(1)}$. These are counts of pairs with common order, whereas a distribution of the indices at individual primes controls single-prime counts.
Proposition \ref{pr:rec} also requires the reciprocal-order weight at this scale. For Erdős problem 893,
Kovač and Luca \cite{KL} prove $\limsup f(2n)/f(n)=\infty$ with $f(n)=\sum_{k\le n}\tau(2^k-1)$. The
construction of two individual primes of the same order $4u$ yields actual products and provides the input
for the lower bound in Section~\ref{sec:inferior}.

\section{Products of primes of equal order}\label{sec:equal}

This section proves Theorem~\ref{intro:over} for equal-order products and overpseudoprimes.

Shevelev \cite[Theorem~3 and Remark~1]{Shev08} proved that the number $\mathrm{Ov}_2(x)$ of overpseudoprimes to
base 2 up to $x$ satisfies $\mathrm{Ov}_2(x)\le x^{3/4}(1+O((\log\log x)^2/\log x))$. In his argument the
overpseudoprimes with exactly $k\ge3$ prime factors contribute $O_k(x^{(k+1)/(2k)}(\log\log x)^{2k})$, and in
\cite[Remark~2]{Shev08} he records Pomerance's conjecture that $\mathrm{Ov}_2(x)=o(x^{1/2+\epsilon})$ for every
$\epsilon>0$. Theorem~\ref{r15ps:equal} below replaces the factor $(\log\log x)^{2k}$ by a saving of more than
$(\log x)^{k}$ for every fixed $k$, and Corollary~\ref{cor:ovtotal} saves more than $(\log x)^2$ over Shevelev's
bound for the total count.

For $k=2$ the power $3/4$ and most of this saving already follow from Theorem~\ref{th:main}: every product of two
distinct odd primes of equal order passes base 2 by Lemma~\ref{le:crt}, so that $\mathcal E_2(Z)\le\Qq_2(Z)$. The
logarithmic factor in Theorem~\ref{r15ps:equal} is smaller than that of Theorem~\ref{th:main} by
$L_Z^{\dd/2}H_Z^{3/4}$, up to a constant. By Theorem~\ref{th:momento}, the equal-order counts and
$\Qq_2$ have the same power growth.

\subsection{\texorpdfstring{Products of $k$ primes of equal order}{Products of k primes of equal order}}

We first bound the products of $k$ distinct primes at which $2$ has a common order. For squarefree odd composites with $k$ prime factors, this characterizes overpseudoprimality to base 2.

Fix \(k\geq2\). If \(n\) is counted by \(\mathcal E_k(Z)\) and \(m\) is the common order of \(2\) modulo its
prime factors, then by the Chinese remainder theorem all divisors \(d>1\) of \(n\) have order \(m\). In particular \(2^{n-1}=1\bmod n\), since \(n\equiv1\bmod m\).
In Shevelev's terminology these integers are overpseudoprimes. The next bound even allows repeated factors
if one requires equality of orders only at the underlying primes. In that case the prime powers dividing
\(n\) may impose an additional condition for the overpseudoprime property; see \cite[Theorems 2 and 9]{Shev12}.

\begin{theorem}[Theorem~\ref{intro:over}, first bound]\label{r15ps:equal}

For sufficiently large \(Z\) in terms of \(k\),
\[
 \mathcal E_k(Z)\ll_k
 \frac{Z^{(k+1)/(2k)}}
 {L_Z^{k+\dd}H_Z^{3/2}} .
\]
This bound holds globally for all relative sizes of the factors.

\end{theorem}

\begin{proof}

Order the primes \(p_1\leq\cdots\leq p_k\) and place them in bands \(P_i<p_i\leq2P_i\), with \(P_i=2^{j_i}\)
and \(P_1\leq\cdots\leq P_k\). We may allow repetitions for an upper bound. Write \(S=\prod_iP_i\leq Z\).
In a band, for fixed \(m\) the number of possible factors \(p_i\) of order \(m\) is
\(O(\min(m,P_i/m))\). The first estimate counts prime divisors of \(2^m-1\) and the second counts integers
\(1+am\leq2P_i\), \(a\geq1\). Once \(p_1\) is given, its order determines \(m\). Using \(\min(m,P_i/m)\leq\sqrt{P_i}\)
for the other factors gives \(O_k(P_1\prod_{i>1}\sqrt{P_i})=O_k(\sqrt{S P_1})\). If \(P_1<Z^{7/(8k)}\),
this is \(O_k(Z^{(k+1)/(2k)-1/(16k)})\). The number of boxes is \(O_k((1+L_Z)^k)\) and is absorbed by the
power saving.

Two ranges of \(m\) can also be handled directly, irrespective of the factor sizes. If \(m\leq Z^{1/(4k)}\),
summing \(O(m^k)\) tuples per modulus gives \(O_k(Z^{(k+1)/(4k)})\). For \(m\geq Z^{3/(4k)}\), writing \(p_i=1+ma_i\)
forces \(\prod_i a_i\leq Z/m^k\). The number of ordered \(k\)-tuples of positive integers with product \(\leq T\),
\(T\geq1\), is
\[
 \leq T\left(\sum_{a\leq T}1/a\right)^{k-1}
 \leq T(1+\log T)^{k-1},
\]
by fixing the first \(k-1\) factors and bounding the last. There are no tuples for \(T<1\). Thus this
range contributes
\[
 O_k\!\left(Z(1+L_Z)^{k-1}
       \sum_{m\geq Z^{3/(4k)}}m^{-k}\right)
 = O_k\!\left(Z^{(k+3)/(4k)}(1+L_Z)^{k-1}\right).
\]
Both exponents are less than \((k+1)/(2k)\).

Now fix a box with \(P_1\geq Z^{7/(8k)}\) and a dyadic interval \(D<m\leq2D\) in \(Z^{1/(4k)}<m<Z^{3/(4k)}\),
truncated at the endpoints. For every index we have \(\log P_i\asymp_k L_Z\), and \(D\) lies between a fixed
positive power of \(P_1\) and \(P_1^{6/7}\), up to cutoff constants. Each prime \(p_1\) of the band with order in this interval determines its exact order \(m\),
and \(m\mid p_1-1\). For any \(i>1\) and this value of \(m\), the number of possible \(p_i\) is
\begin{equation}
 \ll_k L_Z^{-1}\min\{\varphi(m),P_i/\varphi(m)\}
 \leq C_k\sqrt{P_i}/L_Z . \label{r15ps:level}
\end{equation}
The first bound in the minimum uses \(\Phi_m(2)\) with $\log\Phi_m(2)=\varphi(m)\log 2+O(1)$. By Möbius
inversion the error is at most $\sum_{d\ge1}|\log(1-2^{-d})|$. The second bound is Brun--Titchmarsh for
\(p_i\equiv1\bmod m\) up to \(2P_i\), since \(\log(2P_i/m)\gg_k L_Z\). For \(i=3,\ldots,k\) use the last
bound in \eqref{r15ps:level}. For \(i=2\) we keep the dependence on \(p_1\). Since \(\varphi(m)\le2D\) and
\(1/\varphi(m)=(m/\varphi(m))/m\) with \(m/\varphi(m)\le(p_1-1)/\varphi(p_1-1)\), the minimum is at most
\[
 \frac{p_1-1}{\varphi(p_1-1)}\min\{2D,P_2/D\}.
\]
By Lemma~\ref{le:ford} with exponent \(1\) and \(Q=P_1\), the sum of \((p_1-1)/\varphi(p_1-1)\) over the primes
\(p_1\) of the band for which \(p_1-1\) has a divisor in \((D,2D]\) is \(O_k(P_1/(L_Z^{1+\dd}H_Z^{3/2}))\). The
interval \((D,2D]\) therefore contributes
\[
 \ll_k \frac{P_1\min\{D,P_2/D\}\prod_{i>2}\sqrt{P_i}}{L_Z^{k+\dd}H_Z^{3/2}} .
\]
Up to \(\sqrt{P_2}\) the sum over dyadic \(D\) is a geometric progression in \(D\), and beyond it one in
\(1/D\). Both are \(O(\sqrt{P_2})\), so the total contribution of the box is
\[
 \ll_k
   \frac{P_1\prod_{i>1}\sqrt{P_i}}{L_Z^{k+\dd}H_Z^{3/2}}.
\]
For comparable factors \(P_i\asymp P\), its algebraic part has size \(P^{(k+1)/2}\). To count products of
arbitrary shape we must sum the boxes without losing powers of \(L_Z\).

Put \(a_i=j_i-j_1\) for \(i\geq2\), and put \(A=\sum_{i=2}^k a_i\). Then \(a_i\geq0,\ k j_1+A\leq(L_Z/\log 2)\)
and
\[
 P_1\prod_{i>1}\sqrt{P_i}=2^{(k+1)j_1/2+A/2}.
\]
For fixed \(a_2,\ldots,a_k\), summing first over \(j_1\) gives
\[
 \ll_k 2^{(k+1)((L_Z/\log 2)-A)/(2k)+A/2}
       =Z^{(k+1)/(2k)}\,2^{-A/(2k)}.
\]
The sum over the \(a_i\geq0\) is bounded, even without ordering them, by a product of \(k-1\) convergent
geometric series. This completes the passage to the global bound.

\end{proof}

\subsection{Overpseudoprimes}

Products of two distinct primes give the largest exponent in our upper bound for overpseudoprimes. Theorem~\ref{r15ps:equal} with $k=2$ bounds them, and a cruder argument shows that those with at least three prime factors, counted with multiplicity, number $O(Z^{2/3+\varepsilon})$.

\begin{cor}[Theorem~\ref{intro:over}, second bound]\label{cor:ovtotal}

With $L=L_Z$ and $H=H_Z$, one has
\[
 \mathrm{Ov}_2(Z)\ll \frac{Z^{3/4}}{L^{2+\dd}H^{3/2}},\qquad
 \#\{n\text{ counted}:\Omega(n)\ge3\}\ll_\varepsilon Z^{2/3+\varepsilon}.
\]
\end{cor}

\begin{proof}

For two distinct primes apply Theorem~\ref{r15ps:equal} with $k=2$. There are at most $\sqrt Z$ squares.
If $k=\Omega(n)\ge3$ and $m=\ord_n(2)$, then $n\mid2^m-1$ and every prime factor of $n$ is $1\bmod m$. It
follows that $k\log(m+1)\le\min(\log Z,m\log2)$. For $m\le L$ use that $m/\log(1+m)$ is increasing, and
for $m>L$ use the first term. This gives $k\ll L/H$ uniformly. Fix $k$ and put $T=Z^{1/(2k)}$. If $m\le T$,
its factors come from at most $m$ prime divisors of $2^m-1$, and the lists with repetitions contribute $\sum_{m\le T}m^k\le T^{k+1}$.
If $m>T$, each $m$ contributes at most $Z(\sum_{p\le Z,\ p\equiv1\bmod m}1/p)^k$. Up to $2m$ this sum is
$O(1/m)$, and from $2m$ to $Z$ partial summation with Brun--Titchmarsh gives $O((1+H)/\varphi(m))$. The
bound $m/\varphi(m)\ll1+H$ gives $O((1+H)^2/m)$. For $k\ge3$, we have $\sum_{m>T}m^{-k}\le2 T^{1-k}$ by
comparison with the integral starting at the first integer greater than $T$. Since $[C(1+H)^2]^k=Z^{o(1)}$
uniformly for $k\ll L/H$, the two ranges contribute $Z^{(k+1)/(2k)+o(1)}\le Z^{2/3+o(1)}$. Summation over
$k$ preserves the $o(1)$. This proves both assertions.

\end{proof}

\begin{remark}[Truncated order fibers]\label{re:trunc}
For $t\ge1$ and $Y\ge1$ write $\nu_Y(t)=\#\{p>2:\ord_p(2)=t,\ I_p\le Y\}$. Fixing both index values allows a
two-form sieve: for $M\ge2$ the number of distinct products of primes with common order $t\in(M,2M]$ and indices
at most $Y$ satisfies
\[
 \sum_{M<t\le2M}\binom{\nu_Y(t)}2\ll M\Bigl(\frac{Y}{\log(2M)}\Bigr)^2 ,
\]
and each such product is a strong pseudoprime to base 2. Indeed, a product in the count has $p=ut+1$ and $q=vt+1$
with $1\le u<v\le Y$. With both slopes fixed, Lemma~\ref{le:sieve} with $T=2M$ bounds the number of $t$ even after
the exact-order conditions are dropped, and summing the second mean-value estimate of that lemma gives the bound.
Lemma~\ref{le:crt} gives the Fermat condition and Lemma~\ref{le:fuerte} the strong condition. Without the
restriction on the indices the fibers are controlled only through hypotheses such as the one in
Theorem~\ref{th:KLsharp}.
\end{remark}

\section{A base-2 lower bound with an explicit constant}\label{sec:inferior}

The upper bounds of Sections~\ref{sec:superior}--\ref{sec:equal} lie far above the conjectured order $\sqrt X/L_X^2$. In the other direction, unconditional constructions reach only the logarithmic scale, at which Pomerance \cite{Pom82} asked for unboundedness. This section gives explicit constants at that scale and separates the strong pseudoprimes from the products that fail the strong test. Section~\ref{sec:constantes} gives conditional lower bounds of the conjectured order.

Write $\Qq_2^{\,s}(X)$ and $\Qq_2^{\,w}(X)$ for the two parts of $\Qq_2(X)$ consisting respectively of strong
pseudoprimes and of Fermat pseudoprimes that fail the strong test.

\begin{theorem}\label{th:lower}

Let $C_\varphi=\zeta(2)\zeta(3)/\zeta(6)$, where $\zeta$ is the Riemann zeta function. Then
\[
 \liminf_{X\to\infty}\frac{\Qq_2^{\,s}(X)}{\log X}\ge
 \frac{C_\varphi}{6\log 2},\qquad
 \liminf_{X\to\infty}\frac{\Qq_2^{\,w}(X)}{\log X}\ge
 \frac{C_\varphi}{6\log 2}.
\]
In particular, $\liminf\Qq_2(X)/\log X\ge
C_\varphi/(3\log2)=0.9346723136\ldots$.

\end{theorem}

\begin{proof}

We apply the strong criterion of Lemma \ref{le:fuerte} to each of the two constructions.

Fix an odd integer $u\ge7$. The Bang--Zsigmondy theorem and Lehmer's construction \cite{Lehmer,E49,Somer}
give primes of orders $u$ and $2u$. The first prime is odd, and $u$ is odd, so $2u$ divides its predecessor.
The predecessor of the second prime is a multiple of $2u$. Their product satisfies Lemma \ref{le:crt}, and
the orders have valuations $0$ and $1$. The product is at most $\Phi_u(2)\Phi_{2u}(2)$. Indeed, if a prime
has order $t$, then the factorization $2^t-1=\prod_{d\mid t}\Phi_d(2)$ forces it to divide $\Phi_t(2)$ since
each factor indexed by $d<t$ divides $2^d-1$.

By Schinzel's theorem \cite{Schinzel}, in the explicit form given in \cite[Theorem 2.3]{Somer}, there are
two distinct primes of order $4u$ for each odd $u\notin\{1,3,5\}$. Their product divides $\Phi_{4u}(2)$
and passes the strong test because the two orders agree. The exception $u=5$ matters in applying Proposition
3.2 of \cite{LRU}. For $m=20$ one has $\Phi_{20}(2)=205=5\cdot41$, with orders $4$ and $20$. The assertion
used there holds for $m>20$ with $m\equiv4\bmod8$. Changing this endpoint preserves the asymptotic conclusion.

Uniformly in $t\ge1$, Möbius inversion gives
\[
 \left|\log\Phi_t(2)-\varphi(t)\log 2\right|
 \le \sum_{j\ge1}\big|\log(1-2^{-j})\big|=C_0<\infty .
\]
Since $\varphi(2u)=\varphi(u)$ and $\varphi(4u)=2\varphi(u)$, the logarithm of each product constructed
above is $\le 2\varphi(u)\log2+2C_0$. Choose one pair for each $u$ in each family. Unique factorization
recovers the multiset of orders and hence $u$. The strong criterion distinguishes the two families, so their
counts add without multiplicity.

Let $T(y)=\#\{m:\varphi(m)\le y\}$, and let $U(y)$ denote the part counting odd $m$. The classical inverse
distribution theorem for $\varphi$, as stated in \cite{Bateman}, gives $T(y)\sim C_\varphi y$. Separating
the exact power of two in $m$ gives
\[
 T(y)=2U(y)+\sum_{j\ge1}U(y/2^j),\qquad
 2U(y)-U(y/2)=T(y)-T(y/2).
\]
For $b(y)=U(y)/y$, this identity reads $b(y)=C_\varphi/4+o(1)+b(y/2)/4$. Since $b$ is bounded, its upper
and lower limits are both $C_\varphi/3$. Applying this with $y=(\log X-2C_0)/(2\log2)$ and omitting the
three exceptional indices proves the two lower limits.

\end{proof}

The proof in \cite[Proposition 3.2]{LRU} allows any constant $c<1/(8\log2)=0.18033688\ldots$. Theorem \ref{th:lower}
improves the coefficient while preserving the logarithmic order of the lower bound. Neither construction
controls the ratio $q/p$. The two primes of the strong family have the same order, so the first lower limit also
holds for the equal-order products and for the overpseudoprimes:
\[
 \liminf_{X\to\infty}\frac{\mathcal E_2(X)}{\log X}\ge\frac{C_\varphi}{6\log2}=0.4673361568\ldots,
 \qquad \mathrm{Ov}_2(X)\ge\mathcal E_2(X),
\]
in the notation of Section~\ref{sec:equal}.

\section{What base 3 adds on a fixed ray}\label{sec:rayos}

We now turn to Question~\ref{q:local} of the introduction. On a ray $(a,b)$, Lemma~\ref{le:crt} turns each Fermat test into power-residue conditions at the two primes: base $t$ passes exactly when $t$ is an $a$th power modulo $p$ and a $b$th power modulo $q$. Base 2 alone already forces $q$ to be represented by an explicit binary quadratic form, by cubic, quartic or octic reciprocity. Theorem~\ref{t:exacto} shows that on five rays the base-3 test then becomes a congruence on the coordinates of that representation, and Proposition~\ref{pr:autom} shows that on the ray $(1,2)$ it is automatic. Sections~\ref{sec:criba}--\ref{sec:constantes} turn these criteria into proportions and constants, and Section~\ref{sec:datos} verifies them pair by pair up to $10^{16}$.

Let $p=ag+1<q=bg+1$ be primes greater than 3, with $(a,b)=1$ and $g>0$ even. Thus $g=(p-1,q-1)$. For $t=2,3$
write $F_t(pq)$ if $t^{pq-1}\equiv1\bmod pq$. The counting functions $N_{a,b}$, $Q^{(2)}_{a,b}$, $Q^{(3)}_{a,b}$
and $Q^{(2,3)}_{a,b}$ are those of Section~\ref{sec:criterio}. Lemma \ref{le:crt} gives the identity
\begin{equation}
 F_t(pq)\quad\Longleftrightarrow\quad
 t^g\equiv1\pmod p\quad\text{and}\quad t^g\equiv1\pmod q .
 \label{eq:crt}
\end{equation}
In a cyclic group of order $bg$, the solutions of $x^g=1$ form precisely the subgroup of $b$th powers.

Zhang and Tang \cite[\S3, (3.1)--(3.2)]{ZT03} give exact ordinary Fermat criteria when $b/a=3/2$ or $6$.
For any integer test base $t$ coprime to $pq$, both ask for the cubic character of $t$ at $q$ to be $1$.
The ratio $3/2$ also requires $(t/p)=1$, and the ratio $6$ requires $(t/q)=1$. On the ratio $4$,
Zhang \cite[(4.1)]{Zhang01} requires the quartic character of $t$ at $q$ to be $1$. He also used the
cubic character on the ratio $3$ \cite[\S5, (5.1)]{Zhang01}. The following table applies these criteria
and the octic criterion to base 3 in the coordinates determined by the base-2 test.

\begin{theorem}[Theorem~\ref{intro:rays}, criteria]\label{t:exacto}

Suppose that $F_2(pq)$. With $h\in\Z$, $k>0$ and the normalized coordinates in the table, the last column
is equivalent to $F_3(pq)$.
\[
\begin{array}{c|l|l}
(a,b)&q\text{ in base-2 coordinates}&\text{extra condition}\\ \hline
(1,3)&(3h+1)^2+27k^2   &3\mid k\\
(2,3)&(6h+1)^2+108k^2  &3\mid k,\quad 3\mid h\\
(1,4)&(4h+1)^2+64k^2   &3\mid k\\
(1,6)&(6h+1)^2+432k^2  &3\mid k\\
(1,8)&(8h+1)^2+256k^2  &8h+1\equiv c\pmod3
\end{array}
\]
In the last row $q=c^2+2d^2$, where $c\equiv1\bmod4$ and $d>0$. In particular, the last condition implies
$3\mid k$.

\end{theorem}

We give the normalization of both coordinates explicitly. In every case $p=(aq+b-a)/b$.

\subsection{Reciprocity and normalization of the coordinates}

This subsection proves Theorem~\ref{t:exacto}, first the coordinates supplied by base 2 and then the cubic, quartic and octic conditions for base 3.

\paragraph{The coordinates supplied by base 2.}

For $q\equiv1\bmod3$, the cubic residuacity of 2 is equivalent to $q=x^2+27y^2$, and for $q\equiv1\bmod4$
its quartic residuacity is equivalent to $q=x^2+64y^2$. The absolute values of the coordinates are unique
for a represented prime \cite[pp.~973, 977]{Sun}. Take $y>0$. For $(1,3)$ choose $x\equiv1\bmod3$. For $(2,3)$,
the supplementary law at 2 gives that 2 is a square modulo $p=2g+1$ precisely when $4\mid g$. This gives
$q\equiv1\bmod12$. Reducing the cubic representation modulo 4 forces $x$ to be odd and $y$ even. The sign
of $x$ is determined by $x\equiv1\bmod6$. For $(1,4)$ choose $x\equiv1\bmod4$ in the quartic representation.
For $(1,6)$ start with $q\equiv1\bmod12$. Then 2 is a square precisely when $q\equiv1\bmod24$. Reducing
the cubic representation modulo 4 shows that $x$ is odd and $y$ is even. Reduction modulo 8 gives $4\mid y$.
Thus $x=6h+1$ and, after choosing the sign, $y=4k$. Conversely, each of the forms in these four rows implies
the corresponding congruences and power-residue conditions.

On $(1,8)$ we have $q\equiv1\bmod16$. Write $q=x^2+64y^2$ when the quartic condition on 2 holds. The octic
supplement gives $(-2)^{(q-1)/8}\equiv(-1)^y\bmod q$ for these primes \cite[p.~977]{Sun}. Since $(q-1)/8$
is even, this selects precisely $2\mid y$. Writing $y=2k$ gives $x^2\equiv1\bmod16$, hence $x\equiv\pm1\bmod8$.
Exactly one sign gives $x=8h+1$. In every row, if both forms take prime values greater than 3, they determine
a unique $(h,k)$ with $k>0$. A zero second coordinate would give a square.

\paragraph{Cubic conditions.}

For a prime $q\equiv1\bmod3$ write $4q=L^2+27M^2$. The absolute values of $L,M$ are unique. The classical
criterion for 2 requires both to be even, and the criterion for 3 requires $3\mid M$ \cite[Theorem 2.4]{PCubic}.
Thus, if the criterion for 2 already holds and $q=x^2+27y^2$, then $|M|=2|y|$ and
\[
 3\in(\F_q^\times)^3\quad\Longleftrightarrow\quad 3\mid y
 \quad\Longleftrightarrow\quad q=x^2+243z^2\ \text{for some }z\in\Z.
\]
On $(1,3)$ we have $y=k$, and the test at $p$ is automatic. On $(2,3)$ we have $y=2k$. In this row $p=2g+1\equiv1\bmod4$,
so the quadratic criterion $(3/p)=1$ is equivalent to $p\equiv1\bmod12$, that is, to $3\mid g$. Since $g=(q-1)/3=12h^2+4h+36k^2\equiv h\bmod3$,
this gives the second condition. Thus the two tests in \cite[(3.1)]{ZT03} become $3\mid k$ and $3\mid h$ in
this row. On $(1,6)$ we have $q=6g+1\equiv1\bmod12$, so the quadratic condition in \cite[(3.2)]{ZT03}
is automatic for base 3. In a cyclic group, being a sixth power is equivalent to being both a square and a cube.
Since $y=4k$, the cubic condition becomes $3\mid k$.

\paragraph{The quartic condition.}

On $(1,4)$ we have $q\equiv1\bmod8$, so $-1$ is a fourth power. If $q\equiv2\bmod3$, then $(3/q)=-1$ and
3 is excluded. For $q\equiv1\bmod3$ we obtain $q\equiv1\bmod24$. The classical quartic criterion for $-3$
in $q=A^2+B^2$, with $A$ odd and $B$ even, is $3\mid B$ \cite[pp.~362--363]{PW}. Here $B=8k$. Conversely,
$3\mid k$ implies $q\equiv1\bmod3$ because $q>3$ is prime. This proves the equivalence on the entire ray.

\paragraph{The octic condition.}

On $(1,8)$, if $g\equiv2\bmod3$ then $3\mid p$, and if $g\equiv1\bmod3$ then $3\mid q$. Hence $6\mid g$
and $q\equiv1\bmod48$. Put $A=8h+1$ and $B=16k$. For $q=A^2+B^2=c^2+2d^2$ with $A\equiv c\equiv1\bmod4$,
criterion (a) at the beginning of Hudson's paper \cite{Hudson} gives
\[
 3^{(q-1)/8}\equiv1\bmod q \quad\Longleftrightarrow\quad A\equiv c\bmod3
 \qquad (q\equiv1\bmod48).
\]
The representation by $c^2+2d^2$ exists and is unique with the specified signs. This also follows from unique
factorization in $\Z[\sqrt{-2}]$ and from $(-2/q)=1$. Modulo 3, $c^2+2d^2=1$ forces $3\nmid c$ and $3\mid d$.
If $A\equiv c\bmod3$, then $A^2=1\bmod3$ and $3\mid B$. This recovers the quartic condition $3\mid k$ and
completes the proof of Theorem \ref{t:exacto}.\qed

\paragraph{Examples on the cubic and octic rays.}

On $(1,3)$, the product $11\cdot31=341$ passes base 2. The representation $31=(-2)^2+27$ gives $k=1$, so
the product fails base 3. In contrast, $103\cdot307=31\,621$ has $307=(-8)^2+27\cdot3^2$ and passes both.
These examples show how the divisibility of a coordinate distinguishes two represented primes. For $(p,q)=(1453,11617)$
we have $q=49^2+96^2=(-95)^2+2\cdot36^2$ and $49\equiv-95\bmod3$. For $(p,q)=(1531,12241)$ we have $q=(-55)^2+96^2=(-71)^2+2\cdot60^2$
and $-55\not\equiv-71\bmod3$. In both cases $k=6$ and 3 is a fourth power modulo $q$. Only the first product
passes the base-3 Fermat test.

\begin{prop}[The ray $(1,2)$]\label{pr:autom}

Let $p=g+1$ and $q=2g+1$ be primes greater than 3. Every product $pq$ passes the base-3 Fermat test. It
also passes base 2 if and only if $4\mid g$.

\end{prop}

The simultaneous family $(12u+1)(24u+1)$ is the base-2 specialization of Shanks's construction quoted
in \cite[\S4, p.~893]{GP}. The proof below gives the two individual tests on this ray.

\begin{proof}

Primality excludes $g\equiv1,2\bmod3$, since in the first case $3\mid q$ and in the second $3\mid p$. Thus
$6\mid g$, $q\equiv1\bmod12$ and $3$ is a square modulo $q$. At $p$ both tests are automatic since $p-1=g$.
At $q$ the condition is quadratic residuacity. Here $q\equiv1\bmod4$, so $(2/q)=1$ is equivalent to $q\equiv1\bmod8$,
that is, to $4\mid g$.

\end{proof}

\section{Unconditional proportions with a sifted companion}\label{sec:criba}

Theorem~\ref{t:exacto} decides each pair, but counting pairs of primes on a ray is beyond current methods. This section keeps the larger factor $q$ prime and only requires its companion to have no small prime factors. In this population the proportion of base-3 acceptance among base-2 acceptance exists unconditionally and equals $\rho_{a,b}$, the two tests are asymptotically independent (Theorem~\ref{t:criba}), and Proposition~\ref{pr:allrays} extends the proportion to every ray with $a\le2$.

In this section we vary a prime $q$ and write $g=(q-1)/b$ and $m=ag+1$. We require $g$ to be a positive
even integer and $(m,6)=1$. In Theorem~\ref{t:criba} the ray is one of the five in Theorem~\ref{t:exacto}. Define $\mathcal T_2$ by the conditions
\[
 2^g\equiv1\bmod q,\qquad a=2\ \Longrightarrow\ 4\mid g,
\]
and define $\mathcal T_{23}\subset\mathcal T_2$ by also requiring
\[
 3^g\equiv1\bmod q,\qquad a=2\ \Longrightarrow\ 3\mid g.
\]
When $m$ is prime, these are exactly the corresponding Fermat conditions for $mq$. For general $m$ these tests prescribe the order at $q$ and congruences on $g$. For instance, take $(a,b)=(1,3)$, $q=643$ and $g=214$. The companion $m=215=5\cdot43$ survives sifting at $z=3$. Both $2^{214}$ and $3^{214}$ are $1$ modulo $643$, yet $3^{mq-1}\equiv40\pmod{43}$. Thus $mq$ fails Fermat in base 3. Let $\mathcal T_0$ denote the family subject
only to the initial requirements, and let $\mathcal T_3$ denote the subfamily of $\mathcal T_0$ satisfying
the two conditions in the last display. For $i=0,2,3,23$ put
\[
 S_i(T,z)=\#\{T<q\le2T:\ q\in\mathcal T_i,\ P^-(m)>z\},
 \qquad P^-(m)=\min\{\ell\text{ prime}:\ell\mid m\}.
\]

\begin{theorem}[Theorem~\ref{intro:rays}, sifted proportions]\label{t:criba}

Let $z=z(T)\ge3$, $\log z=o(\log T)$ and
\[
 V_{a,b}(z)=\prod_{\substack{5\le\ell\le z,\ \ell\ {\rm prime}\\ \ell\nmid b-a}}
              \left(1-\frac1{\ell-1}\right).
\]
For each of those five rays and $i\in\{0,2,3,23\}$, unconditionally, $S_i(T,z)\sim \delta_i\,T V_{a,b}(z)/\log T$. In particular, both $S_{23}(T,z)/S_2(T,z)$
and $S_3(T,z)/S_0(T,z)$ tend to the same proportion $\rho_{a,b}$. Here $\delta_i$ denotes the density of $\mathcal T_i$ relative to all primes. The two local indicators are asymptotically
independent on each sifted ray. Before the additional sifting, the densities relative to all primes are
\[
\begin{array}{c|ccccc}
 (a,b)&(1,3)&(2,3)&(1,4)&(1,6)&(1,8)\\ \hline
 \delta_0&1/3&1/3&1/4&1/6&1/16\\
 \delta_3&1/9&1/18&1/16&1/18&1/64\\
 \delta_2&1/9&1/18&1/8&1/36&1/64\\
 \delta_{23}&1/27&1/108&1/32&1/108&1/256\\
 \rho_{a,b}&1/3&1/6&1/4&1/3&1/4 .
\end{array}
\]

\end{theorem}

Sections~\ref{ssec:frob} and~\ref{ssec:sift} prove Theorem~\ref{t:criba}.

\subsection{Frobenius densities}\label{ssec:frob}

We first express each condition as a Frobenius condition in a fixed normal field and compute its density.

For a Galois extension $M/K$, write $\Gal(M/K)$ for its Galois group. In the three cases involving a cubic condition take $K=\Q(\zeta_{24})$. In $K^\times/K^{\times3}$ the classes
of $2,3,\zeta_3$ are independent. Indeed, their valuations at 2 and 3 give $v(2)=4$ and $v(3)=2$ respectively,
both coprime to 3. Once the first two classes are eliminated, $\zeta_3$ is still not a cube since $\zeta_9\notin K$
by ramification at 3. By Kummer theory the extension
\[
 K(\sqrt[3]2,\sqrt[3]3,\zeta_9)/K
\]
has group $(\Z/3)^3$. It is normal over $\Q$. For each allowed class $q\bmod24$ with $q\equiv1\bmod3$, its
lifts over $\Q$ give three independent choices of phase on the two cube roots and on $\zeta_9$. Indeed,
a lift $\sigma$ maps $\sqrt[3]t$ to $\zeta_3^{e_t}\sqrt[3]t$ for $t=2,3$ and maps $\zeta_9$ to one of the
three cube roots of $\zeta_3$. Multiplying by the kernel of order 27 realizes the full Cartesian product
of these choices. For $q\equiv1\bmod3$, fixing one of the three roots of $X^3-t$ fixes all three. Triviality
of each phase on the roots of 2 and 3 is invariant under conjugation and is equivalent to the corresponding
integer being a cube modulo $q$. The action on $\zeta_9$ determines the class of $q\bmod9$. Counting these
lifts gives the proportions by Chebotarev \cite[\S3.2.3]{SerreNX}.

On $(1,3),(2,3),(1,6)$, the underlying progressions for $q$ are respectively $1\bmod3$, $1\bmod12$ and $1\bmod24$.
In the second row, $4\mid g$ comes from the quadratic condition on 2 at the companion. In the third row,
the quadratic condition on 2 at $q$ restricts $q\equiv1\bmod12$ to $q\equiv1\bmod24$. In each row, $3\nmid m$
retains two of the three lifts modulo 9. The underlying progression here incorporates the quadratic base-2
conditions before the cubic ones are imposed. Requiring 2 to be a cube contributes an independent factor
$1/3$. The densities of $\mathcal T_2$ are therefore $(1/2,1/4,1/8)\cdot (2/3)(1/3)$, in the order given.
The cubic condition on 3 contributes a further factor $1/3$. On $(2,3)$ we also require $q\equiv1\bmod9$,
one of the two admissible classes. On $(1,6)$ the quadratic condition on 3 at $q$ is automatic.

On $(1,4)$ we require complete splitting in $E=\Q(\zeta_8,\sqrt[4]2)$. This field has degree 8. Indeed,
$\sqrt2\in\Q(\zeta_8)$, whereas $\sqrt[4]2$ generates over $\Q$ a nonnormal field and therefore cannot belong
to an abelian extension. We automatically have $3\nmid m$ for $q>3$ since $4m=q+3$. Requiring 3 to be a
fourth power adjoins $\sqrt[4]3$ to $E$. The extension $E/\Q$ is unramified at 3, and $X^4-3$ is Eisenstein
over $E$ at 3. The resulting normal field has degree $4[E:\Q]$.

On $(1,8)$ the corresponding field is $E=\Q(\zeta_{48},\sqrt[8]2)$. The cyclotomic field of degree 16 contains
$\sqrt2$ and excludes $\sqrt[4]2$ by the same nonnormality argument. Since it contains all fourth roots of unity, a fourth
root of $\sqrt2$ has degree 4. At a smaller degree its square would be fixed by the entire group. Thus $[E:\Q]=64$.
The ramification index at 3 is 2, since the subfield generated by $\zeta_{16},\sqrt[8]2$ ramifies only at
2, and $\Q(\zeta_3)/\Q$ has ramification index 2 at 3. Moreover, $\sqrt3\in\Q(\zeta_{12})$. For a normalized
valuation above 3 we have $v(\sqrt3)=1$, so $X^4-\sqrt3$ is Eisenstein at $E$. Adjoining $\sqrt[8]3$ gives
a field of degree 256 over $\Q$ containing all the conjugate roots. Complete splitting in the two fields
characterizes $\mathcal T_2$ and $\mathcal T_{23}$, respectively. To compute $\delta_0,\delta_3$, remove
the base-2 restrictions. In the three rows with a cubic condition, the initial progressions become $1\bmod3,1\bmod3,1\bmod12$,
respectively. The condition $3\nmid m$ retains two classes out of three in each row. The cubic character
of 3 is still uniform. The quadratic condition at $m$ on $(2,3)$ retains one of those two classes, and quadratic
residuacity at $q$ on $(1,6)$ remains automatic. On $(1,4)$ and $(1,8)$ we start from $\Q(\zeta_8)$ and
$\Q(\zeta_{48})$, respectively. The same Eisenstein argument gives relative degree 4 for the radical of
3 over these fields. This proves all the densities in the table and the identity $\delta_3/\delta_0=\delta_{23}/\delta_2$.

\subsection{Sifting the companion of a prime}\label{ssec:sift}

We now sift the companion $m=ag+1$, using the distribution of these Frobenius sets in arithmetic progressions.

For each row and each $i$, the set $\mathcal T_i$ is a union of Frobenius classes in a fixed normal field
ramified only at 2 and 3. Each of the fields just used for an individual row has degree at most 256. Write $\Pi_i(T)=\#\{T<q\le 2T:q\in\mathcal T_i\}$.
Then $\Pi_i(T)\sim\delta_i T/\log T$. Apply the theorem of Murty and Murty \cite{MM} in the form of \cite[Lemma 1 and equation (12)]{Thorner} on the dyadic interval at level $D=T^{1/160}$. Changing its endpoint convention from $[T,2T)$ to $(T,2T]$ contributes $O(D)$ to the sum of absolute remainders, which the stated error absorbs. For a group of order $h\ge4$ the result
allows any $\theta<2/h$ for moduli coprime to the discriminant. In the lower-degree abelian cases, the Bombieri--Vinogradov
theorem suffices for the fixed progressions imposed at 2 and 3. A finite sum handles unions of classes.

For a prime $\ell\ge5$ with $\ell\nmid b-a$, the event $\ell\mid m$ selects the invertible class $q\equiv(a-b)/a\bmod\ell$.
For $\ell\mid b-a$ the event would require $q=\ell$, which is impossible for large $T$. By the Chinese remainder
theorem, for squarefree $d$ supported on the former primes, $d\mid m$ selects an invertible class modulo
$d$. It follows that
\[
 \#\{T<q\le2T:q\in\mathcal T_i,\ d\mid m\}
 =\Pi_i(T)/\varphi(d)+r_{i,d},\qquad
 \sum_{\substack{d\le D\\(d,6(b-a))=1\\ \mu^2(d)=1}}|r_{i,d}|
 \ll_A T/(\log T)^A
\]
for every $A>0$. The local sieve density is $1/(\ell-1)$ and satisfies the dimension-one condition by Mertens'
theorem. The corresponding product is $V_{a,b}(z)\asymp1/\log z$. Since $s=\log D/\log z\to\infty$, the
fundamental lemma of sieve theory \cite[Chapter~11]{FI} gives
\[
 S_i(T,z)=\Pi_i(T)V_{a,b}(z)(1+o(1))+O_A(T/(\log T)^A).
\]
The error is $o(\Pi_i(T)V_{a,b}(z))$, for example with $A=4$. If $z$ stays bounded, inclusion--exclusion
over its finitely many primes is uniform among the finitely many products involved. For an unbounded sequence
apply the fundamental lemma along subsequences with $z\to\infty$ and the finite calculation along the remaining
bounded subsequences. This proves Theorem \ref{t:criba}.\qed

The choice $z=\exp(\sqrt{\log T})$ allows the threshold
to grow faster than every fixed power of $\log T$.

On the five rays the same proof with $z=T^\epsilon$, for a sufficiently small fixed $\epsilon>0$, gives a lower bound $\gg T/(\log T)^2$
both for the classes satisfying the additional condition and for those failing it. Choose $s=1/(160\epsilon)$
large enough that the relative error in the fundamental lemma is less than $1/2$, and apply it separately
to $\mathcal T_{23}$ and $\mathcal T_2\setminus\mathcal T_{23}$, both of positive density. These companions
have a bounded number of prime factors counted with multiplicity, since $m\ll T$. An asymptotic formula
for prime pairs requires the additional primality of $m$.

\subsection{\texorpdfstring{All rays with $a\le2$}{All rays with a at most 2}}

The coordinates of Theorem~\ref{t:exacto} are specific to the five rays. The proportion is not. Let $a\in\{1,2\}$
and let $b>a$ be coprime to $a$. Define $\mathcal T_0,\mathcal T_2,\mathcal T_3,\mathcal T_{23}$ and $S_i(T,z)$ as at
the beginning of this section, with the requirement $(m,6)=1$ replaced by $(m,6b)=1$; for $b\in\{3,4,6,8\}$ the two
requirements coincide. The extra coprimality requirement excludes at most finitely many prime companions for each fixed $b$. Among the remaining prime companions, the conditions defining $\mathcal T_2$ and $\mathcal T_{23}$ are
the Fermat conditions for $mq$ in base 2 and in both bases: for $a=1$ the conditions at $m=g+1$ hold automatically,
and for $a=2$ they are $(2/m)=1$ and $(3/m)=1$, that is, $4\mid g$ and $3\mid g$ for even $g$.

\begin{prop}[All rays with $a\le2$]\label{pr:allrays}
Let $a\in\{1,2\}$, let $b>a$ be coprime to $a$, and let $z=z(T)\ge3$ with $\log z=o(\log T)$. There are positive
rational numbers $\delta_0,\delta_2,\delta_3,\delta_{23}$, depending only on $(a,b)$, such that
\[
 S_i(T,z)\sim\delta_i\,\frac{T}{\log T}
 \prod_{\substack{5\le\ell\le z,\ \ell\ {\rm prime}\\ \ell\nmid b(b-a)}}\Bigl(1-\frac1{\ell-1}\Bigr)
 \qquad(i=0,2,3,23).
\]
Moreover $\delta_{23}/\delta_2=\delta_3/\delta_0=\rho_{a,b}$, where
\[
 \rho_{1,b}=\begin{cases}1/b,& b\ \text{odd or}\ b\equiv1\pmod3,\\[2pt] 2/b,&\text{otherwise,}\end{cases}
 \qquad
 \rho_{2,b}=\begin{cases}1/b,& b\equiv1\pmod3,\\[2pt] 1/(2b),&\text{otherwise.}\end{cases}
\]
\end{prop}

For $(a,b)=(1,3),(2,3),(1,4),(1,6),(1,8)$ the formula gives $1/3,1/6,1/4,1/3,1/4$, in agreement with
Theorem~\ref{t:criba}, and for $(1,2)$ it gives $1$, in agreement with Proposition~\ref{pr:autom}. For example,
$\rho_{1,5}=1/5$, $\rho_{1,10}=1/10$, $\rho_{1,12}=1/6$, $\rho_{1,14}=1/7$ and $\rho_{2,5}=1/10$.

\begin{proof}
Put $N_b=72b^2$ and $K_b=\Q(\zeta_{N_b})$. For a prime $q\nmid N_b$ with $q\equiv1\bmod b$, the residue of $g=(q-1)/b$ modulo
$72b$ is determined by $q\bmod N_b$. Hence the class of $q$ modulo $N_b$ determines whether $g$ is even, whether
$(m,6b)=1$, whether $4\mid g$ and whether $3\mid g$, and it determines the symbols $(2/q)$ and $(3/q)$.

Let $t\in\{2,3\}$ and let $q\equiv1\bmod b$ be a prime with $q\nmid6b$. To count the primes for which $t^g\equiv1\bmod q$, put $M_b=K_b(\sqrt[b]2,\sqrt[b]3)$, which is normal over $\Q$, and let $\sigma$ be a
Frobenius element of $M_b/\Q$ at a prime above $q$. Since $q\equiv1\bmod b$, $\sigma$ fixes $\zeta_b$, and there are residues
$e_2,e_3$ modulo $b$ with $\sigma(\sqrt[b]t)=\zeta_b^{e_t}\sqrt[b]t$, and reducing
$\sigma(\sqrt[b]t)\equiv(\sqrt[b]t)^q$ gives $\zeta_b^{e_t}\equiv t^{(q-1)/b}$. As $q\nmid b$, this shows that
$t^g\equiv1\bmod q$ if and only if $e_t=0$, a condition which is invariant under conjugation.

Write $\epsilon_b=2$ if $b$ is even and $\epsilon_b=1$ if $b$ is odd. We claim that
$\kappa\mapsto(e_2(\kappa),e_3(\kappa))$ maps $\Gal(M_b/K_b)$ isomorphically onto $(\epsilon_b\Z/b\Z)^2$. Let $x=2^i3^j$
with integers $i,j$ and suppose that $x=\gamma^b$ with $\gamma\in K_b$. Since $\zeta_b\in K_b$, the positive real root
$x^{1/b}$ belongs to $K_b$. The field it generates is real, and it is normal over $\Q$ because $K_b/\Q$ is abelian. As
in the proof of Corollary~\ref{co:grh}, this forces $b\mid2i$ and $b\mid2j$. Conversely
$\sqrt2,\sqrt3\in\Q(\zeta_{24})\subset K_b$, so that $2^{b/2}$ and $3^{b/2}$ are $b$th powers in $K_b$ when $b$ is even.
Thus the subgroup of $K_b^\times/K_b^{\times b}$ generated by $2$ and $3$ has order $(b/\epsilon_b)^2$, and Kummer
theory gives $[M_b:K_b]=(b/\epsilon_b)^2$. For even $b$, every $\kappa\in\Gal(M_b/K_b)$ fixes
$\sqrt t=(\sqrt[b]t)^{b/2}$, so that $e_t(\kappa)$ is even. The claim follows by comparing orders.

Fix a class $c$ modulo $N_b$ with $c\equiv1\bmod b$, and let $\sigma_c\in\Gal(K_b/\Q)$ be the corresponding
automorphism. Its lifts to $M_b$ form a coset of $\Gal(M_b/K_b)$, and $e_t(\sigma\kappa)=e_t(\sigma)+e_t(\kappa)$ for a lift
$\sigma$ and $\kappa\in\Gal(M_b/K_b)$. By the claim, $(e_2,e_3)$ is uniformly distributed over a coset of
$(\epsilon_b\Z/b\Z)^2$ as the lift varies. For even $b$ the parity of $e_t$ is the same for all lifts: since
$\sigma(\sqrt t)=(-1)^{e_t(\sigma)}\sqrt t$, it is even exactly when $\sigma_c$ fixes $\sqrt t$, that is, when
$(t/q)=1$ for the primes $q\equiv c\bmod N_b$. By the Chebotarev density theorem, among
the primes $q\equiv c\bmod N_b$ the proportion with $2^g\equiv1\bmod q$ is $\epsilon_b/b$ if $b$ is odd or $(2/q)=1$, and
$0$ otherwise. The same holds for $3$, and the proportion with both congruences is $(\epsilon_b/b)^2$ if $b$ is odd
or $(2/q)=(3/q)=1$, and $0$ otherwise.

Each $\mathcal T_i$ is therefore a union of Frobenius classes of $M_b/\Q$, and its density $\delta_i$ is
$(\epsilon_b/b)^{k_i}/\varphi(N_b)$ times the number of classes $c$ that are admissible for $\mathcal T_i$, with
$k_0=0$, $k_2=k_3=1$ and $k_{23}=2$. The classes $c\equiv1\bmod b$ correspond to the residues of $g$ modulo $72b$.
The class is admissible for $\mathcal T_0$ when $g$ is even, $bg\not\equiv2\bmod3$, and $(ag+1,6b)=1$. For even $b$
and even $g$ we have $q\equiv1\bmod4$, so that $(2/q)=1$ means $8\mid bg$ and $(3/q)=1$ means $3\mid bg$. Hence
$\mathcal T_2$ adds conditions on $g$ modulo a power of 2 only, whereas $\mathcal T_3$ adds the conditions
\[
 3\mid bg\quad\text{if $b$ is even},\qquad 3\mid g\quad\text{if }a=2,
\]
which concern $g$ modulo 3 only. By the Chinese remainder theorem,
$\delta_{23}/\delta_2=\delta_3/\delta_0=(\epsilon_b/b)\,\vartheta$, where $\vartheta$ is the proportion of the
residues $g$ modulo 3 with $bg\not\equiv2$ and $ag\not\equiv2$ that satisfy the displayed conditions. For $a=1$ the
allowed residues are $\{0,1\}$ when $b\equiv0,1\bmod3$ and $\{0\}$ when $b\equiv2\bmod3$, so that for even $b$ we
get $\vartheta=1,\frac12,1$ for $b\equiv0,1,2\bmod3$, and $\vartheta=1$ for odd $b$. For $a=2$ the allowed residues
are $\{0,2\}$, $\{0\}$, $\{0,2\}$ for $b\equiv0,1,2\bmod3$, so that $\vartheta=\frac12,1,\frac12$. This is the
stated formula. The residue $g\equiv0\bmod72b$ is admissible for every $\mathcal T_i$, so all four densities are
positive.

The sifting argument is that of Theorem~\ref{t:criba}. The field $M_b$ is ramified only at primes dividing $6b$ and its
degree depends only on $b$, so the theorem of Murty and Murty \cite{MM}, in the form of \cite{Thorner}, gives a
level of distribution $T^{c}$ with some $c=c(b)>0$ for moduli coprime to $6b$. For a prime $\ell\ge5$ with
$\ell\nmid b(b-a)$, the event $\ell\mid m$ selects the invertible class $q\equiv(a-b)/a\bmod\ell$, and for
$\ell\mid b-a$, $\ell\nmid b$, it would require $q=\ell$. The primes dividing $b$ are excluded from the sieve by the
requirement $(m,6b)=1$. The fundamental lemma gives the asymptotic formula as before.
\end{proof}

For $a\ge3$ and $m$ prime, the test at $m$ asks for $a$th-power residuacity modulo $m$. The proof above uses explicit congruences on $g$ for $a=1,2$ and does not supply a corresponding sieve model for general $a$. At the end of Section~\ref{ssec:sixrays} we compare the formula with a census of 125 rays with $b\le300$ and all
primes $q\le2^{31}$, which also covers the case in which the companion $m$ is prime.

\section{Proportions when both factors are prime}\label{sec:pares}

Section~\ref{sec:criba} avoids conjectures on prime values. To compare the proportions with the actual pseudoprimes, whose two factors are prime, we now assume explicit instances of the Bateman--Horn conjecture. Under them the proportions $\rho_{a,b}$ hold for prime pairs, and Section~\ref{sec:constantes} computes the corresponding constants.

\begin{theorem}[Theorem~\ref{intro:rays}, prime pairs]\label{t:pares}

Assume the multivariable Bateman--Horn conjecture for prime values in the form of Conjecture A.3 of Destagnol
and Sofos \cite{DS}, for the polynomial pairs and congruence classes specified in this section. For each
of the five rays,
\[
 \frac{Q^{(2,3)}_{a,b}(X)}{Q^{(2)}_{a,b}(X)}
 \longrightarrow \rho_{a,b},
 \qquad
 \frac{Q^{(3)}_{a,b}(X)}{N_{a,b}(X)}
 \longrightarrow \rho_{a,b}.
\]
Hence $N_{a,b}Q^{(2,3)}_{a,b}/(Q^{(2)}_{a,b}Q^{(3)}_{a,b})\to1$ on each ray. The first four rays require
only pairs of quadratic polynomials in two variables. For the fifth we also use a pair of quartic polynomials
in four variables.

\end{theorem}

The rest of this section proves Theorem~\ref{t:pares}.

\paragraph{Prime-value hypotheses.}

For two distinct irreducible polynomials $f_1,f_2\in\Z[x_1,\ldots,x_n]$ write
\[
 \nu_\ell=\#\{x\in\F_\ell^n:f_1(x)f_2(x)=0\},\qquad
 \mathfrak S(f_1,f_2)=\prod_\ell\frac{1-\nu_\ell/\ell^n}{(1-1/\ell)^2}.
\]
For the fixed polynomial pairs specified below, we assume the asymptotic formula of Conjecture A.3 in \cite{DS} on each dilation of a fixed box $\mathcal B\subset\mathbb R^n$ on which both leading
forms are positive. If $f_i^\circ$ denotes the homogeneous part of highest degree, the main term is
\[
 \mathfrak S(f_1,f_2)\int_{T\mathcal B}
  \frac{d\mathbf x}{\log f_1^\circ(\mathbf x)\log f_2^\circ(\mathbf x)} .
\]
Fix a box $\mathcal B$ with $f_i^\circ>1$, rescaling if necessary. On the compact annuli used below the
leading forms are positive and bounded away from zero. After dilation by $T$, replacing $f_i^\circ$ by $f_i$
in the logarithms changes the integrand by a uniform factor $1+o(1)$. The constants below may therefore
be computed with the full polynomials. On a congruence class one applies the conjecture after substituting
the progressions for the coordinates. The families, the moduli and the boxes fixed before taking the limit
are independent of the cutoff. We justify approximation by boxes for regions with boundary of measure zero
below. This hypothesis is used from Theorem \ref{t:pares} onwards, specifically for the conditional prime-value
asymptotics and their constants.

\paragraph{Four counts in the plane.}

Start with the polynomials $F=q$ in the table of Theorem \ref{t:exacto} and with $P=(aF+b-a)/b$. In each
class $(h,k)\bmod3$ substitute $h=3u+i,\ k=3v+j$ in both forms. For every prime $\ell\ne3$ this change is
bijective over $\F_\ell^2$. At 3 all the classes with $FP\ne0\bmod3$ give the same local factor because
both polynomials take only unit values modulo 3 on the class. The real areas and the leading coefficients
are independent of $(i,j)$. The restriction $k>0$ selects a unique representation and occupies the same
real half-plane in every class. Counting the classes in the table gives the quotient of the constants.
\[
\begin{array}{c|c|c|c}
(a,b)&\text{condition for }FP\ne0\bmod3&\#\text{admissible}&
\#\text{passing base 3}\\ \hline
(1,3)&2h+1\ne0&6&2\\
(2,3)&2h+1\ne0&6&1\\
(1,4)&(h+1,k)\ne(0,0)&8&2\\
(1,6)&2h+1\ne0&6&2
\end{array}
\]
In the three cubic rows $F=1\bmod3$ and $P=1+2h\bmod3$. To compare with \emph{all} pairs in these rows,
take $4q=L^2+27M^2$, $L=1+3u$, $M>0$ and $L\equiv M\bmod2$. The normalized representation is unique. The
remaining conditions at 2 are congruences on $(u,M)$ modulo a power of 2, both for the full set and for
the base-2 subset. For example, cubic residuacity of 2 requires both $L,M$ to be even. The relevant quadratic
conditions on 2 depend on $q\bmod8$ or on $g\bmod4$. Subdivide into classes modulo a sufficiently large
power of 2, for instance $32$, to clear denominators and decide these conditions. At 3 we have $q\equiv7+6u\bmod9$.
Two values of $u\bmod3$ and all three values of $M\bmod3$ are allowed. Cubic residuacity of 3 selects $M=0\bmod3$.
On $(2,3)$ we also select the unique class giving $q=1\bmod9$. The count is the same in each allowed class
at 2. On $(1,4)$, to count all pairs we start from $q=A^2+B^2$, $A\equiv1\bmod4$, $B>0$ even and $q\equiv1\bmod8$.
The extra base-2 condition is $8\mid B$. At 3, $FP$ is a unit in eight of the nine classes of $(A,B)$, and
$3\mid B$ retains two, whether or not the base-2 condition is imposed. The quartic criterion proves the
second ratio.

The polynomials in these pairs are irreducible over $\Q$ because their definite quadratic leading forms
have no rational linear factors. After clearing the denominators each polynomial retained is primitive.
A prime dividing its content would divide all its values, which is excluded at 2 and 3 by the allowed classes
and elsewhere by the existence of a local point where both forms are units. The polynomials are distinct
and relatively prime by the affine relation with nonzero constant term. The sublattices and classes have
fixed index. The cutoff $pq\le X$ describes an ellipse $F\le T_X$, where $T_X\sim\sqrt{bX/a}$. On an annulus
away from the origin, approximate the domain from inside and outside by boxes as in Section \ref{sec:constantes}.
Each admissible class contributes the same asymptotic constant.

\subsection{\texorpdfstring{The octic count via the norm from $\Q(\zeta_8)$}{The octic count via the cyclotomic norm}}

The four rays above use one quadratic representation of $q$. The octic ray needs two.

Comparing $8h+1$ and $c$ in Theorem \ref{t:exacto} requires two quadratic norms of the same prime. To count
them together take $K=\Q(\zeta_8)=\Q(i,\sqrt{-2})$. This field has discriminant 256 and ring of integers
$\Z[\zeta_8]$. The Minkowski bound is $24/\pi^2<3$ and the unique ideal of norm 2 is generated by $1-\zeta_8$,
so the class number is one. Each prime $q\equiv1\bmod16$ splits completely in $K$. Each of its four prime
ideals has a generator $\alpha$ of norm $q$. Conversely, an element of prime norm $q>3$ gives one of these
ideals if $q\equiv1\bmod16$.

Let $\beta=\Nrm_{K/\Q(i)}\alpha$ and $\gamma=\Nrm_{K/\Q(\sqrt{-2})}\alpha$. Multiply $\beta$ by a unit $\varepsilon\in\{\pm1,\pm i\}$
chosen to give real part $A\equiv1\bmod4$ and even imaginary part $B$. Multiply $\gamma$ by $\eta\in\{\pm1\}$
to give real coefficient $c\equiv1\bmod4$. Both normalizing units depend only on $\alpha\bmod4\Ocal_K$.
Their normalized elements have absolute quadratic norm $q$. Each class modulo $16\Ocal_K$ fixes the units,
the condition $q\equiv1\bmod16$ and, when base 2 is imposed, the condition $16\mid B$. The last is the octic
base-2 criterion on this ray proved above using the supplement in \cite{Sun}.

At 3 the identity underlying the count is
\[
 \Ocal_K/3\Ocal_K\simeq\F_9\times\F_9 .
\]
Indeed, $i$ has irreducible polynomial over $\F_3$, and the two images of $\sqrt{-2}$ are $1,-1$. The order
$\Z[i,\sqrt{-2}]$ has index a power of 2 in $\Ocal_K$, since $\zeta_8=(1+i)/\sqrt2$ is recovered after inverting
2. Choose the same image of $i$ in both factors. Writing $\alpha\bmod3=(u,v)$ gives
\[
 \beta\bmod3=uv,\qquad \gamma\bmod3=(u^4,v^4),\qquad
 q\bmod3=(uv)^4 .
\]
For the first norm the conjugation interchanges the two factors. For the second it acts by $x\mapsto x^3$
in each factor. Primality of $q$ and $p=(q+7)/8$, both greater than 3, requires exactly $(uv)^4=1$ at 3.
For each $u\in\F_9^\times$ there are four choices of $v$ with $v^4=u^4=s\in\{\pm1\}$, giving 32 admissible
classes.

Now $c=\eta s\bmod3$, and $\varepsilon uv$ is one of $1,-1,i,-i$. The equality $A=c\bmod3$ holds precisely
when $\varepsilon uv=\eta s$. Thus for each $u$ the condition is exactly
\[
 v=\eta\varepsilon^{-1}u^3 .
\]
This choice satisfies $v^4=u^4$ and retains eight of the 32 classes. The calculation applies to \emph{every}
admissible class at 2 and to all eight choices of $\varepsilon,\eta$. This gives both octic ratios.

\paragraph{The prime-value hypothesis for the octic parametrization.}

In the integral basis $1,\zeta_8,\zeta_8^2,\zeta_8^3$ fix a class $\alpha_0$ modulo $48\Ocal_K$ and substitute
$\alpha=\alpha_0+48\sum_{j=0}^3t_j\zeta_8^j$. Apply Bateman--Horn to
\[
 \Nrm_{K/\Q}(\alpha),\qquad \frac{\Nrm_{K/\Q}(\alpha)+7}{8}.
\]
The second polynomial has integer coefficients on the admissible classes. The first is irreducible over
$\Q$ because its four linear factors over $\mathbb C$ are distinct conjugates in a single orbit. In independent
complex linear coordinates the other polynomial is proportional to $x_1x_2x_3x_4+7$, which is irreducible
since it is primitive and linear in $x_1$. The two polynomials are coprime. For every prime $\ell>3$ the
change of class is bijective modulo $\ell$, and specialization at $\alpha=1\bmod\ell$ gives classes with
no obstruction. At 2 and 3 the forms are units on each retained class. Thus for every fixed choice at 2
all the classes counted have the same product of local factors.

To determine the multiplicity, let $\sigma_1,\sigma_2$ be one embedding from each complex conjugate pair.
By the unit theorem the values of $\log|\sigma_1 e/\sigma_2 e|$ on the units $e$ form $R\Z$, with $R>0$
and kernel consisting of the eight roots of unity. Choose generators satisfying $0\le\log|\sigma_1\alpha/\sigma_2\alpha|<R$.
Each prime ideal then has eight representatives, and each $q$ counted gives exactly 32 generators in both
the numerator and the denominator. On the annulus $T/2<\Nrm(\alpha)\le T$ the resulting region is the dilation
by $T^{1/4}$ of a bounded set with boundary of measure zero. Inner and outer box approximations yield the
same leading integral because the boundary can be covered by boxes of arbitrarily small total volume. Rescaling
ensures positivity of the leading forms as required by the box formulation. Summing dyadic annuli recovers
the lower end of the range. For annuli of large outer norm $Y$, a fixed covering gives $O(Y/\log^2Y)$. Below
$\sqrt T$ the count of elements in the fundamental domain is $O(\sqrt T)$. Translations arising from congruence
classes preserve the leading integrals. Under the stated hypotheses, the ratio $8/32$ therefore carries
over from residue classes to prime pairs. Dividing the first ratio in Theorem~\ref{t:pares} by the second gives
$N_{a,b}Q^{(2,3)}_{a,b}/(Q^{(2)}_{a,b}Q^{(3)}_{a,b})\to1$, which completes the proof of Theorem~\ref{t:pares}.\qed

\section{Constants and windows for comparable prime factors}\label{sec:constantes}

Under the hypotheses of Section~\ref{sec:pares}, the count on each ray has the order $\sqrt X/L_X^2$ of Galway's conjecture. This section computes the constants (Proposition~\ref{pr:c2}), shows that six rays contribute a positive proportion of that order for bases 2 and 3 (Corollary~\ref{co:seis}), and gives conditional lower bounds in a dyadic window (Section~\ref{ssec:window}), which complement the conditional upper bounds of Theorem~\ref{th:KLsharp}.

Write $\Gw(X)=\sqrt X/L_X^2$ for the scale of Galway's conjecture, $B_\ell=\ell(\ell-2)/(\ell-1)^2$ and $\Ctw=\prod_{\ell>2}B_\ell$, with
the product taken over primes.

\begin{prop}[Base-2 constants for the five families]\label{pr:c2}

Under the bivariate instances of Bateman--Horn for the stated forms, one has $Q^{(2)}_{a,b}(X)\sim c_{a,b}^{(2)}\Gw(X)$
for the forms $F=q=(mh+1)^2+Dk^2$ and $P=(aF+b-a)/b$ in the classification. The area column includes the
multiplicity and the conversion from $\log q$ to $\log X$.
\[
\begin{array}{c|rr|l|l}
(a,b)&m&D&A_{a,b}&c^{(2)}_{a,b}/\Ctw\\ \hline
(1,3)&3&27&2\pi/9&16/(3\sqrt3)\\
(2,3)&6&108&\pi/(18\sqrt2)&8/(3\sqrt6)\\
(1,4)&4&64&\pi/8&4\\
(1,6)&6&432&\pi\sqrt2/36&32/(9\sqrt6)\\
(1,8)&8&256&\pi/(16\sqrt2)&6/(5\sqrt2)
\end{array}
\]

\end{prop}

\begin{proof}

All the polynomials $F,P$ are irreducible over $\mathbb Q$ since their positive definite quadratic leading
forms have no rational linear factors. They are distinct and coprime because $bP-aF=b-a\ne0$. Use the multivariable
prime-value conjecture \cite{DS}. In this setting it predicts the local product times the integral of $1/(\log F\log P)$
over the region $F\leq Q$, where $Q$ is the positive root of $Q(aQ+b-a)=bX$, with $Q\sim\sqrt{bX/a}$. For
$F=(mh+1)^2+Dk^2$ the area of the region is $\pi Q/(m\sqrt D)$. Outside a fixed region of small values,
the mean of the integrand is asymptotic to $4/\log^2 X$. Indeed, one may integrate with respect to $F$,
since the area of $F\le t$ is linear in $t$. Dividing by the multiplicity two gives $A_{a,b}=2\pi\sqrt{b/a}/(m\sqrt D)$.
To pass from the boxes in Conjecture A.3 to the ellipse, approximate a compact elliptic annulus away from
the origin in the scaled coordinates by boxes. The boundary has measure zero, and rescaling the annulus
gives the required positivity of the leading forms. The shift $1$ in $mh+1$ disappears in the scaled coordinates
and is absorbed between slightly larger and smaller annuli. The same bound $O(T/\log^2 T)$ on dyadic annuli
$T<F\le2T$ follows from a fixed covering and can be summed towards the origin. This justifies passage to
the full ellipse.

To evaluate the singular series, let $\chi$ be the primitive quadratic character of discriminant $-3$ in
the three cubic rows and of discriminant $-4$ in the two Gaussian rows. For odd primes dividing neither
the coefficients of the change of variables nor $b-a$, the resulting nondegenerate norm form over $\mathbb F_\ell$
has character $\chi(\ell)$. Its zero level has $\ell+(\ell-1)\chi(\ell)$ solutions and every nonzero level
has $\ell-\chi(\ell)$. In the split case, count pairs with prescribed product. In the inert case, use the
norm from $\mathbb F_{\ell^2}$ to $\mathbb F_\ell$. The two levels for $F=0$ and $P=0$ are distinct. The
Bateman--Horn factor is therefore
\[
 \frac{1-\{2\ell+(\ell-2)\chi(\ell)\}/\ell^2}{(1-1/\ell)^2}
       = B_\ell\left(1-\frac{\chi(\ell)}{\ell}\right).
\]
For the real quadratic character $\chi$ occurring in a row, write $L(1,\chi)=\sum_{n\ge1}\chi(n)/n$ for its Dirichlet $L$-value. Explicitly, $\chi_{-3}$ is periodic modulo 3 with values $0,1,-1$ at $0,1,2$, and $\chi_{-4}$ is periodic modulo 4 with values $0,1,0,-1$ at $0,1,2,3$. The table records the exceptional local factors and the singular series $\mathcal S=\mathfrak S(F,P)$ after extraction of $\Ctw/L(1,\chi)$.
A dash indicates that the generic formula applies.
\[
\begin{array}{c|cccc|c}
 (a,b)&\ell=2&\ell=3&\ell=5&\ell=7&\mathcal S L(1,\chi)/\Ctw\\ \hline
 (1,3)&2&3/2&-&-&8/3\\
 (2,3)&4&3/2&-&-&16/3\\
 (1,4)&4&2&-&-&8\\
 (1,6)&4&3/2&3/2&-&64/9\\
 (1,8)&4&-&-&4/3&24/5
\end{array}
\]
On $(1,3)$, the two polynomials modulo $2$ have a common zero at exactly two of the four pairs. In the remaining
rows both polynomials are always odd. In the three cubic rows $F\equiv1\pmod3$ and $P\equiv 2h+1\pmod3$.
At the collisions $(1,4)$ on $3$, $(1,6)$ on $5$ and $(1,8)$ on $7$, the character equals $-1$, and the
origin is the only common point in the norm plane. These data give the stated factors. With the primes ordered
by size the product of $1-\chi(\ell)/\ell$ is $1/L(1,\chi)$. Substituting $L(1,\chi_{-3})=\pi/(3\sqrt3)$
and $L(1,\chi_{-4})=\pi/4$ and multiplying by $A_{a,b}$ gives the five coefficients in the table.

\end{proof}

\begin{cor}[Contribution of six rays]\label{co:seis}

Under the hypotheses of Theorem \ref{t:pares} and Proposition \ref{pr:c2}, the limit $c^{(2,3)}_{a,b}$ of $Q^{(2,3)}_{a,b}(X)/\Gw(X)$
exists and equals $\Ctw$ times the corresponding entry in the table below.
\[
\begin{array}{c|ccccc}
 (a,b)&(1,3)&(2,3)&(1,4)&(1,6)&(1,8)\\ \hline
 c^{(2,3)}_{a,b}/\Ctw&16/(9\sqrt3)&4/(9\sqrt6)&1&32/(27\sqrt6)&3/(10\sqrt2).
\end{array}
\]
The five coefficients, including the factor $\Ctw$, sum to $1.9169950604\ldots$. Assuming in addition Hardy--Littlewood
for the family $(12t+1,24t+1)$ on the ray $(1,2)$ adds $2\sqrt2\,\Ctw$. On this ray the condition at $q=2p-1$
is quadratic, and $q\equiv1\bmod12$ for $p,q>3$, so base 3 is automatic. Since the rays are disjoint by
uniqueness of $(a,b)$, the total simultaneous count for semiprimes with distinct factors satisfies the conditional
lower bound
\[
 \liminf_{X\to\infty}\frac{\Qq_{23}(X)}{\Gw(X)}
 \ \ge\ \Ctw\left(2\sqrt2+\frac{16}{9\sqrt3}+\frac{4}{9\sqrt6}
             +1+\frac{32}{27\sqrt6}+\frac{3}{10\sqrt2}\right)
 =3.7842146471\ldots .
\]
The sum is over these six families.

\end{cor}

\begin{proof}

Proposition \ref{pr:c2} gives the base-2 coefficient. The first ratio in Theorem \ref{t:pares} multiplies
it by $\rho_{a,b}$ on each of the five slopes. To check the contribution of $(1,2)$ directly, base 2 requires
$4\mid g$ and primality with both factors greater than 3 requires $3\mid g$. The forms become $12t+1,24t+1$.
At 2 and 3 the local factors are $4$ and $9/4$, respectively. For $\ell>3$ they are $B_\ell$. Since $B_3=3/4$,
the singular series is $12\Ctw$. The upper endpoint is $t\sim\sqrt{X/288}$ and the two logarithms contribute
$4/(\log X)^2$. Multiplication gives $2\sqrt2\,\Ctw$. Uniqueness of the reduced ratio $(q-1)/(p-1)$ allows
us to add the six counts without overlap.

\end{proof}

\begin{cor}\label{cor:ventana}

Under the stated instances of Bateman--Horn on $(2,3)$, the number of pairs on this ray passing both bases
in the window $z<p<q\le2z$ is asymptotic to
\[
 \frac{\Ctw}{54}\frac{z}{L_z^2}.
\]
In particular, this expression gives a lower bound for $\A_{23}(z)$ in the sense of a lower limit.

\end{cor}

\begin{proof}

If $p=2g+1,\ q=3g+1$, the upper endpoint $g\le G$ corresponds to the product cutoff $X=(2G+1)(3G+1)\sim6G^2$.
Thus the coefficient of $G/(\log G)^2$ is $[4\Ctw/(9\sqrt6)]\sqrt6/4=\Ctw/9$. The window requires $(z-1)/2<g\le(2z-1)/3$.
Subtracting the asymptotic counts at these endpoints, whose ratios to $z$ tend to $1/2$ and $2/3$, gives
the factor $1/6$. On the other five rays, including $(1,2)$, we have $b/a\ge2$. For $b>2a$ the allowed interval
is eventually empty. For $b=2a$ it contains at most one integer $g$. This accounts for the difference between
the window constant and the sum of the six global contributions.

\end{proof}

\subsection{Lower bounds in a dyadic window for one and two bases}\label{ssec:window}

Theorem~\ref{th:lower} is unconditional but of logarithmic size. Under prime-value hypotheses we now obtain lower bounds of the conjectured order in every dyadic window, for one base and for two.

The family $p=4k+1$, $q=8k+1=2p-1$ in \cite[Proposition 3.1]{LRU} satisfies $p<q<2p$. However, $z<p$ and
$q\le2z$ force $z<p\le z+1/2$. For any $z$ this interval contains at most one integer $p$, and for integral
$z$ it is empty. Dickson's qualitative conjecture for those two forms therefore implies infinitude in the
region $p<q<2p$. A growing lower bound for $\A_2(z)$ in every window requires a family in the interior of
the window and a quantitative estimate.

\begin{prop}\label{pr:lowerbox}

The univariate Bateman--Horn conjecture for the two quadratics $24h^2+8h+73$ and $(6h+1)^2+108$ implies
$\A_2(z)\gg z^{1/2}/L_z^2$. The two-variable formulation of Destagnol and Sofos \cite[Conjecture A.3]{DS},
applied to an explicit pair of quadratics, implies $\A_2(z)\gg z/L_z^2$.

\end{prop}

\begin{proof}

For integers $h,k$ take
\[
 q=(6h+1)^2+108k^2,\qquad
 p=(2q+1)/3=24h^2+8h+1+72k^2 .
\]
If both values are primes greater than 3, row $(2,3)$ of Theorem \ref{t:exacto} shows that their product
passes base 2. Indeed, $q=(6h+1)^2+27(2k)^2$ makes 2 a cube modulo $q$, and $p\equiv1\bmod8$ makes it a
square modulo $p$, so both orders divide $(q-1)/3$. If $3z/2<q<2z$, both factors belong to $(z,2z]$. For
example, $h=0$ and $k=1$ give $73$ and $109$.

Setting $k=1$ gives the stated quadratics. They are distinct, positive and irreducible over $\Q$ by their
negative discriminants. For $\ell=2,3$, taking $h=0$ avoids both zeros. For $\ell\ge5$ their product is
a nonzero polynomial of degree four, so it cannot vanish on all of $\F_\ell$. Apply Bateman--Horn \cite{BH}
on an interval of $h$ between two positive multiples of $\sqrt z$ chosen so that $3z/2<q<2z$. Subtracting
the asymptotic counts at the endpoints gives $\gg\sqrt z/L_z^2$ distinct products by monotonicity of
$q$.

In two variables $p$ is strictly positive on $\mathbb R^2$ and $q$ vanishes at just one real point. Neither
has a real linear factor. They are primitive quadratics irreducible over $\Q$. At $(h,k)=(0,0)$ both have
value $1$, so there is no local obstruction. More precisely, for $\ell>3$ write $\chi(\ell)=\big(\frac{-3}{\ell}\big)$.
The equations $q=0$ and $q=-1/2$ over $\F_\ell$ are disjoint. They have $\ell+(\ell-1)\chi(\ell)$ and $\ell-\chi(\ell)$
solutions respectively. In the split case the form is a product of two independent coordinates, and in the
inert case it is the norm from $\F_{\ell^2}$. Thus the normalized local factor is $1-\chi(\ell)/\ell+O(1/\ell^2)$.
Its product converges to a positive number by the prime number theorem in arithmetic progressions modulo
$3$, or by the Euler product of $L(1,\chi)$. Choose a compact box with nonempty interior in the first quadrant
where $3/2<36h^2+108k^2<2$, with a margin at both endpoints. At scale $\sqrt z$, Conjecture A.3 gives $\gg z/L_z^2$
points with both values prime. For each prime $q>3$ the number of relevant representations is at most twelve
because $6h+1+6k\sqrt{-3}$ has norm $q$ in the Eisenstein integers, which have unique factorization, six
units and at most two ideals of norm $q$. Each element determines $(h,k)$. Dividing by this constant proves
the claim.

\end{proof}

Assume both the Kovač--Luca conjecture and this bivariate instance of Bateman--Horn. Then Theorem~\ref{th:KLsharp} gives the upper bounds and Proposition~\ref{pr:lowerbox} the lower bounds in
\[
 \lim_{z\to\infty}\frac{\log\A_2(z)}{\log z}=1,\qquad
 \lim_{X\to\infty}\frac{\log\Qq_2(X)}{\log X}=\tfrac12 .
\]
For the second limit take $z=\sqrt X/2$. This argument leaves Galway's asymptotic formula and its constant
open.

\begin{cor}\label{co:ventdoble}

Under univariate Bateman--Horn for $216t^2+24t+649$ and $(18t+1)^2+972$, we have $\A_{23}(z)\gg \sqrt z/L_z^2$.
Under the order-fiber hypothesis of Proposition \ref{pr:KL} and the bivariate instances of Corollary \ref{cor:ventana},
we also have
\[
 \lim_{z\to\infty}\frac{\log\A_{23}(z)}{\log z}=1,\qquad
 \lim_{X\to\infty}\frac{\log\Qq_{23}(X)}{\log X}=\tfrac12 .
\]

\end{cor}

\begin{proof}

Set $h=3t,\ k=3$ in row $(2,3)$ of Theorem \ref{t:exacto}. The forms $p,q$ are the stated quadratics, and
base 3 meets both divisibility conditions. The forms are primitive, positive, distinct and irreducible by
their negative discriminants. For $2$ and $3$, the choice $t=0$ avoids both zero sets. For $\ell\ge5$ the
leading coefficients are units, and the product has degree four, less than $\ell$. The singular series is
positive for this admissible pair of irreducible polynomials. Here this also follows from the quadratic
characters of the discriminants and convergence of the corresponding Dirichlet products. Subtracting the
main terms on a positive interval $t\asymp\sqrt z$ with $3z/2<q\le2z$ proves the first assertion. For the
limits, $\A_{23}\le\A_2$ and $\Qq_{23}\le\Qq_2$ allow us to apply Theorem \ref{th:KLsharp}. Corollary \ref{cor:ventana}
gives the window lower bound, and $\Qq_{23}(X)\ge\A_{23}(\sqrt X/2)$ gives the global one.

\end{proof}

\section{Interpretation for RSA-type moduli}\label{sec:rsa}

Theorems~\ref{intro:rank} and~\ref{intro:over} count products that pass fixed bases. In applications the modulus is often the product of two primes of the same size, as for RSA moduli, and the bases may be chosen at random. This section normalizes the fixed-base counts as probabilities (Corollary~\ref{co:probabilidades}) and compares them with random bases (Theorem~\ref{th:random}): for moduli whose two factors lie in $(z,2z]$, the fixed bases 2 and 3 together are passed with probability at most $z^{-2/3+o(1)}$, and random bases with probability $z^{-1+o(1)}$.

Equip $\mathcal R(z)=\{pq:z<p<q\le 2z\}$ with the uniform measure. Equivalently, choose two distinct primes
uniformly without replacement from $(z,2z]$ and forget their order. The prime number theorem gives
\[
 |\mathcal R(z)|=\binom{\pi(2z)-\pi(z)}2\sim \frac{z^2}{2L_z^2}.
\]
For $t=2,3$ we write $F_t$ also for the event that the sampled product passes base $t$, in the notation of Section~\ref{sec:rayos}.

\begin{cor}\label{co:probabilidades}

The probability that the product of the sampled factors passes Fermat for both fixed bases satisfies the
unconditional bound
\[
 \Pr_{\mathcal R(z)}(F_2\cap F_3)
 \ll \frac{1}{z^{2/3}L_z^{\dd/3}H_z^{1/2}}.
\]
For base 2 alone,
\begin{equation}
 \Pr_{\mathcal R(z)}(F_2)\ll z^{-1/2}\wt_1(z).
 \label{eq:fixed}
\end{equation}
Under the hypotheses of Corollary \ref{cor:ventana}, one has $\liminf_{z\to\infty} z\,\Pr_{\mathcal R(z)}(F_2\cap F_3)\ge\Ctw/27$.

\end{cor}

\begin{proof}

Divide the dyadic bound in Corollary \ref{cor:23}, or the asymptotic formula for the ray $(2,3)$ in Corollary
\ref{cor:ventana}, by the cardinality above. For a single base divide the first bound of Theorem \ref{th:main}
in rank one by the same cardinality.

\end{proof}

The sampling model takes as input a modulus that is already the product of two primes. In RSA key generation
the candidates for the individual prime factors undergo their own primality tests \cite[Appendices A.1 and B.3]{NIST}.
The exceptional moduli of \cite{LRU}, under the assumptions of that paper's result, belong to the set passing
the fixed base 2 and therefore satisfy its probability bound. The probability above concerns a congruence
for composite moduli drawn from the stated measure. An error probability for another procedure depends on
its input distribution and the tests it applies.

\subsection{Sampling the bases}

For a fixed modulus, the proportion of bases that pass the Fermat or the strong test has a closed form.

Fix $n=pq$ and choose $t$ uniformly in $(\Z/n\Z)^\times$. Let $\kappa_F(n)$ denote the fraction of bases passing
Fermat and $\kappa_S(n)$ the fraction passing the strong probable-prime test. For the latter write, as in Lemma~\ref{le:fuerte}, $n-1=2^s e$
with $e$ odd. A unit passes if $t^e=1\bmod n$ or if $t^{2^j e}=-1\bmod n$ for some $0\le j<s$. If $p=ug+1,\ q=vg+1$,
$(u,v)=1$ and $w=v_2(g)$, then
\begin{equation}
 \kappa_F(pq)=\frac{g^2}{(p-1)(q-1)}=\frac1{uv},
 \qquad \kappa_S(pq)=\frac{1+2/4^w}{3}\,\kappa_F(pq).\label{eq:liars}
\end{equation}
Indeed, in each component the units passing Fermat form the subgroup of order $g$. Writing $g=2^wg_0$ with
$g_0$ odd, there are $g_0$ choices of odd order and $2^{j-1}g_0$ of order with valuation $j$ for $1\le j\le w$.
All these local orders divide $n-1$. By Lemma~\ref{le:fuerte}, a base passes the strong test when both valuations are zero or both
equal the same $j>0$. In the latter case both components reach $-1$ at the same step. The number of strong
bases is therefore $g_0^2(1+\sum_{j=1}^w4^{j-1})$. The second formula is the semiprime case of Monier's
count \cite{Monier}, given explicitly in \cite[Theorem 4]{IMR}. Zhang and Tang applied the corresponding
per-modulus fraction to fixed ratios in \cite[\S5 and Table 6]{ZT03}. If one chooses from all residue
classes modulo $pq$ and tests $t^{pq}\equiv t\bmod pq$, there are $(g+1)^2$ solutions,
since in each component the zero class is added to the $g$ units. For two independent random units of the
same modulus the conditional probability is $\kappa_F(pq)^2$. On the ray $(1,2)$, for example, this probability
is $1/4$, whereas every prime pair with both factors greater than 3 that passes the fixed base 2 also passes
the fixed base 3.

\subsection{Averaging over the moduli}

Averaging this proportion over the moduli gives the probability that random bases pass.

\begin{theorem}\label{th:random}

For each fixed integer $j\ge1$, let $M_j(z)$ be the probability that $j$ bases chosen uniformly and independently
from the units of the \emph{same} uniformly sampled modulus $n\in\mathcal R(z)$ all pass Fermat. There are
an absolute integer $K\ge2$, a constant $c_j>0$ depending only on $j$, and an absolute constant $C>0$ such that, for all sufficiently large $z$,
\[
 \frac{c_j}{zL_z^{K-2}}\ \le\ M_j(z)\ \le\ \frac{C}{z}.
\]
In particular, $\log M_j(z)/\log z\to-1$. The same exponent holds for $j$ strong tests with the same sampling
rule. Under the quantitative Hardy--Littlewood conjecture for admissible pairs of linear forms, $M_j(z)\sim C_j/z$,
with
\begin{equation}
 C_j=2\!\!\sum_{\substack{u<v<2u\\(u,v)=1}}
   \frac{\mathfrak S(u,v)}{(uv)^j}
   \left(\frac{2}{v}-\frac1u\right)>0 ,
 \label{eq:HLconstant}
\end{equation}
where $\mathfrak S(u,v)$ is the singular series of $ug+1,vg+1$.

\end{theorem}

\begin{proof}

Conditional independence gives $M_j(z)=|\mathcal R(z)|^{-1}\sum_n\kappa_F(n)^j$. The representation by $g,u,v$
is unique if $(u,v)=1$. Since $q\le2z<2p$, we have $v\le2u$ by $g(v-2u)<1$ and integrality. For each $u<v\le2u$
the number of contributing values of $g$ is at most $2z/v$. If $v>\sqrt z$, their total contribution to
the unnormalized sum is
\[
 \ll_j z\sum_{v>\sqrt z}v^{-2j}\ll_j z^{3/2-j}.
\]
For $v\le\sqrt z$ apply Lemma \ref{le:sieve}. The forms are admissible since $g=0$ avoids both roots for
every prime. With $\Delta=uv(v-u)$ and any $\epsilon>0$, the number of $g\le2z/v$ is
\[
 \ll \frac{z}{vL_z^2}
        \left(\frac{\Delta}{\varphi(\Delta)}\right)^2
 \ll_\epsilon \frac{z v^\epsilon}{vL_z^2},
\]
where only values of $g$ making both forms prime are counted. The second inequality uses $\Delta\le v^3$
and Mertens' theorem. For each $v$ we have $\sum_{v/2\le u<v}(uv)^{-j}\ll_j v^{1-2j}$. The sum over $v$
converges even for $j=1$ on taking $\epsilon=1/2$. Thus $\sum_n \kappa_F(n)
\ll z/L_z^2$. Since $0\le \kappa_F(n)\le1$,
the same majorant works for every $j\ge1$ with the same constant.

For the reverse inequality use the quantitative form of Maynard's theorem for linear forms in an interval
\cite[Theorem 3.4]{Maynard}. Fix a value of $K$ that guarantees two prime values. Take the $K$ coefficients
$c_i=4K+i$, $1\le i\le K$. The family $c_i t+1$ is admissible. For large $z$ and $z/(4K)\le t\le z/(3K)$
all the forms lie strictly between $z$ and $2z$. The interval has length comparable to $z$ and meets the
length hypothesis of the theorem. In the short-interval version take the left endpoint as the parameter
and $y$ as the length, which is comparable to the parameter. With $m=2$ and, for example, $\epsilon=1/6$,
one has $x^{7/12+\epsilon}\le y\le x$ for large $x$. The coefficients $c_i$ are fixed and the constant terms
are 1, as required. Thus $\gg_K z/L_z^K$ values of $t$ yield at least two primes. For any selected pair,
$g=(p-1,q-1)\ge t$ and $\kappa_F(pq)\ge1/(25K^2)$. Each product admits at most $\binom K2$ such representations
since, for fixed slopes, the ordered primes determine $t$. The use of slopes taking two simultaneous prime
values appears in the proof of Li and Pomerance \cite[Theorem 2]{LiP}. Here we restrict the slope ratio
and count in every window. Dividing by $|\mathcal R(z)|$ proves the lower bound. By \eqref{eq:liars}, $\kappa_F(n)/3\le \kappa_S(n)\le \kappa_F(n)$,
so the means of $\kappa_S(n)^j$ lie between $3^{-j}M_j(z)$ and $M_j(z)$.

For the last assertion let $\rho_\ell$ be the number of roots of $(ug+1)(vg+1)$ modulo $\ell$ and let
$\mathfrak S(u,v)=\prod_\ell(1-\rho_\ell/\ell)(1-1/\ell)^{-2}$. For fixed $u,v$ with $v<2u$, the interval
for $g$ is $((z-1)/u,(2z-1)/v]$ and has length $z(2/v-1/u)+O(1)$. Hardy--Littlewood gives the leading prime-pair
term with factor $\mathfrak S(u,v)/L_z^2$. The boundary case $v=2u$ forces $(u,v)=(1,2)$ and contributes
$O(1)$. After multiplication by $L_z^2/z$, the tail with $V<v\le\sqrt z$ is $O_j(\sum_{v>V}v^{-2j+1/2})$
uniformly in $z$. The remaining tail tends to zero. We may thus first take the limit at fixed $V$ and then
let $V\to\infty$, without assuming conjectural uniformity in the slopes. The factor $2$ in \eqref{eq:HLconstant}
comes from $|\mathcal R(z)|$. The sum is already positive at $(u,v)=(2,3)$.

\end{proof}

\begin{remark}[Fixed powers of a base]\label{re:powers}
The construction in the proof of Theorem \ref{th:random} also gives an unconditional lower bound of exponent $1$ in the window
for bases that are suitable fixed powers. Take $D=\lcm(c_1,\ldots,c_K)$. For factors $p=c_it+1$ and $q=c_jt+1$, both $p-1$
and $q-1$ divide $Dt$. Once $z$ exceeds the primes dividing a fixed integer $a\ge2$, $(a^D)^t=1$ holds in both components.
Since $t\mid(p-1,q-1)$, Lemma \ref{le:crt} shows that the $\gg z/L_z^K$ distinct products of that construction pass the
base $a^D$. Hence $\A_{a^D}(z)\gg z/L_z^K$ for all sufficiently large $z$, simultaneously for the $D$th powers of any fixed
list of integer bases greater than one; $D$ and $K$ are absolute.
\end{remark}

\section{Numerical verification and coverage of the rays}\label{sec:datos}

The computations below give the small tables of pseudoprime counts, check the criteria of Theorem~\ref{t:exacto} pair by pair on the six rays, the proportions of Proposition~\ref{pr:allrays} on further rays, the share of the six rays in the global count, and the decimals of the constants.

\subsection{Direct modular powers and the strong test}\label{ssec:direct}

We first count the small products that pass one, two or three fixed bases, and the base-2 strong test.

The table counts products of two distinct odd primes. The Fermat columns were computed by modular
exponentiation modulo the product, and the strong column by successive squaring modulo the product. The last two columns require simultaneous
acceptance for the indicated bases.

\begin{center}\small
\begin{tabular}{@{}r rrrr@{}}
\toprule
Cutoff $pq\le X$ & Fermat $2$ & Strong $2$ & Fermat $2,3$ & Fermat $2,3,5$\\
\midrule
$10^4$ & 11 & 5 & 1 & 0\\
$10^5$ & 34 & 12 & 6 & 0\\
$10^6$ & 107 & 37 & 13 & 2\\
$10^7$ & 311 & 134 & 51 & 4\\
\bottomrule
\end{tabular}
\end{center}

The window $10^4<p<q\le2\cdot10^4$ contains $533\,028$ pairs. Of these, $37$ pass the base-2 Fermat test,
$16$ are strong base-2 pseudoprimes, $2$ pass Fermat for both 2 and 3, and none of the $533\,028$ passes
all three bases 2, 3 and 5.

\subsection{\texorpdfstring{Six rays up to a product of $10^{16}$}{Six rays up to a product of 10 to the power 16}}\label{ssec:sixrays}

Next we verify Theorem~\ref{t:exacto} and Proposition~\ref{pr:autom} pair by pair on the six rays.

The sieve of Eratosthenes supplies primes in both progressions. The first enumeration computes $t^{pq-1}\bmod pq$
for $t=2,3$ directly modulo the product, with 128-bit intermediate products. The second runs through the
quadratic forms with $k>0$ and both signs of $h$, retaining each prime $q$ with its companion. For $(1,8)$
it computes $c$ by a separate enumeration of $q=c^2+2d^2$. Comparing the sorted sets of $q$ on each ray
verifies the identity pair by pair, for base 2 and for the two bases together, up to $pq\le10^{16}$. Integer
binary search determines the admissible interval for $g$. The counts exclude $p\le3$ and all squares.

\begin{center}\small
\begin{tabular}{@{}crrrrrr@{}}\toprule
Ray& Prime pairs& Base 2& Base 3& Both& $100Q^{(2,3)}/Q^{(2)}$& $100Q^{(3)}/N$\\\midrule
$(1,2)$&311\,373&155\,594&311\,373&155\,594&100.00&100.00\\
$(1,3)$&508\,295&169\,472&169\,448&56\,377&33.27&33.34\\
$(2,3)$&359\,285&59\,900&60\,142&10\,057&16.79&16.74\\
$(1,4)$&440\,125&220\,096&109\,977&55\,156&25.06&24.99\\
$(1,6)$&480\,686&79\,983&160\,329&26\,793&33.50&33.35\\
$(1,8)$&187\,167&46\,423&46\,765&11\,520&24.82&24.99\\
\bottomrule\end{tabular}
\end{center}

The following table gives examples for each ray. Each entry is the first pair by increasing product in its
column. The direct test modulo $pq$ was compared with the coordinate criterion of Section \ref{sec:rayos}.

\begin{center}\small
\begin{tabular}{@{}ccc@{}}\toprule
Ray & Passes both bases 2 and 3 & Passes base 2, fails base 3\\\midrule
$(1,2)$ & $(37,73)$ & impossible by Proposition \ref{pr:autom}\\
$(1,3)$ & $(103,307)$ & $(11,31)$\\
$(2,3)$ & $(2593,3889)$ & $(73,109)$\\
$(1,4)$ & $(151,601)$ & $(19,73)$\\
$(1,6)$ & $(677,4057)$ & $(73,433)$\\
$(1,8)$ & $(1453,11617)$ & $(43,337)$\\\bottomrule
\end{tabular}
\end{center}

Of the $46\,423$ base-2 candidates on the last ray, $23\,192$ pass the quartic condition on 3 and $11\,520$
its octic condition. The five rays studied here contribute $159\,903$ simultaneous products at the cutoff.
Including $(1,2)$ brings the total to $315\,497$. The verification files give the intermediate counts $10^8,10^{10},10^{12},10^{14}$
and an independent enumeration of the finite local classes.

\paragraph{Rays outside the table.}
For Proposition~\ref{pr:allrays} we counted, for every prime $q$ with $2^{16}<q\le2^{31}$, the four families
$\mathcal T_0,\mathcal T_2,\mathcal T_3,\mathcal T_{23}$ on 125 rays with $a\le2$: every $b\le64$ (odd $b$ when
$a=2$), eighteen further values $72\le b\le300$ for $a=1$ and thirteen further odd values $75\le b\le225$ for
$a=2$. The counts are taken by direct modular exponentiation in the intervals $T<q\le2T$ with $T=2^k$,
$16\le k\le30$, for three conditions on the companion $m=ag+1$: only $(m,6b)=1$, which gives $S_i(T,3)$;
$P^-(m)>100$, which gives $S_i(T,100)$; and $m$ prime, which counts the prime pairs on the ray. The table gives
the ratio $\#\mathcal T_{23}/\#\mathcal T_2$ in the top interval $2^{30}<q\le2^{31}$, which contains
$50\,697\,537$ primes, for ten rays.
\begin{center}\small
\begin{tabular}{@{}crcccc@{}}\toprule
Ray & $S_2(T,3)$ & $(m,6b)=1$ & $P^-(m)>100$ & $m$ prime & $\rho_{a,b}$\\\midrule
$(1,5)$ & $1\,014\,326$ & $0.19988$ & $0.20004$ & $0.20153$ & $1/5$\\
$(1,7)$ & $1\,035\,791$ & $0.14264$ & $0.14343$ & $0.14254$ & $1/7$\\
$(1,10)$ & $507\,314$ & $0.10006$ & $0.10015$ & $0.10036$ & $1/10$\\
$(1,12)$ & $704\,663$ & $0.16622$ & $0.16700$ & $0.16555$ & $1/6$\\
$(1,14)$ & $129\,358$ & $0.14331$ & $0.14306$ & $0.14279$ & $1/7$\\
$(1,16)$ & $396\,104$ & $0.06261$ & $0.06204$ & $0.06283$ & $1/16$\\
$(1,60)$ & $28\,234$ & $0.03248$ & $0.03161$ & $0.03194$ & $1/30$\\
$(1,100)$ & $10\,191$ & $0.00942$ & $0.00862$ & $0.00856$ & $1/100$\\
$(2,5)$ & $1\,013\,602$ & $0.09977$ & $0.10012$ & $0.10010$ & $1/10$\\
$(2,7)$ & $258\,898$ & $0.14229$ & $0.14191$ & $0.14208$ & $1/7$\\\bottomrule
\end{tabular}
\end{center}
We measure the deviation of $\#\mathcal T_{23}$ from $\rho_{a,b}\,\#\mathcal T_2$ in binomial standard deviations,
on the rays other than $(1,2)$ with $\rho_{a,b}(1-\rho_{a,b})\,\#\mathcal T_2\ge25$. In the top interval these are
106 rays for $(m,6b)=1$, with mean $-0.09$, standard deviation $0.92$ and largest absolute value $2.57$; 98 rays
for $P^-(m)>100$, with mean $-0.14$, standard deviation $0.90$ and largest absolute value $2.69$; and 88 rays for
prime pairs, with mean $-0.18$, standard deviation $1.03$ and two values beyond $3$, the largest $3.24$. The
counts $S_2(T,3)$ and $S_{23}(T,3)$ also agree with $\delta_2N$ and $\delta_{23}N$, where $N$ is the number of primes
in the interval and $\delta_i$ are the exact densities of Proposition~\ref{pr:allrays}: the standardized
deviations have standard deviation $1.03$ and $0.91$. We see no loss of accuracy as $b$ grows. At smaller heights
the ratio tends to lie below $\rho_{a,b}$. For $k\le24$ the mean deviation is negative in eight of the nine
intervals, between $-0.65$ and $-0.29$, while for $k=29$ and $k=30$ it is $-0.04$ and $-0.09$. The cumulative
counts up to $q\le2\cdot10^8$ on 23 rays, also included in the replication package, show the same effect: there the
ratio lies below $\rho_{a,b}$ on 21 of the 22 rays other than $(1,2)$. The binomial scale is a descriptive
normalization for deterministic counts, and we do not claim a lower-order term.

\subsection{A census of all simultaneous semiprime pseudoprimes}

To place the families in the global count we filtered an initial segment of Feitsma and Galway's factored
catalog \cite{FG}. Each line lists the prime factors with multiplicity of a base-2 pseudoprime. We reconstructed
each product and checked strict increase up to the first entry exceeding $10^{16}$. We retained products
with exactly two distinct prime factors $p<q$, each of multiplicity one. For each product we compared $3^{pq-1}\bmod pq$
with the two powers $3^g$ modulo $p$ and $q$, with $g=(p-1,q-1)$. We classified rays from the factors themselves,
independently of the representations by quadratic forms. Completeness of the global census up to the cutoff
relies on completeness of the base-2 catalog \cite{FG}. We also performed an independent global enumeration
up to $10^9$. Run through $p\le\sqrt X$. For $p<37$ factor $2^{p-1}-1$, and for $p\ge37$ sieve the primes
$p<q\le X/p$ in $q\equiv1\bmod\operatorname{ord}_p 2$. Then apply the base-2 Fermat test and the common-order
criterion for base 3. Lemma \ref{le:crt} guarantees coverage. At all six powers of ten from $10^4$ to $10^9$,
this procedure recovers both global columns. The last cutoff gives $2\,455$ and $376$. The preceding table
provides a further independent direct check of the counts on the six rays.

\begin{center}
\begin{tabular}{@{}rrrr@{}}\toprule
$X$ & Base 2, all rays & Bases 2 and 3, all rays & Both bases, six rays\\\midrule
$10^{4}$ & 11 & 1 & 1\\
$10^{6}$ & 107 & 13 & 13\\
$10^{8}$ & 880 & 140 & 124\\
$10^{10}$ & 6\,502 & 1\,010 & 833\\
$10^{12}$ & 46\,071 & 7\,135 & 5\,790\\
$10^{14}$ & 334\,565 & 51\,869 & 41\,768\\
$10^{16}$ & 2\,520\,501 & 394\,131 & 315\,497\\
\bottomrule\end{tabular}
\end{center}

The last catalog entry included was $9\,999\,999\,191\,658\,001$ and the first excluded was $10\,000\,006\,286\,491\,369$.
We removed the squares $1093^2$ and $3511^2$. The six families account for $315\,497/394\,131$ of the simultaneous
products at the cutoff, about $80.05\%$. The remaining $78\,634$ have other slopes. Among all base-2 candidates,
the observed base-3 acceptance fraction is $394\,131/2\,520\,501$, about $15.64\%$. These are finite-sample
proportions. The limits in Theorem \ref{t:pares} fix the ray before sending $X$ to infinity.

\subsection{Rational enclosures for the constants}

The exact expressions in the theorems determine the constants. To check the decimals quoted in the text,
the supplemental program produces the intervals below with outward-rounded endpoints. The displayed decimal endpoints are rational numbers.
\[
\begin{array}{l|cc}
 &\text{lower}&\text{upper}\\ \hline
 \Ctw &0.660161815846869&0.660161815846870\\
 \sum_{\text{five rays}} c^{(2,3)}_{a,b}
 &1.916995060445685&1.916995060445686\\
 \sum_{\text{six rays}} c^{(2,3)}_{a,b}
 &3.784214647108670&3.784214647108671\\
 C_\varphi/(3\log2)&0.934672313630256&0.934672313630257\\
 \lambda&0.910549855132844&0.910549855132845
\end{array}
\]
We give the error bounds for the infinite product, which converges slowly in its original form. For $k\ge2$
let $b_k=k^{-1}\sum_{d\mid k}\mu(d)2^{k/d}$. The integer $kb_k$ counts binary words of least period $k$,
so $0\le b_k\le 2^k/k$. Möbius inversion and comparison of logarithmic coefficients give
\[
 \frac{1-2x}{(1-x)^2}=\prod_{k\ge2}(1-x^k)^{b_k},\qquad |x|<1/2.
\]
For $P=31$ and $K=16$ we compute
\[
 U=\prod_{2<p\le P}B_p\ \prod_{k=2}^K
     \left(\zeta(k)\prod_{p\le P}(1-p^{-k})\right)^{-b_k}.
\]
Each omitted factor in approximating $\Ctw$ by $U$ lies between zero and one. Since $-\log(1-t)\le2t$ for
$0\le t\le1/2$ and $\sum_{n>P}n^{-k}\le P^{1-k}/(k-1)$, we have
\[
 0\le\log(U/\Ctw)
 \le2\sum_{k>K}\frac{2^k}{k}\frac{P^{1-k}}{k-1}
 \le E,\qquad E=\frac{2P(2/P)^{K+1}}{1-2/P}.
\]
Hence $U(1-E)\le\Ctw\le U$. For each $\zeta(k)$ we apply Euler--Maclaurin starting at $N=64$ through $B_{16}$
inclusive. The remainder has absolute value at most
\[
 \frac{|B_{16}|}{16!}\,k(k+1)\cdots(k+14)\,N^{-k-15}.
\]
This follows by integrating the derivative of order 16 of $t^{-k}$ against the periodic Bernoulli polynomial,
whose absolute value is at most $|B_{16}|$. For $\log2$ sum the first 50 terms of $2\sum_{j\ge0}3^{-2j-1}/(2j+1)$
and bound the tail by a geometric series. Square roots are enclosed by scaled integer square roots, and
$\lambda$ by enclosing the minimizer via bisection of the strictly increasing derivative numerator of Section \ref{sec:fibras}. For the bisection
enclose $\exp y$, $0\le y\le2$, using the Taylor polynomial of degree 48 with remainder at most $2^{49}/(49!(1-2/50))$.
All operations use integers, fractions or directed rational rounding to a scale of $10^{64}$. For Corollary
\ref{co:seis} the decimals verify the exact coefficients. Their interpretation as asymptotic coefficients
for prime counts retains the hypotheses stated there.

\paragraph{Reproducibility.}

The ancillary files contain the C++17 program for the form enumeration and modular tests, the output for each
ray, the Python program for the rays of Proposition~\ref{pr:allrays} up to $q\le2\cdot10^8$, the C++17 census of 125 rays up to $q\le2^{31}$ with its summary program, the Python program counting classes over the finite fields, and the Feitsma--Galway catalog reader with its
cumulative CSV file. A further program enumerates the small products of Section~\ref{ssec:direct}. Counts and congruence decisions use integer arithmetic. Some summaries use floating-point ratios or standardized deviations. The interval program uses directed rational arithmetic and reproduces the five
rational enclosures in this section. The programs and reference outputs are listed in the accompanying \texttt{README.md}, and the same package is archived on Zenodo \cite{Zenodo}. The finite census and the decimal certificates can be reproduced from them. The theorem proofs use the exact formulas and the stated mathematical hypotheses.

\section{Scope and open questions}\label{sec:fin}

We summarize what the three theorems settle and what remains open.

The rank-two argument counts all ratios of the factors and gives an upper bound even for very unbalanced
semiprimes. The local criterion answers an individual question for each prime pair of a fixed ratio. On
$(1,2)$ every prime pair with factors greater than 3 satisfies the base-3 test. On $(1,8)$ the quartic condition
leaves two octic possibilities, distinguished by the normalized sign of a coordinate modulo 3.
Proposition~\ref{pr:allrays} gives the proportion on every ray with $a\le2$. For $a\ge3$ our argument gives no general reduction at the smaller factor to a congruence on $g$; coordinates such as those of
Theorem~\ref{t:exacto} are available to us only on the five rays.

Recovering the five proportions with a prime in place of the sifted companion, without assuming a conjecture
on prime values, remains open. A growing lower bound in \emph{every} window for the fixed base 2, or for
the pair of bases 2 and 3, also remains open; for suitable fixed powers $a^D$ of any base, Remark~\ref{re:powers} gives an unconditional lower bound of exponent $1$ in every window. The logarithmic construction in Section \ref{sec:inferior}
allows ratios $q/p$ without a prescribed upper bound. Deriving a global asymptotic from constants on fixed
rays would further require uniform control of growing slopes and their tail. The empirical proportion $80.05\%$
and the lower limits in Section \ref{sec:constantes} concern different cutoffs and different assertions.
The single-base Galway problem and its global constant are studied in \cite{Galway,GP}. The simultaneous
problem also requires the interaction of two fixed bases for the same product.

Several bounds remain open. For base 2 alone no power saving below $X^{3/4}$ is known, even under GRH, since
Corollary~\ref{co:grh} requires rank at least two. The conjectured order is $\sqrt X/(\log X)^2$ \cite{Galway,GP},
and for the overpseudoprimes Pomerance conjectured $\mathrm{Ov}_2(x)=o(x^{1/2+\epsilon})$ \cite[Remark~2]{Shev08};
Theorem~\ref{intro:over} saves only powers of $\log x$. In the critical configuration for the exponent $3/4$ the order
of $2$ modulo each factor is close to the square root of that factor (Section~\ref{sec:fibras}). By
Corollary~\ref{co:erdos}, a power saving for $\Qq_2$, for $\mathcal E_2$ or for $\mathrm{Ov}_2$ is equivalent to a
bound $O(y^{1-\eta})$ for the number of primes $p\le y$ with $\ord_p(2)\le p^{1/2+\delta}$, for some
$\delta,\eta>0$. This is a power-saving form of a conjecture of Erdős \cite{Erdos76}. Erdős and Murty \cite{EM}
proved that, for every positive function $\varepsilon(p)$ tending to $0$, one has $\ord_p(2)>p^{1/2+\varepsilon(p)}$
for almost all primes $p$. Corollary~\ref{co:erdos} asks for more in two respects: a fixed $\delta>0$ in place of
$\varepsilon(p)$, and a bound $O(y^{1-\eta})$ for the number of exceptions up to $y$. For the bases 2 and 3
together the exponents are $2/3$ unconditionally and $13/20+\varepsilon$ under the stated GRH for each $\varepsilon>0$, and our growing lower bounds for $\Qq_{23}$ assume simultaneous prime-value conjectures. Finally, Pomerance \cite{Pom82} observed that he could not show that the number of
pseudoprimes up to $x$ with at most $k$ distinct prime factors, divided by $\log x$, is unbounded. For $k=2$,
Theorem~\ref{th:lower} reaches this scale with an explicit constant; to our knowledge unboundedness remains open.

\section*{Use of artificial intelligence}
OpenAI ChatGPT and Codex assisted with mathematical exploration, drafting, reference checks and computational verification.


\begin{thebibliography}{99}

\bibitem{Lehmer} D.~H. Lehmer, \emph{On the converse of Fermat's theorem}, Amer. Math. Monthly \textbf{43} (1936), 347--354.

\bibitem{E49} P.~Erdős, \emph{On the converse of Fermat's theorem}, Amer. Math. Monthly \textbf{56} (1949), 623--624. \doi{10.2307/2304732}.

\bibitem{BZ} K.~Zsigmondy, \emph{Zur Theorie der Potenzreste}, Monatsh. Math. Phys. \textbf{3} (1892), 265--284. \doi{10.1007/BF01692444}.

\bibitem{Schinzel} A.~Schinzel, \emph{On primitive prime factors of $a^n-b^n$}, Proc. Cambridge Philos. Soc. \textbf{58} (1962), 556--562. \doi{10.1017/S0305004100040561}.

\bibitem{Somer} L.~Somer, \emph{On superpseudoprimes}, Math. Slovaca \textbf{54} (2004), 443--451. \url{http://hdl.handle.net/10338.dmlcz/130094}.

\bibitem{Bateman} P.~T. Bateman, \emph{The distribution of values of the Euler function}, Acta Arith. \textbf{21} (1972), 329--345. \doi{10.4064/aa-21-1-329-345}.

\bibitem{PSW} C.~Pomerance, J.~L. Selfridge and S.~S. Wagstaff, Jr., \emph{The pseudoprimes to $25\cdot 10^9$}, Math. Comp. \textbf{35} (1980), 1003--1026. \doi{10.1090/S0025-5718-1980-0572872-7}.

\bibitem{Jaeschke} G.~Jaeschke, \emph{On strong pseudoprimes to several bases}, Math. Comp. \textbf{61} (1993), 915--926. \doi{10.1090/S0025-5718-1993-1192971-8}.

\bibitem{Pom81} C.~Pomerance, \emph{On the distribution of pseudoprimes}, Math. Comp. \textbf{37} (1981), 587--593. \doi{10.1090/S0025-5718-1981-0628717-0}.

\bibitem{Wagstaff} S.~S. Wagstaff, Jr., \emph{Pseudoprimes and a generalization of Artin's conjecture}, Acta Arith. \textbf{41} (1982), 141--150. \doi{10.4064/aa-41-2-141-150}.

\bibitem{Hooley} C.~Hooley, \emph{On Artin's conjecture}, J. Reine Angew. Math. \textbf{225} (1967), 209--220. \doi{10.1515/crll.1967.225.209}.

\bibitem{EM} P.~Erdős and M.~R. Murty, \emph{On the order of $a \pmod p$}, in: Number Theory (Ottawa, ON, 1996), CRM Proc. Lecture Notes \textbf{19}, Amer. Math. Soc., Providence, RI, 1999, 87--97. \doi{10.1090/crmp/019/10}.

\bibitem{Matt} C.~R. Matthews, \emph{Counting points modulo $p$ for some finitely generated subgroups of algebraic groups}, Bull. London Math. Soc. \textbf{14} (1982), 149--154. \doi{10.1112/blms/14.2.149}.

\bibitem{GM} R.~Gupta and M.~R. Murty, \emph{A remark on Artin's conjecture}, Invent. Math. \textbf{78} (1984), 127--130. \doi{10.1007/BF01388719}.

\bibitem{Pappa} F.~Pappalardi, \emph{On the order of finitely generated subgroups of $\Q^\times \pmod p$ and divisors of $p-1$}, J. Number Theory \textbf{57} (1996), 207--222. \doi{10.1006/jnth.1996.0044}.

\bibitem{SW} J.~Sorenson and J.~Webster, \emph{Strong pseudoprimes to twelve prime bases}, Math. Comp. \textbf{86} (2017), 985--1003. \doi{10.1090/mcom/3134}.

\bibitem{Ford} K.~Ford, \emph{The distribution of integers with a divisor in a given interval}, Ann. of Math. (2) \textbf{168} (2008), 367--433. \doi{10.4007/annals.2008.168.367}.

\bibitem{Galway} W.~F. Galway, \emph{The density of pseudoprimes with two prime factors}, preliminary report, Abstracts Amer. Math. Soc., abstract 962-11-1393, 2000. \url{https://jointmathematicsmeetings.org/meetings/national/jmm-archive/962-11-1393.pdf}.

\bibitem{GP} A.~Granville and C.~Pomerance, \emph{Two contradictory conjectures concerning Carmichael numbers}, Math. Comp. \textbf{71} (2002), 883--908. \doi{10.1090/S0025-5718-01-01355-2}.

\bibitem{FG} J.~Feitsma and W.~F. Galway, \emph{Tables of pseudoprimes and related data}, factored catalog of base-2 Fermat pseudoprimes below $2^{64}$, Centre for Experimental and Constructive Mathematics, Simon Fraser University. \url{https://www.cecm.sfu.ca/Pseudoprimes/}.

\bibitem{DJU10} L.~Dieulefait and J.~Jiménez Urroz, \emph{Small primitive roots and malleability of RSA moduli}, J. Comb. Number Theory \textbf{2} (2010), 171--179; preprint version \doi{10.48550/arXiv.0801.0030}. The estimate is also stated in \cite[\S1]{LRU}.

\bibitem{DJU} L.~V. Dieulefait and J.~Jiménez Urroz, \emph{Factorization and malleability of RSA moduli, and counting points on elliptic curves modulo $N$}, Mathematics \textbf{8} (2020), art.~2126. \doi{10.3390/math8122126}.

\bibitem{LRU} F.~Luca, D.~Ralaivaosaona and J.~Jiménez Urroz, \emph{On pseudoprime RSA moduli}, Afr. Mat. \textbf{36} (2025), art.~118. \doi{10.1007/s13370-025-01338-1}.

\bibitem{LiP} S.~Li and C.~Pomerance, \emph{Fermat pseudoprimes}, Acta Arith. \textbf{222} (2026), 191--196. \doi{10.4064/aa250411-9-10}.

\bibitem{KL} V.~Kovač and F.~Luca, \emph{On the number of divisors of Mersenne numbers}, Exp. Math. (2026), \doi{10.1080/10586458.2026.2636060}; preprint version \doi{10.48550/arXiv.2506.04883}.

\bibitem{Pollack} P.~Pollack, \emph{Nonnegative multiplicative functions on sifted sets, and the square roots of $-1$ modulo shifted primes}, Glasg. Math. J. \textbf{62} (2020), 187--199. \doi{10.1017/S0017089519000041}. Corrected author version: \url{https://www.pollack-math.net/brunmultiplicative4.pdf}; errata: \url{https://www.pollack-math.net/errata-nonnegative.pdf}.

\bibitem{Zhang01} Z.~Zhang, \emph{Finding strong pseudoprimes to several bases}, Math. Comp. \textbf{70} (2001), 863--872. \doi{10.1090/S0025-5718-00-01215-1}.

\bibitem{ZT03} Z.~Zhang and M.~Tang, \emph{Finding strong pseudoprimes to several bases.~II}, Math. Comp. \textbf{72} (2003), 2085--2097. \doi{10.1090/S0025-5718-03-01545-X}.

\bibitem{Hudson} R.~H. Hudson, \emph{Diophantine determinations of $3^{(p-1)/8}$ and $5^{(p-1)/4}$}, Pacific J. Math. \textbf{111} (1984), 49--55. \doi{10.2140/pjm.1984.111.49}.

\bibitem{MM} M.~R. Murty and V.~K. Murty, \emph{A variant of the Bombieri--Vinogradov theorem}, in: Number Theory (Montreal, QC, 1985), CMS Conf. Proc. \textbf{7}, Amer. Math. Soc., Providence, RI, 1987, 243--272. Author's text: \url{https://mrammurty-ai.github.io/variant.pdf}.

\bibitem{Thorner} J.~Thorner, \emph{Bounded gaps between primes in Chebotarev sets}, Res. Math. Sci. \textbf{1} (2014), art.~4. \doi{10.1186/2197-9847-1-4}.

\bibitem{Dickson} L.~E. Dickson, \emph{A new extension of Dirichlet's theorem on prime numbers}, Messenger Math. \textbf{33} (1904), 155--161. \href{https://oeis.org/wiki/File:A_new_extension_of_Dirichlet%27s_theorem_on_prime_numbers.pdf}{Scanned text}.

\bibitem{HL} G.~H. Hardy and J.~E. Littlewood, \emph{Some problems of `Partitio numerorum'; III: On the expression of a number as a sum of primes}, Acta Math. \textbf{44} (1923), 1--70. \doi{10.1007/BF02403921}.

\bibitem{BH} P.~T. Bateman and R.~A. Horn, \emph{A heuristic asymptotic formula concerning the distribution of prime numbers}, Math. Comp. \textbf{16} (1962), 363--367. \doi{10.1090/S0025-5718-1962-0148632-7}.

\bibitem{DS} K.~Destagnol and E.~Sofos, \emph{Rational points and prime values of polynomials in moderately many variables}, Bull. Sci. Math. \textbf{156} (2019), art.~102794. \doi{10.1016/j.bulsci.2019.102794}.

\bibitem{Monier} L.~Monier, \emph{Evaluation and comparison of two efficient probabilistic primality testing algorithms}, Theoret. Comput. Sci. \textbf{12} (1980), 97--108. \doi{10.1016/0304-3975(80)90007-9}.

\bibitem{IMR} S.~T. Ishmukhametov, B.~G. Mubarakov and R.~G. Rubtsova, \emph{On the number of witnesses in the Miller--Rabin primality test}, Symmetry \textbf{12} (2020), art.~890. \doi{10.3390/sym12060890}.

\bibitem{EP} P.~Erdős and C.~Pomerance, \emph{On the number of false witnesses for a composite number}, Math. Comp. \textbf{46} (1986), 259--279. \doi{10.1090/S0025-5718-1986-0815848-X}.

\bibitem{Maynard} J.~Maynard, \emph{Dense clusters of primes in subsets}, Compos. Math. \textbf{152} (2016), 1517--1554. \doi{10.1112/S0010437X16007296}.

\bibitem{Shev12} V.~Shevelev, G.~García-Pulgarín, J.~M. Velásquez-Soto and J.~H. Castillo, \emph{Overpseudoprimes, and Mersenne and Fermat numbers as primover numbers}, J. Integer Seq. \textbf{15} (2012), art.~12.7.7. \url{https://cs.uwaterloo.ca/journals/JIS/VOL15/Castillo/castillo2.html}.

\bibitem{Shev08} V.~Shevelev, \emph{An upper estimate for the overpseudoprime counting function}, preprint, 2008. \doi{10.48550/arXiv.0807.1975}.

\bibitem{Erdos76} P.~Erdős, \emph{Bemerkungen zu einer Aufgabe in den Elementen}, Arch. Math. (Basel) \textbf{27} (1976), 159--163. \doi{10.1007/BF01224655}.

\bibitem{MV} H.~L. Montgomery and R.~C. Vaughan, \emph{The large sieve}, Mathematika \textbf{20} (1973), 119--134. \doi{10.1112/S0025579300004708}.

\bibitem{PinedaRI} S.~Pineda Montoya, \emph{Residual indices of fixed multiplicative subgroups over number fields: power savings and moment obstructions}, preprint, Hexagon 2611.01069v1, 2026. \url{https://hexagonmath.org/2611.01069v1}.

\bibitem{SerreLocal} J.-P. Serre, \emph{Local Fields}, Grad. Texts in Math. \textbf{67}, Springer, New York, 1979. \doi{10.1007/978-1-4757-5673-9}.

\bibitem{GMo} L.~Grenié and G.~Molteni, \emph{Explicit versions of the prime ideal theorem for Dedekind zeta functions under GRH}, Math. Comp. \textbf{85} (2016), 889--906. \doi{10.1090/mcom3031}.

\bibitem{FI} J.~Friedlander and H.~Iwaniec, \emph{Opera de Cribro}, Amer. Math. Soc. Colloq. Publ. \textbf{57}, Amer. Math. Soc., Providence, RI, 2010. \doi{10.1090/coll/057}.

\bibitem{MurtySrinivasan} M.~R. Murty and S.~Srinivasan, \emph{Some remarks on Artin's conjecture}, Canad. Math. Bull. \textbf{30} (1987), 80--85. \doi{10.4153/CMB-1987-012-5}.

\bibitem{FM} A.~T. Felix and M.~R. Murty, \emph{On a conjecture of Erdős}, Mathematika \textbf{58} (2012), 275--289. \doi{10.1112/S0025579311008205}.

\bibitem{Pom82} C.~Pomerance, \emph{A new lower bound for the pseudoprime counting function}, Illinois J. Math. \textbf{26} (1982), 4--9. \doi{10.1215/ijm/1256046896}.

\bibitem{Sun} Z.-W. Sun, \emph{The tangent function and power residues modulo primes}, Czechoslovak Math. J. \textbf{73} (2023), 971--978. \doi{10.21136/CMJ.2023.0395-22}.

\bibitem{PCubic} P.~Pollack, \emph{Not Always Buried Deep: A Second Course in Elementary Number Theory}, Amer. Math. Soc., Providence, RI, 2009. Author's text: \url{https://www.pollack-math.net/NABDofficial.pdf}.

\bibitem{PW} L.~Pehlivan and K.~S. Williams, \emph{Some new evaluations of the Legendre symbol $\bigl(\frac{a+b\sqrt q}{p}\bigr)$}, Acta Arith. \textbf{170} (2015), 361--380. \doi{10.4064/aa170-4-4}.

\bibitem{SerreNX} J.-P. Serre, \emph{Lectures on $N_X(p)$}, Chapman \& Hall/CRC Res. Notes Math. \textbf{11}, CRC Press, Boca Raton, FL, 2012. \doi{10.1201/b11315}.

\bibitem{NIST} National Institute of Standards and Technology, \emph{Digital Signature Standard (DSS)}, FIPS PUB 186-5, 2023. \doi{10.6028/NIST.FIPS.186-5}.

\bibitem{Zenodo} S.~Pineda Montoya, J.~H. Rúa Muñoz, J.~Porras Barrera and E.~Gallego González, \emph{Replication materials for ``Simultaneous Fermat pseudoprimes with two prime factors: rank-dependent logarithmic savings and sifted base-3 proportions''}, version 1.0.0, Zenodo, 2026. \doi{10.5281/zenodo.24574022}.

\end{thebibliography}
\end{document}